%% file: main.tex
\documentclass[11pt]{article}

 \usepackage{lscape}
\usepackage{pdflscape}
\usepackage{url}

 \usepackage{blkarray,stackengine,tikz,trimclip}
\usetikzlibrary{arrows.meta}
\usetikzlibrary{arrows.meta}

\usetikzlibrary{shapes.geometric, arrows}

\tikzstyle{startstop} = [rectangle, rounded corners, 
minimum width=4cm, 
minimum height=1cm,
text centered, 
draw=black, 
fill=red!30]

\tikzstyle{io} = [trapezium, 
trapezium stretches=true, 
trapezium left angle=70, 
trapezium right angle=110, 
minimum width=4cm, 
minimum height=1cm, text centered, 
draw=black, fill=blue!30]

\tikzstyle{process} = [rectangle, 
minimum width=5cm, 
minimum height=1.0cm, 
text centered, 
text width=4.5cm, 
draw=black, 
fill=orange!30]

\tikzstyle{decision} = [diamond, 
minimum width=5cm, 
minimum height=1.0cm, 
text centered, 
text width=2.5cm, 
draw=black, 
fill=green!30]
\tikzstyle{arrow} = [thick,->,>=stealth]

\usepackage{amsthm}

\usepackage{etex}
\usepackage{amsmath,latexsym,amssymb}
\usepackage{amsfonts}
\usepackage{bbm}
\usepackage{color}
\usepackage{fancyhdr}
\usepackage{mathdots}
\usepackage{soul}

\newcommand*{\llbrace}{\{\mskip-5mu\{}
\newcommand*{\rrbrace}{\}\mskip-5mu\}}
\usepackage{ulem}

\definecolor{DarkRed}{rgb}{0.65,0.00,0.05} 
\definecolor{DarkBlue}{rgb}{0,0.08,0.45}   

\usepackage{comment}
\usepackage{stmaryrd}
\usepackage[T1]{fontenc}
\DeclareUrlCommand\email{\urlstyle{rm}}
\def\medskip{\vspace{0.125in}\noindent}

\small{
\markboth{}{}
}

\newcommand{\norm}[1]{\left\lVert#1\right\rVert}

\usepackage{mathtools}
\DeclarePairedDelimiter{\ceil}{\lceil}{\rceil}

\title{High-Order Entropy Correction with SIAC Filters\footnote{Distribution Statement A: Approved for Public Release; Distribution is Unlimited. AFRL-2023-6002.}}

\author{Matthew J. Picklo{\footnotemark[2] } \and
Ayaboe K. Edoh{ \footnotemark[3]} 
}

\footnotetext[2]{Department of Applied Mathematics and Statistics, Colorado School of Mines, Golden, CO 80401, USA. \email{mpicklo@mines.edu}}

\footnotetext[3]{Jacobs Engineering Group, Inc. – Edwards Air Force Base, California 93524, USA \email{Ayaboe.Edoh@jacobs.com}}

\date{}

\begin{document}
\maketitle
\begin{abstract}
This article considers the application of Smoothness-Increasing Accuracy-Conserving (SIAC) filtering for the non-linear stabilization of discontinuous Galerkin (DG) discretizations via entropy correction. Upon constructing discrete filters from continuous convolution SIAC kernels, the schemes are made to be conservative and are then appended to the DG method in a semi-discrete fashion. Performance of these tunable SIAC filters is compared to the local averaging typically employed in the entropy correction of finite element methods, and their capabilities are demonstrated for energy conservation as well as a shock regularization strategy based on an artificial viscosity estimate. Relaxation Runge-Kutta time integration methods are further employed in order to ensure a fully-discrete energy preserving procedure, with impacts of the overall solution accuracy being investigated for calculations of the one- and two-dimensional Burgers equation.

\end{abstract}


\input{sections/intro}

  \input{sections/DGbackground}
\input{sections/entropy_correction_for_DG}

\input{sections/SIACbackground}

\input{sections/results}

\input{sections/conclusions}
\input{sections/acknowledgments}
\bibliographystyle{siam}
\bibliography{main.bib}

\end{document}

%% file: sections/intro.tex
\section{Introduction}\label{sec:intro}

In this paper we consider scalar hyperbolic equations of the form
\begin{equation}
     \frac{d}{dt}u+\nabla \cdot \mathbf{f}(u)=0, \label{eq:cons_law}
\end{equation}
subject to appropriate initial and boundary conditions. It is well known that even for smooth initial data, classical solutions of equation \eqref{eq:cons_law} may not exist for all time, as discontinuities can develop in finite time when characteristics cross \cite{LeVeque1993}. As such, typically we enlarge the set of admissible solutions to include those that satisfy equation \eqref{eq:cons_law} in a variational, or "weak" sense. Unfortunately, these weak solutions are not unique, and furthermore when attempting to model natural phenomena, one expects a single physically-relevant solution. A selection criterion for finding this physically-relevant weak solution is to look for weak solutions satisfying entropy conditions. Let $U(u)$ be a convex function of $u$, and let $\mathbf{F}=[F_1,F_2,\hdots,F_d]$ be corresponding entropy fluxes satisfying $U'(u) f_j'(u)=F'_j(u),$ $j=1,\hdots,d$ where $d$ is the number of spatial dimensions. We say that $u$ satisfies an entropy condition for a given entropy $U$ provided
\begin{equation}
    \frac{d}{dt}U(u)+\nabla \cdot \mathbf{F} \leq 0 \label{eq:econs_law}
\end{equation}
is satisfied in the sense of distributions (see \cite{LeVeque1993} for details). The entropy solution to the original conservation law is then defined to be the weak solution that satisfies equation \eqref{eq:econs_law} for every convex entropy. Note that equality holds in equation \eqref{eq:econs_law} when the solution $u$ is still a classical solution to equation \eqref{eq:cons_law}. 

When constructing numerical methods to solve hyperbolic conservation laws, it is desirable that these methods satisfy a discrete equivalent of equation \eqref{eq:econs_law}. Such a relation can often be used to demonstrate the energy stability of a numerical discretization, for example see \cite{Gassner_Skew}. Unfortunately it is impossible for high order numerical methods to satisfy such a relation for every convex entropy, and such a restriction limits schemes to first order \cite{Lax_Harten}. As such, a single convex $U$ is generally considered.
\newpage 

In \cite{Jiang1994}, the authors demonstrated that the analytic Discontinuous Galerkin (DG) method, i.e. the DG method under exact integration, is entropy stable for the square entropy $U=u^2/2$. Unfortunately, exact integration of flux terms is in general either not possible or prohibitively expensive computationally. This has motivated alternative means of ensuring a discrete entropy (in)equality by the use of entropy conservative (dissipative) fluxes \cite{Tadmor2003} in conjunction with DG methods exhibiting the summation-by-parts (SBP) property \cite{Gassner_Skew}. The SBP characteristic enables discrete proofs for the square entropy despite inexact integration, with extensions to high-order made possible by flux-differencing implementations \cite{LeFloch:2002, Fisher:2013a}.

Another approach for recovering discrete entropy  consistency is the entropy correction approach of Abgrall and collaborators. Their correction procedure is useful for enforcing auxiliary dynamics to a baseline discretization, and has been utilized in compressible flows with respect to various secondary quantities such as angular momentum, kinetic energy, and entropy \cite{ABGRALL_Ang_mom_corr,Abgrall_Thermo_Corr}. This enforcement of additional dynamical/physical consistencies within the discretization has been shown to improve solution quality and robustness. 
The current work contributes to the correction formalism of Abgrall, which may serve as a more accessible avenue to energy stability compared to the Tadmor-type energy conserving/stable fluxes.

The framework provided by Abgrall detailed above allows for great flexibility in the choice of operators (filtering schemes) that are used for the correction. The key requirements are that the chosen operators must be conservative (assuming the underlying system is conservative) and that one possesses a measure of the volumetric growth/decay of the secondary quantities used to enforce their target dynamics.
The specific choice of filter dissipation operator that is employed, however, affects the scales at which the correction alters the solution.
One may also improve local properties of the secondary dynamics that may not be explicitly addressed by the formulation \cite{Edoh_Corr}. The element-wise averaging associated with the original correction procedure relegates the forcing to a length scale associated with the element size; meanwhile, general diffusion operators  \cite{Mattsson:2004, RANOCHA_SBP_filt, Hicken:2020} provide additional focus to subcell resolutions.

 The current work extends the use of SIAC filters for the correction procedure. These filter operators are based on a continuous convolution principle of filtering and are thus mesh agnostic and very flexible. Additionally, SIAC filters are tunable in their damping properties. When used within the finite element framework, they then provide the ability to prescribe smoothness across a range of scales that span between multiple elements and down to the subcell level. 
 Namely, the performance of SIAC-based corrections is studied relative to enforcing energy consistency for a ``vanilla" Burgers DG space discretization. To make SIAC filters compatible with the procedure of Abgrall, we first introduce a volume-weighted correction that enforces conservation properties for the operator, while maintaining its nominal damping properties. Comparison with the local-element averaging procedure of Abgrall shows that the new approach performs just as well in the target auxiliary energy discretization; this is further confirmed by a proposed LS procedure that assigns a convex weighting between the correction operators. Further, the procedure is efficiently extended to 2D via efficient line filtering \cite{LSIAC}, and fully discrete conservative properties are achieved by employing relaxation Runge-Kutta (rRK) procedures in time. Meanwhile, regularization of sharp gradients is achieved via the same correction approach, wherein the target energy decay due to shock mollification is approximated according to an artificial viscosity premise. While the inviscid Burgers equation is the PDE model considered in this work, it is expected that SIAC entropy correction will prove useful in more complex scenarios such as under-resolved high Reynolds number compressible flows.

The paper is organized as follows.
In Section 2 the baseline DG method used in this study is described. Next in Section 3 we describe a) a novel adaptation of Abgrall's original entropy correction scheme with multi-element filters, b) introduce a least-squares procedure for minimizing subcell entropy discrepancies amongst different filter-based corrections, c) describe a method for shock regularization that does not require computation of an explicit diffusive term, and d) construct a fully discrete entropy conservative scheme by means of Relaxation Runge-Kutta methods. Meanwhile Section 4 provides an overview to SIAC filters and their various parameterizations and associated properties with details on the construction of the discrete filtering operators. In Section 4 we present numerical results, detailing the subcell entropy conservation and solution quality under regularization of different filter parameterizations and combinations of filters via the least-squares procedure. We also extend these results to two dimensions, and confirm that the higher order accuracy of the underlying discretization remains unaffected by the correction.
Lastly, Section 5 provides conclusions of the current work as well as further lines of investigation.


%
%
%
%

%% file: sections/DGbackground.tex
\section{Baseline Numerical Method}\label{sec:DGSEM}

In this work, we focus on the discontinuous Galerkin (DG) method. DG methods are a type of finite element method where continuity is not enforced at element boundaries. In that sense they can also be thought of as a high order extension of finite volume techniques. These family of methods are well known for their computational efficiency and as high order methods applicable in complex geometries, see \cite{Hes08W} for an in-depth introduction to the method.

Specifically, the DG spectral element methods (DGSEM) using Legendre-Gauss-Lobatto (LGL) nodes is considered, which enable a summation-by-parts property to be satisfied by the spatial discretization. This has been shown by Gassner \cite{Gassner_Skew} to enable energy conservation, under suitable numerical flux choices, despite under-integration of volumetric terms. The scheme is formulated in one dimension as follows. Consider equation \eqref{eq:cons_law} with supplied initial data $u(x,0)=u_0(x)$ on some domain $\Omega$.
Subdivide the domain into elements $\Omega_e=[x_{e-1/2},x_{e+1/2}]$ s.t. $\cup_{e=1}^{N_{elm}}\Omega_e=\Omega$ and introduce the global to reference element mappings $\xi_e(x)=\frac{2}{\Delta x_e}(x-x_e)$ so that $\xi_e(\Omega_e)=[-1,1]$. 
Multiplication of equation \eqref{eq:cons_law} by a test function $v$ and integration over element $\Omega_e$ yields
\begin{equation*}
    \int_{\Omega_e}v\frac{d}{dt}u +v\frac{d}{dx}f(u)\;d\Omega_e=0 \ ,
\end{equation*}
next integration by parts on the flux term yields
\begin{equation}
    \int_{\Omega_e}v\frac{d}{dt}u \;d\Omega_e-\int_{\Omega_e}f(u)\frac{d}{dx}v\;d\Omega_e+[vf(u)]_{\Gamma_e}=0.
\end{equation}
On each element the solution $u$ is approximated by a polynomial of degree $p$ given by
\begin{equation*}
u_h\Big|_{\Omega_e}(x)=\sum_{k=0}^p u^e_k\ell_k(\xi_e(x)),
\end{equation*}
where $\{\ell_k\}_{k=0}^p$ are Lagrange polynomials given by
\begin{equation*}
    \ell_k(\xi)=\prod_{m=0,m\neq k}^p\frac{\xi-\xi_m}{\xi_k-\xi_m}
\end{equation*}
defined at the $p+1$ point Legendre-Gauss-Lobatto nodes $-1=\xi_0 <\hdots <\xi_{p}=1$. The values $u^e_k$ are the values of the approximation at $x^e_k=\frac{\Delta x_e}{2}\xi_k+x_e$ since $\ell_i(\xi_j)=\delta_{ij}$. 
Continuity of the approximation at element interfaces is not enforced in this method, so for the boundary fluxes one introduces a numerical flux term which solves the Riemann problem at the element boundary \cite{leveque_2002}. Choosing the test function $v=\ell_j(\xi_e)$ yields the form
\begin{equation}
   \frac{d}{dt}\int_{\Omega_e}u_h\ell_j(\xi_e(x)) \;d\Omega_e-\int_{\Omega_e}f(u_h)\frac{d}{dx}\ell_j(\xi_e(x))\;d\Omega_e+\ell_j(1) f^{num}_{e+1/2}-\ell_j(-1) f^{num}_{e-1/2}=0,\;\;\; j=0,\hdots,p. \label{eq:pre_quad}
\end{equation}

Associated with the LGL nodes is the LGL quadrature rule. Denoting the weights and nodes by $\{\omega_k,\xi_k\}_{k=0}^p$ we have
\begin{equation}  
\int_{-1}^1h(\xi)\;d\xi\approx\sum_{k=0}^p\omega_kh(\xi_k),
\end{equation}
which is exact for polynomials of degree $2p-1$. 
Approximating the integral volume integrals in equation $\eqref{eq:pre_quad}$ using the $p+1$ node LGL quadrature results in under-integration or ``Mass-lumping" of the first integral, though it produces a diagonal mass matrix. Upon substitution of the element-specific expression for $u_h$ we obtain the DGSEM scheme
\begin{equation}
    \frac{\Delta x_e}{2}\omega_j\frac{d}{dt}u^e_j-\sum_{k=0}^p\omega_k f^e_k \frac{d\ell_j}{d\xi}(\xi_k)+\ell_j(1)f^{num}_{e+1/2}-\ell_j(-1)f^{num}_{e-1/2},\;j=0,\hdots,p, \;\;\;e=1,\hdots,N_{elm}.
\end{equation}
A matrix-vector interpretation is allowed by the scheme given above. Letting $\mathbf{u}^e=[u^e_0,\hdots,u^e_p]^T$, $\mathbf{f}^e=[f^e_0,\hdots,f^e_p]^T$, and $\mathbf{f}_{num}^e=[f^{num}_{e-1/2},0,\hdots,0,f^{num}_{e+1/2}]^T$, we obtain the expression 
\begin{equation}
    \frac{\Delta x_e}{2}\mathbf{M}\left(\frac{d}{dt}\mathbf{u}^e\right)=\mathbf{D}^T\mathbf{M}\mathbf{f}^e-\mathbf{B}\mathbf{f}^e_{num}\label{eq:ref_DG_local},
\end{equation}
where $\mathbf{M}=\text{diag}([\omega_0,\hdots,\omega_p])$, $(\mathbf{D})_{ij}=\frac{d\ell_j}{d\xi}(\xi_i)$, and $\mathbf{B}=\text{diag}([-1,0,\hdots,0,1])$.
As will become important in what follows, the local scheme can be expressed in a global format. To match conventions applied later, we choose here to include the mapping factor $\Delta x_e/2$ in the quadrature matrix i.e. $\mathbf{M}_e=\frac{\Delta x_e}{2}\mathbf{M}$. The subscript $G$ used here will denote the global appending of vectors
\begin{equation*}
    \mathbf{u}_G=[\mathbf{u}^T_1,\hdots,\mathbf{u}^T_{N_{elm}}]^T,
\end{equation*}
or the global block concatenation of matrices 
\begin{equation*}
    \mathbf{M}_G=\begin{bmatrix}
    \mathbf{M}_1&&\\
    & \ddots&\\
    & & \mathbf{M}_{N_{elm}}
    \end{bmatrix}.
\end{equation*}
Using this convention, the global scheme is expressible as
\begin{equation}
    \mathbf{M}_G\left(\frac{d}{dt}\mathbf{u}_G\right)=\mathbf{D}^T_G(\mathbf{I}_{N_{elm}}\otimes \mathbf{M})\mathbf{f}_G-\mathbf{B}_G\mathbf{f}_{G,num},
\end{equation}
or directly as the semi-discrete ODE
\begin{equation}
\begin{aligned}
    \frac{d}{dt}\mathbf{u}_G&=\mathbf{M}^{-1}_G(\mathbf{D}^T_G(\mathbf{I}_{N_{elm}}\otimes \mathbf{M})\mathbf{f}_G-\mathbf{B}_G\mathbf{f}_{G,num}),\\
    &=\mathbf{r}_G. \label{eq:ref_DGSEM}
\end{aligned}
\end{equation}
When the scale is apparent from the context, the $G$ subscript is omitted.

%% file: sections/entropy_correction_for_DG.tex
\section{Entropy Correction}\label{sec:EC}


Consider now the one-dimensional inviscid Burgers equation for which $f(u)=u^2/2$. Integration of the original PDE over some domain $\Omega$ coupled with integration by parts allows an expression for the volumetric time rate of change of $u$ in the domain:
\begin{equation}
    \frac{d}{dt}\int_{\Omega} u\;d\Omega=-f|_{d\Omega}.
\end{equation}
Similarly, considering as an entropy the square entropy $U=u^2/2$, we can multiply both sides by the entropy variable $w=\frac{d}{du}U=u$ and perform the same procedure to obtain
\begin{equation}
    \frac{d}{dt}\int_{\Omega} U\;d\Omega=-F|_{d\Omega}
\end{equation}
where $F(u)=u^3/3$. In the continuous case, one can determine the volumetric time rate of change of the conserved quantity $u$ and its energy by the fluxes at the boundary of the domain. In the following sections we are concerned with the construction of correction terms to our numerical schemes with periodic boundary conditions that will enable conservation 
\begin{align}
    \frac{d}{dt}\int_{\Omega} u\;d\Omega&=0,
    \intertext{and entropy conservation/stability}
    \frac{d}{dt}\int_{\Omega} U\;d\Omega&\leq0 \label{eq:Ent_diss}
\end{align}
of the resulting discretization in a discrete sense. When considering the numerical methods described previously equation \eqref{eq:ref_DGSEM}, it can be shown  (see similar derivation in \cite{offner2023approximation}) that on a single element
\begin{align}
    \mathbf{1}^T\mathbf{M}_e\left(\frac{d}{dt}\mathbf{u}_e\right)&=f^{num}_{e-1/2}-f^{num}_{e+1/2}.\\
    \intertext{Hence, summing over all the elements yields} 
    \mathbf{1}^T\mathbf{M}_G\left(\frac{d}{dt}\mathbf{u}_G\right)&=f^{num}_{-1/2}-f^{num}_{N+1/2},
    \intertext{that is to say that conservation is determined entirely by inflow/outflow at the boundary of the domain. For periodic boundary conditions, one recovers strict conservation}
\mathbf{1}^T\mathbf{M}_G\left(\frac{d}{dt}\mathbf{u}_G\right)&=0.
\end{align} 

Stability in energy requires additional assumption on the scheme at hand. In \cite{Gassner_Skew}, the author considered skew formulations of the Burgers DGSEM and demonstrated how to choose numerical fluxes to ensure energy conservation/stability. Here, we instead consider the general approach of Abgrall \cite{Abgrall2018} and apply a correction term to ensure energy stability. Dropping the global indicators here for readability, we have 
\begin{align}
    \frac{1}{2}\mathbf{u}^T\mathbf{M}\left(\frac{d}{dt}\mathbf{u}\right)&=\frac{1}{2}\mathbf{u}^T\mathbf{M}(\mathbf{r}+\mathbf{c})\\
&=\frac{1}{2}\mathbf{u}^T\mathbf{M}\mathbf{r}+\frac{1}{2}\mathbf{u}^T\mathbf{M}\mathbf{c}.
\end{align}
As will be detailed later, $\mathbf{c}$ can be chosen to cancel out the first term resulting in energy conservation. Similarly, it can be augmented with shock regularization capabilities to ensure energy dissipation.

More generally, the correction can be applied to ensure a discretization can satisfy spatial entropy conservation condition. The aforementioned entropy correction formulation is as follows. Given the element-level semi-discrete ODE of equation \eqref{eq:ref_DG_local} derived above (where we drop the $e$ subscripts for readability), append a correction term to ensure conservation of discrete entropy. The new equation is given by
\begin{equation}
    \frac{d}{dt} \mathbf{u}=\mathbf{r}+\mathbf{c},
\end{equation}
where $\mathbf{c}$ is a correction to the degrees of freedom of the local element of the form
\begin{equation}
\mathbf{c}=\alpha\Big(\mathbf{w}-\frac{\mathbf{1}^T\mathbf{M}\mathbf{w}}{\mathbf{1}^T\mathbf{M}\mathbf{1}}\mathbf{1}\Big).\label{eq:corr_ab}
\end{equation}
Here $\mathbf{w}$ is an entropy variable for the given entropy, and
\begin{align*}
\alpha=\frac{\phi_{loc}}{\mathbf{w}^T\mathbf{M}\mathbf{w}-\frac{(\mathbf{1}^T\mathbf{M}\mathbf{w})^2}{\mathbf{1}^T\mathbf{M}\mathbf{1}}},
\end{align*}
where $\phi_{loc}$ is selected to cancel out the entropy flux on the local element. This term is defined as 
\begin{equation}
\phi_{loc}=-\int_{\partial \Omega_e}\mathbf{n}\cdot F^{num}\;d\Omega_e-\mathbf{w}^T\mathbf{M}\mathbf{r},
\end{equation}
where $F_{num}$ is the numerical entropy flux associated with the predetermined numerical flux $f_{num}$ by the relation
\begin{equation}
   F^{num}_{e+1/2}=\llbrace w_{e+1/2}\rrbrace f^{num}_{e+1/2}-\llbrace \psi_{e+1/2}\rrbrace,
\end{equation}
where $\psi_{e+1/2}=w_{e+1/2}f_{e+1/2}-F_{e+1/2}$. The correction as written satisfies conservation 
\begin{equation}
\mathbf{1}^T\mathbf{M}\mathbf{c}=0 \label{eq:cons}
\end{equation}
and the desired correction property
\begin{equation}
\mathbf{w}^T\mathbf{M}\mathbf{c}=-\int_{\partial \Omega_e}\mathbf{n}\cdot F^{num}\;d\Omega_e-\mathbf{w}^T\mathbf{M}\mathbf{r}.\label{eq:econs}
\end{equation}
Note that typically a small constant is added to denominator of $\alpha$ to prevent division by zero in the case of a constant solution in an element.
This form of correction was shown by the authors of \cite{Abgrall2022} to be the solution to the minimization problem $\min \norm{\mathbf{c}}_{\mathbf{M}}$ subject to local the constraints of equations \eqref{eq:cons} and \eqref{eq:econs}.

Having discussed this previous approach and in order to motivate the choice of alternative corrections, note that equation \eqref{eq:corr_ab} uses an average of the local data that may be Defined by the filter $\mathbf{K}_{loc. avg}=\frac{\mathbf{1}\mathbf{1}^T\mathbf{M}}{\mathbf{1}^T\mathbf{M}\mathbf{1}}$. We can therefore rewrite equation \eqref{eq:corr_ab} as
\begin{equation}
    \mathbf{c}=(\phi_{loc})\frac{(\mathbf{I}-\mathbf{K}_{loc. avg})\mathbf{w}}{\mathbf{w}^T\mathbf{M}(\mathbf{I}-\mathbf{K}_{loc. avg})\mathbf{w}}.
\end{equation}
That is to say, this correction employs local element-wise averaging which can be interpreted as a specific type of filter. When comparing to alternative formulations, we will refer to this element local correction approach of Abgrall as the ``local" approach or ``local" method.

The present work expands the correction approach to include multi-element filters such as the SIAC operators previously mentioned. In this section the $G$ subscript is once again used to denote a matrix or vector that is a global concatenation of single element components. Supposing a global target dynamic, define 
\begin{equation}
    \phi_{glob.}=\int_{\Omega}\frac{d}{dt} U_{targ.}\;d\Omega-\mathbf{w}^T\mathbf{M}_{G}\mathbf{r}_G.
\end{equation}
Then one can construct a global correction term as 
\begin{align}
    \mathbf{c}_G&=\alpha_{glob.}(\mathbf{I}_G-\mathbf{K}_{G})\mathbf{w}_G,
    \intertext{where}
{\alpha}_{glob.}&=\frac{{\phi}_{glob.}}{\mathbf{w}_G^T\mathbf{M}_G(\mathbf{I}_G-\mathbf{K}_{G})\mathbf{w}_G},\label{eq:global_corr}
\end{align}
and $\mathbf{K}_{G}$ is an unspecified conservative filter ($\mathbf{1}^T\mathbf{M}_G\mathbf{K}_G\mathbf{v}=\mathbf{1}^T\mathbf{M}_G\mathbf{v}$ for arbitrary $\mathbf{v}$.)
The augmented global scheme $\frac{d}{dt}\mathbf{u}_G=\mathbf{r}_G+\mathbf{c}_G$ will remain conservative and will satisfy the target global dynamic: $\mathbf{1}_G^T\mathbf{M}_G\mathbf{c}_G=0$ and $\mathbf{w}_G^T\mathbf{M}_G(\mathbf{c}_G+\mathbf{r}_G)=-\int_{\Omega}\frac{d}{dt} U_{targ.}\;d\Omega$, respectively. Note that while the new correction does not to the authors knowledge satisfy an optimally condition like that of the local correction, the numerical results to follow show that the global correction can better reduce the relative magnitude of the correction: $\norm{\mathbf{c}}_{\mathbf{M}_G}/\norm{\mathbf{r}}_{\mathbf{M}_G}$. This is not a contradiction as allowing for multi-element filters affords more flexibility in allocation of the correction.

\subsection{An Adaptive Filter Selection Strategy}\label{subsec:LS}
To further improve the global correction procedure, we consider a constrained least-squares (LS) problem for selecting an optimal convex combination of two correction terms. Equivalently, this provides a means of differentiating between filtering operators and their implied scales of influence. Since every correction considered maintains global conservation, a subcell metric of entropy violation is used to differentiate between them.
Such notions of local entropy consistency suggest greater solution fidelity. The metric used here comes directly from equation \eqref{eq:econs_law} and is given by
\begin{align}
D_h(x^e_i)&=w^e_i\frac{d}{dt}u^e_i+\frac{d}{dx}F(u_h(x^e_i))\label{eq:EC_residual}\\
&=w^e_i\mathbf{r}(x^e_i)-\mathbf{rF}(x^e_i)\label{eq:entropy_metic},
\end{align}
where $x^e_i$ is the $i$th node of element $e$, $\mathbf{r}$ is the DGSEM discretization of the RHS of the conservation law previously described, and $\mathbf{rF}$ is the discretization via DGSEM of the entropy flux term $F(u)_x$. The interpretation here is that if the local entropy is being well preserved at a given node, $D_h$ should be near 0. We note that while this was the choice of measure for subcell entropy conservation employed in this work, other options are also valid.

With a discrepancy metric in hand, the constrained LS problem is as follows: Given two entropy correction terms $\mathbf{c}_1$ and $\mathbf{c}_2$ resulting from conservative filters $\mathbf{K}_1$ and $\mathbf{K}_2$, choose a correction 
\begin{align}    
\mathbf{c}_{\theta}&=\theta \cdot \mathbf{c}_1+(1-\theta) \cdot\mathbf{c}_2\\
&=\mathbf{c}_2+\theta \cdot(\mathbf{c}_1-\mathbf{c}_2)
\end{align}
parameterized by $\theta\in[0,1]$ s.t. $\norm{\mathbf{D}_h}_{\mathbf{M}_G}$ is minimized. The form of the vector $\mathbf{D}_h$ in terms of $\theta$ is found by substitution of the discrete scheme into equation $\eqref{eq:EC_residual}$ ($\odot$ denoting elementwise multiplication):
\begin{align}
\mathbf{D}_h&=\mathbf{w}\odot\frac{d}{dt}\mathbf{u}+\frac{d}{dx}\mathbf{F(u)} \nonumber\\
    &=\mathbf{w}\odot(\mathbf{r}+\mathbf{c}_{\theta})-\mathbf{rF}\nonumber\\
    &=\mathbf{w}\odot(\mathbf{r}+\mathbf{c}_2)-\mathbf{rF}+\theta (\mathbf{w})\odot(\mathbf{c}_1-\mathbf{c}_2)\nonumber\\
    &=\mathbf{a}+\theta\mathbf{b} \label{eq:LS}
\end{align}
where $\mathbf{a}=\mathbf{w}\odot(\mathbf{r}+\mathbf{c}_2)-\mathbf{rF}$ and $\mathbf{b}=(\mathbf{w})\odot(\mathbf{c}_1-\mathbf{c}_2).$ The optimization problem is then solved numerically at each step of the time-evolution of $\mathbf{u}$ to choose an appropriate correction term. Without a constraint on $\theta$, the least-squares solution is given in a straight forward manner as 
\begin{align}
\theta=-(\mathbf{b}^T\mathbf{M}_G\mathbf{b})^{-1}(\mathbf{b}^T\mathbf{M}_G\mathbf{a}).
\end{align}
Note that this is just the standard least-square solution procedure for $\norm{\mathbf{D_h}}_{2}$, but the structure of $\mathbf{M}_G$ makes it easy to find the least-squares solution in its associated norm. As the quadrature matrix is positive definite there exists a matrix $\mathbf{R}$ s.t. $\mathbf{R}^T\mathbf{R}=\mathbf{M}_G$, and so it can be shown that the 2-norm LS solution to  \eqref{eq:LS} multiplied by $\mathbf{R}$ is the same as the $\mathbf{M}_G-$norm LS solution. 

The explicit expression for $\theta$ provided above does not require any numerical optimization solvers, however it unfortunately does not enforce the convex-combination constraint $\theta\in[0,1]$ as part of the least squares problem. Here the authors employed the MATLAB \text{lsqlin} \cite{MATLABlsqlin} routine to numerically solve for $\theta$ satisfying the constrained LS minimization problem, which uses an active-set quadratic programming algorithm \cite{Active_set}, see documentation in \cite{MATLABquad} for details.


\subsection{Shock Regularization based on Artificial Viscosity}\label{subsec:reg}
Referencing the entropy inequality given in equation \eqref{eq:econs}, it is expected that entropy is conserved for smooth solutions, and decreases for non-smooth ones. Unfortunately, it is difficult to determine numerically when a discontinuity has formed in the numerical solution, and a conservative correction should be replaced with a dissipative one. This is compounded in that DGSEM already has only weak continuity at element interfaces. 
Additionally, high order numerical methods such as DGSEM often result from numerical instability in the form of Gibb's oscillations when trying to resolve discontinuities. As a result stabilization techniques are necessary to prevent simulations from crashing post shock-formation. Various methodologies have been proposed for stabilization such as filtering, special entropy dissipative fluxes, artificial dissipation, and (subcell)limiters \cite{Glaubitz2019,GUERMOND,Gassner_Skew,subcell_DG}, see \cite{hest_CLbook} chapter 12 for an overview of some common methods. The approach described here embeds artificial dissipation within the correction technique.

The idea of artificial viscosity is motivated by an equivalent definition of an entropy solution where the entropy solution is obtained as the limiting solution as $\epsilon\rightarrow 0$ of the viscous PDE
\[\frac{d}{dt}u_{\epsilon}+\nabla \cdot \mathbf{f}(u_{\epsilon})= \nabla \cdot (\epsilon\nabla u).\]
For this parabolic problem the solution remains smooth, with the addition of diffusion causing the sharp features of a discontinuity to be smeared. Adding artificial viscosity to a numerical scheme can then stabilize the underlying method and reduce spurious oscillations. Unfortunately, discretizing the high-order derivative terms in the parabolic equation often places more restrictive time-stepping constraints on explicit methods. 

To avoid these issues, the corrections previously discussed are instead augmented to produce an entropy dissipative scheme by utilizing an artificial dissipation mechanism, whereby dissipation is ensured by augmenting the magnitude of the spatial correction term. Note that we term this artificial dissipation rather than viscosity as we do not explicitly append a diffusive term to the underlying PDE.


The change appears in the definition of the magnitude term ${\phi}_{glob.}$. For entropy conservation in the global filter case, a spatially constant magnitude $\phi_{glob.}=-(\mathbf{w}^T\mathbf{M}\mathbf{r})$ was chosen to cancel out the global entropy production of the scheme. To ensure dissipation, we simply need to augment $\phi_{glob.}$ with a non-positive constant ${\phi}_{diss}$. That is replace $\phi_{glob.}$ in equation \eqref{eq:global_corr} with
\begin{equation}
{\phi}^{\epsilon}_{glob.}=\phi_{glob.}+{\phi}_{diss}.
\end{equation}
Note that incorporating dissipation in this way will not affect the global conservation of the method as only the scaling of the correction changes, not its global conservation properties.
We discuss the choice of the constant $\phi_{diss}$ below specifically for the one-dimensional Burgers equation, though we note that the procedure extends to more general conservation laws and has been applied by Edoh \cite{Edoh_Corr} with respect to regularizing the Euler system in the presence of shocks. 

Consider the time change of energy for the viscous Burgers equation 
\[\frac{d}{dt}u+\frac{d}{dx}(u^2/2)= \frac{d}{dx}\left(\nu_{AV}\frac{d}{dx}\right)u.\] 
Following an approach similar to that producing \eqref{eq:Ent_diss}, we obtain using integration by parts
\begin{align}
    \frac{d}{dt}\int_{\Omega} U\;d\Omega&= -F\Big|_{d\Omega}+\left(\nu_{AV}\frac{d}{dx}U\right)\Big|_{d\Omega}-\int_{\Omega}\nu_{AV}\left(\frac{d}{dx}u\right)^2\;d\Omega.
    \intertext{Supposing $\nu_{AV}\geq 0$ makes the third term non-negative and so}
    \frac{d}{dt}\int_{\Omega} U\;d\Omega&\leq -F\Big|_{d\Omega}+\left(\nu_{AV}\frac{d}{dx}U\right)\Big|_{d\Omega}.
\end{align}
In order to enforce entropy dissipation, we use the third term integrated over a single element as an indicator of the magnitude of the dissipation,
\begin{equation}
    {\phi}^e_{diss}=\frac{2}{\Delta x_e} (\mathbf{D}\mathbf{u}_e)^T\mathbf{M}[\nu_{AV}^e]\mathbf{D}\mathbf{u}_e,
\end{equation}
where $\mathbf{D}$ is the differentiation operator in equation \eqref{eq:ref_DG_local}, and the inverse Jacobian term is resultant from switching the derivatives from the reference element $[-1,1]$ to the physical element $\Omega_e$. In what follows we pursue an approach such that the viscosity term $\nu_{AV}$ is constant within an element. Unfortunately, the global correction term approach requires a constant dissipative term be applied globally. To that end, $\phi_{diss}$ is chosen to be the sum of each individual element's dissipative magnitude: $\phi_{diss}=\sum_e \phi^e_{diss}$. Meanwhile for the local approach, a piecewise constant viscosity is acceptable as changing the magnitude of the local correction term element to element does not impact conservation.

In determining $\nu_{AV}$, a variant of entropy viscosity inspired by the approach of Guermond et al. \cite{GUERMOND} is considered. Defining for each element the minimum distance between LGL nodes on that element $h_e=\frac{\Delta x_e}{2}\min_{i\neq j}|\xi_i-\xi_j|$, then one makes use of the entropy residual $D_h(x^e_i)$ previously defined in equation \eqref{eq:entropy_metic} to define the entropy viscosity as
\begin{equation}
    \nu^e_{E}=c_E( h_e^2 \max_k|D_h(x^e_k)|)/\norm{\mathbf{U}-\Bar{\mathbf{U}}}_{\infty}, 
\end{equation}
where for Burgers equation with the square entropy, $\bar{\mathbf{U}}=\frac{1}{2}\mathbf{u}_G^T\mathbf{M}_G\mathbf{u}_G/(\mathbf{1}^T\mathbf{M}_G\mathbf{1})$. The tuning parameter $c_E$ is user supplied to specify the amount of dissipation. The maximum viscosity $\nu_{max}$ is given by 
\begin{equation}
    \nu^e_{max}=c_{max}(\Delta x^e \max_k|f'(u^e_k)|),
\end{equation}
where again $c_{max}$ is a tuning parameter. The element-wise viscosity is then selected as the minimum of the two quantities: $\nu^e_{AV}=\min\{\nu^e_E,\nu^e_{\max}\}$. The performance of this technique and the problem-dependent tuning parameters will discussed in the results section. Note that the magnitude of the dissipation for this approach is determined by the error metric $D_h$ and the volume integral of the gradient of the solution squared; however, the spatial allocation of this dissipation is determined by the kernel filter applied. 
\newline

\subsection{Fully-Discrete EC via relaxation Runge-Kutta time-stepping}\label{subsec:rRK}
Having described a semi-discrete numerical scheme which is entropy conserving in space, we now demonstrate how to create fully discrete energy conserving schemes by applying the conservative time-stepping relaxation Runge-Kutta (rRK) schemes of Ketcheson \cite{Ketchesonrrk}. Note that applying this approach for conservation of more general entropies than the square entropy leads to root finding problems whose efficient solution we do not address here (see \cite{KetchesonEuler}). 

Consider the ODE IVP
\begin{equation}
    \frac{d}{dt}u=f(u(t),t),\;\;\;u(t_0)=u_0,
\end{equation}
and suppose that the system is dissipative with respect to some inner-product:
\begin{equation}
    \frac{d}{dt}\norm{u(t)}^2=2\langle u, f(u,t)\rangle \leq0.
\end{equation}
The goal is to ensure that discrete time-stepping yields a solution that is also dissipative \begin{equation}
\norm{u^{n+1}}\leq \norm{u^n}.
\end{equation}
Similarly, for conservative systems ($\leq\; \rightarrow\; =$), we would like conservation in some discrete norm to hold. For our purposes, the inner-products and norms considered are those defined by global s.p.d. quadrature matrices. 

For time-stepping methods, Ketcheson considered Runge-Kutta methods which are given by
\begin{equation}
\begin{aligned}
    y_i&=u^n+\Delta t \sum_{j=1}^s a_{ij}f(y_j,t_n+c_j\Delta t) \\
    u^{n+1}&=u^n+\Delta t \sum_{j=1}^sb_jf(y_j,t_n+c_j\Delta t),
\end{aligned} \label{eq:RK}
\end{equation}
where $\mathbf{A}\in \mathbb{R}^{s\times s}$, and $\mathbf{b},\mathbf{c}\in \mathbb{R}^{s}$. Further, it is assumed that $c_j=\sum_{i=1}^s a_{ij}$. Following \cite{Ketchesonrrk} we introduce the notational convenience 
\begin{equation*}
    f_i=f(y_i,t_n+c_i\Delta t).
\end{equation*}
The change in energy over a single time step can be shown to satisfy
\begin{equation}
    \norm{u^{n+1}}^2-\norm{u^{n}}^2=2\Delta t\sum_{j=1}^sb_j\langle y_j,f_j\rangle-2(\Delta t)^2\sum_{i,j=1}^sb_ia_{ij}\langle f_j,f_i\rangle+(\Delta t)^2\sum_{i,j=1}^sb_ib_{j}\langle f_i,f_j\rangle. \label{eq:energy_diff_RK}
\end{equation}
It is the approach of rRK methods to introduce an adaptive time step $\gamma_n \Delta t$ such that for 
\begin{equation}
    u^{n+1}_{\gamma}=u^n+\Delta t \gamma_n\sum_{j=1}^sb_jf(y_j,t_n+c_j\Delta t)\label{eq:adap_uRK}
\end{equation}
the change of energy over that adapted time step is 0. To see how the parameter $\gamma_n$ is determined, note that the energy change between $t_n$ and $t_n+\gamma_n\Delta t$ is given by 
\begin{align*}
    \norm{u^{n+1}_{\gamma_n}}^2-\norm{u^{n}}^2=2\Delta t\gamma_n\sum_{j=1}^sb_j\langle y_j,f_j\rangle-2\Delta t^2\gamma_n\sum_{i,j=1}^sb_ia_{ij}\langle f_j,f_i\rangle+(\Delta t\gamma_n)^2\sum_{i,j=1}^sb_ib_{j}\langle f_i,f_j\rangle.
\end{align*}
In a conservative system, the first term is zero, whereas in dissipative systems it is non-positive provided $b_j>0$ for $j=1,\hdots,s$. The baseline RK schemes employed herein will satisfy this later requirement. The goal now is to choose $\gamma_n$ s.t. the two remaining terms cancel out. Algebra manipulations yield a solution 
\begin{equation}
    \gamma_n=\begin{cases}
\frac{2\sum_{i,j=1}^sb_ia_{ij}\langle f_i,f_j\rangle}{\sum_{i,j=1}^sb_ib_{j}\langle f_i,f_j\rangle},& \norm{\sum_{j=1}^sb_jf_j}^2\neq0,\\
1,&\text{Else}
    \end{cases}.
\end{equation}
This last condition is imposed to prevent division by $0$. 

In relating the relaxation RK approach to the entropy-correction approach for Burgers equation with the squared entropy, we note that given a conservative/dissipative semi-discrete system (e.g., $\frac{d}{dt}\mathbf{u}=\mathbf{r}+\mathbf{c}$) 
one can evolve the system in time via the rRK method and preserve the conservative/dissipative nature. In the following section this capability is put to use in a host of test problems.

%% file: sections/SIACbackground.tex
\section{SIAC Filter Overview}\label{sec:SIAC}

In constructing filters to be applied in the entropy correction context, we investigate the Smoothness-Increasing Accuracy-Conserving (SIAC) family of filters already prevalent in the DG community owing to their geometric flexibility and past successes enabling superconvergence and improved data visualization in discontinuous data \cite{Cockburn2002,Mirzaee2013,LSIAC,XLiOne,Steffen2008}. SIAC filters are a convolution kernel filter composed of a linear combination of shifted B-splines, where the coefficients in that combination are chosen so that certain moment, or equivalently polynomial reproduction conditions are satisfied.

Letting $u_h$ denote data of a continuous spatial variable, the filtered data is given by
\begin{align}
    u_h^{\star}(x)&=K_H\star u_h\\
    &=\int_{\mathbb{R}}K_H(x-y)u_h(y)\;dy,
\end{align}
where $K_H(x)=\frac{1}{H}K_{\mathbf{T}}(x/H)$ is a scaled SIAC kernel function. The SIAC kernel can be defined as 
\begin{equation}
    K_{\mathbf{T}}(x)=\sum_{\gamma=0}^r c_{\gamma}B_{\mathbf{T_{\gamma}}}(x),
\end{equation}
where $\mathbf{T}$ is a corresponding knot matrix whose $(\gamma+1)$th rows contain the knot sequence associated with the B-spline $B_{\mathbf{T_{\gamma}}}$ (see de Boor \cite{deBoor} for details on computing these). In general the knot matrix can vary spatially over the domain which is very important when dealing with non-periodic boundary conditions. In this work we will instead consider a simplified parameterization for periodic data which focuses on smoothness, number of moments, and scaling of the kernel in which case the unscaled kernel is given by
\begin{align*}
    K^{(r+1,\ell)}=\sum_{\gamma=0}^r c_{\gamma}B^{\ell}(x+r/2-\gamma).
\end{align*}
Here $B^{\ell}$ is the $\ell$th order central B-spline which can be computed by following recurrence relation:
\begin{align*}
    B^1(t)&=\chi_{[-1/2,1/2)}(t)\\
    B^{\ell}(t)&=(B^{\ell-1}\star B^1)(t)\\
  &=\frac{1}{\ell-1}\Big[\Big(\frac{\ell-1}{2}+t\Big)B^{\ell-1}(t+1/2)+\Big(\frac{\ell-1}{2}-t\Big)B^{\ell-1}(t-1/2)\Big].
\end{align*}
Importantly, the splines have compact support which is then inherited by the kernel itself. The number of splines $r+1$ is tied to the number of consistency + moment conditions that we wish the kernel to satisfy i.e. 
\begin{align}
        \int_{\mathbb{R}}K^{(r+1,\ell)}(x)x^k&=\delta_{0k},\;\;\;\text{for}\;\;\;k=0,\hdots,r.
\end{align}
In the context of obtaining superconvergence, typically for an approximation of degree $p$, $r$ is chosen as $2p+1$ to that the moment information of the data is retained, while $\ell=p+1$ is chosen to ensure the kernel and the resulting data is sufficiently smooth for various error estimates to apply \cite{Cockburn2002}. The scaling is then historically chosen as the maximum edge length of the mesh, though improved performance has been shown using adaptive scalings on nonuniform meshes, see \cite{Jallepalli2019B}.

In higher dimensions, extensions of the SIAC methodology have been developed and have proven effective over complicated geometries \cite{LSIAC,Mirzaee2013}. Perhaps the most straightforward extension is the tensor-product kernel:
\begin{equation}
K^{(\mathbf{r}+1,\mathbf{\ell})}_{\mathbf{H}}(\mathbf{x})=K^{(r_1+1,\ell_1)}_{H_1}({x_1})\otimes K^{(r_2+1,\ell_2)}_{H_2}({x_2}).
\end{equation}
Unfortunately, the dimension of the quadrature necessary to realize the convolution also grows with the dimension and can become computationally expensive. A more computationally efficient alternative is the Line-SIAC (LSIAC) filter  \cite{LSIAC} given by
\begin{align}
u^{\star}_h(\mathbf{x})&=K^{(r+1,\ell)}_{\Gamma,H}\star u_h\\
&=\int_{\Gamma}K^{(r+1,\ell)}_H(t)u_h(\Gamma(t))\;dt,\\
\intertext{with}
\Gamma(t)&=\mathbf{x}+t(\cos(\theta),\sin(\theta)).
\end{align}
For the line filter, the support is one-dimensional and so too are the quadratures. Typically for uniform quadrilateral meshes, the orientation $\theta$ is chosen to align with the element diagonals. In two-dimensional tests to follow, we will employ the LSIAC filter. For an example of the LSIAC procedure, consider the depiction in Figure \ref{fig:LSIAC_Example} (Left) where the $K^{(1,1)}_{\Gamma,\sqrt{2}\Delta x}$ filter is applied to a $p=5$ approximation.  

\input{figures/files/LSIAC_Example}

\subsection{Filter Matrix Construction}
For applications in later sections, we detail here the construction of SIAC filtering matrices for periodic DG data in one dimension. We exclude the details for higher dimensions as they follow the same approach but with more complicated indexing. When applying the matrix form of the SIAC filter, we are effectively truncating the filtered expression so that it remains in the same polynomial space as the original data. That is, in this case we are sampling the filtered approximation at the node locations of the original Lagrange basis.

Before building a global filtering matrix, consider the evaluation of the filtered expression at a node. We have
\begin{align*}
    u^{\star}_h(x^e_{k})&=\int_{\Omega} K_H(x^e_{k}-y)u_h(y)\;dy\\
&=\sum_{e=1}^{N_{elm}}\sum_{k=0}^pu^e_k\int_{\Omega_e}K_H(x^e_{k}-y) L^e_k(y)\;dy\\
&=(\mathbf{k}_{e,k})^T\mathbf{u},
\end{align*}
where
\begin{equation}
    (\mathbf{k}_{e,k})_{z}=\int_{\Omega_j}K_H(x^e_{k}-y) L^j_i(y)\;dy,\;\;\;z(i,j)=i+(j-1)(p+1), 
\end{equation}
for $i=1,\hdots,p+1$, $j=1,\hdots,N_{elm}$. Introducing an index $z(e,k)=k+(e-1)(p+1)$ for a global ordering of the row vectors, we assemble the global filtering matrix $\mathbf{K}=[\mathbf{k}_{1},\hdots,\mathbf{k}_{N_{elm}(p+1})]^T$. Multiplication of this matrix on the nodal vector returns the values of the filtered approximation at those nodes:
\begin{equation}
    \mathbf{u}^{\star}=\mathbf{K}\mathbf{u}.
\end{equation}
The compact support of the kernel will generally allow for a sparse filtering matrix as can be seen for the LSIAC filter in Figure \ref{fig:LSIAC_Example} (Right). As an example of the different discrete filtering operators produced by this procedure, the local element stencils for the $K^{{1,1}}_{\Delta x}$ and $K^{(3,2)}_{\Delta x}$ kernels for $p=1$ data in one dimension on a uniform grid are given by
\[\textbf{stencil}^{(1,1)}_{\Delta x}=\frac{1}{8}\begin{bmatrix}
\begin{array}{cc|cc|cc} 
1&3& \textcolor{red}{\mathbf{3}}& \mathbf{1}&0&0 \\ \hline
0&0&\mathbf{1} & \textcolor{red}{\mathbf{3}}& 3 &1
\end{array}
    \end{bmatrix},
\]
and
\[\textbf{stencil}^{(3,2)}_{\Delta x}=\frac{1}{72}\begin{bmatrix}
\begin{array}{cc|cc|cc|cc|cc}
-1&-2&12 & 27& \textcolor{red}{\mathbf{27}} &\mathbf{12}&-2&-1&0&0\\ \hline
0&0&-1&-2&\mathbf{12} & \textcolor{red}{\mathbf{27}}& {{27}} &{12}&-2&-1\\
\end{array}
    \end{bmatrix}.
\]
Here the red values denotes the location of the node being filtered by the filtering procedure, and the boldface values denote the node of the reference element. The vertical lines separate the nodal coefficients of different elements. The whole discrete filtering matrix in the periodic case is then obtained by periodic repetition of these stencils. Note that for higher polynomial degrees, while the stencils for boundary nodes are equal up to a shift, stencils for the interior nodes can vary.

\subsection{Conservation of the Filtering Procedure}
An important requirement of filters being applied in the EC context is that they are globally conservative, i.e. 
\begin{equation}
\mathbf{1}^T\mathbf{M}\mathbf{u}^{\star}=\mathbf{1}^T\mathbf{M}\mathbf{u},
\end{equation}
where here $\mathbf{M}$ represent the global quadrature rule accounting for Jacobians of local to reference mappings: $\mathbf{1}^T\mathbf{M}\mathbf{1}=|\Omega|$. Following the lead of Vreman \cite{vreman}, we describe one way of constructing a conservative correction for a non-conservative filter. For arbitrary initial data, we wish 
\begin{align*}
    \mathbf{1}^T\mathbf{M}\mathbf{u}^{\star}-\mathbf{1}^T\mathbf{M}\mathbf{u}&= \mathbf{1}^T\mathbf{M}(\mathbf{K}-\mathbf{I})\mathbf{u}=0.
\end{align*}
Now for a non-conservative filter this last equality may not hold, hence subtracting this discrepancy in the form of a properly scaled $\mathbf{1}^T\mathbf{M}(\mathbf{K}-\mathbf{I})$ term from the original matrix will yield conservation: 
\begin{equation}
\mathbf{K}_{corr}=\mathbf{K}-\frac{\mathbf{1}\mathbf{1}^T}{\mathbf{1}^T\mathbf{M}\mathbf{1}}  \mathbf{M}(\mathbf{K}-\mathbf{I}).
\end{equation}
Note that here
\begin{align*}
    \mathbf{1}^T\mathbf{M}\mathbf{K}_{corr}&= \mathbf{1}^T\mathbf{M}\left(\mathbf{K}-\frac{\mathbf{1}\mathbf{1}^T}{\mathbf{1}^T\mathbf{M}\mathbf{1}}  \mathbf{M}(\mathbf{K}-\mathbf{I})\right)\\
    &=\mathbf{1}^T\mathbf{M}\mathbf{K}-\mathbf{1}^T\mathbf{M}(\mathbf{K}-\mathbf{I})\\
    &=\mathbf{1}^T\mathbf{M},
    \end{align*}
    hence
    \[\mathbf{1}^T\mathbf{M}\mathbf{K}_{corr}\mathbf{u}=\mathbf{1}^T\mathbf{M}\mathbf{u}\]
    for any $\mathbf{u}$ and conservation is achieved. 

    Note that there are other forms that the correction could take such as, for example,
    \[
\mathbf{K}_{corr,2}=\mathbf{M}^{-1}\left(\mathbf{K}-\frac{\mathbf{1}\mathbf{1}^T}{\mathbf{1}^T\mathbf{1}} \cdot (\mathbf{K}-\mathbf{M})\right).
    \]
This particular choice has deleterious effects on the spectrum of the corrected operator, so we do not pursue it here. As an example of the issue this second correction causes, consider the constant mode $\mathbf{1}$. Owing to the imposed consistency condition, our filters will reproduce constants i.e. $\mathbf{K}\mathbf{1}=\mathbf{1}$, and this is a desirable property in the corrected filter. Put another way, we wish $(1,\mathbf{1})$ to be an eigenpair of $\mathbf{K}_{corr}$ when it is an eigenpair of $\mathbf{K}$. We show here that the original $\mathbf{K}_{corr}$ preserves this property while $\mathbf{K}_{corr,2}$ does not in general.
\begin{proof}
   Supposing that $\mathbf{K}\mathbf{1}=\mathbf{1}$, note that $\mathbf{K}_{corr}$ is consistent as
    \begin{align*}
\mathbf{K}_{corr}\mathbf{1}&=\mathbf{K}\mathbf{1}-\frac{\mathbf{1}\mathbf{1}^T}{\mathbf{1}^T\mathbf{M}\mathbf{1}}  \mathbf{M}(\mathbf{K}-\mathbf{I})\mathbf{1}\\
&=\mathbf{K}\mathbf{1}-\frac{\mathbf{1}\mathbf{1}^T}{\mathbf{1}^T\mathbf{M}\mathbf{1}}  \mathbf{M}(\mathbf{K}\mathbf{1}-\mathbf{1})\\
&=\mathbf{K}\mathbf{1}-\frac{\mathbf{1}\mathbf{1}^T}{\mathbf{1}^T\mathbf{M}\mathbf{1}}  \mathbf{M}(\mathbf{1}-\mathbf{1})\\
&=\mathbf{K}\mathbf{1}=\mathbf{1}.
    \end{align*}
    Considering $\mathbf{K}_{corr,2}$ however, we have
\begin{align*}
\mathbf{K}_{corr,2}\mathbf{1}&=\mathbf{M}^{-1}\left(\mathbf{K}\mathbf{1}-\frac{\mathbf{1}\mathbf{1}^T}{\mathbf{1}^T\mathbf{1}}  (\mathbf{K}-\mathbf{M})\mathbf{1}\right)\\
&=\mathbf{M}^{-1}\left(\mathbf{1}-\frac{\mathbf{1}\mathbf{1}^T}{\mathbf{1}^T\mathbf{1}}  (\mathbf{1}-\mathbf{M}\mathbf{1})\right)\\
&=\mathbf{M}^{-1}\frac{\mathbf{1}\mathbf{1}^T}{\mathbf{1}^T\mathbf{1}}\mathbf{M}\mathbf{1}.
\intertext{Unless $\mathbf{M}$ is a constant multiple of the identity, the above expression does not equal $\mathbf{1}$.}
\end{align*}
\end{proof}
The discussion above does not mean to suggest that $K_{corr,2}$ will not have eigenvalues equal to one, simply that the corresponding eigenvector will not be a constant mode as hoped. By means of numerical justification, consider the global averaging filter $\mathbf{K}_{glob.avg}=\frac{\mathbf{1}\mathbf{1}^T}{\mathbf{1}^T\mathbf{1}}.$ Note that this filter is consistent $(\mathbf{K}_{glob.avg}\mathbf{1}=\mathbf{1}$), but not conservative for arbitrary quadrature rules, so a correction is appropriate. Here we consider as quadrature rule a 20 element LGL quadrature rule with 6 points per element. The eigenvalues of the filter remain unchanged under correction, with one eigenvalue remaining $1$ while the rest are $0$, but the eigenvector corresponding to $\lambda=1$ does change for the second correction see Figure $\ref{fig:eigenvctor_corr}$. As the second correction violates the consistency requirement, it is unsuitable for present purposes.

\input{figures/files/Global_filt_corr_comp}

\subsection{Fourier and Spectral Response}

For periodic boundary conditions, the spectral effects of the continuous SIAC filtering operation can be assessed directly from the Fourier transform of the kernel given by {\cite{SIAC_Fourier}}
\begin{equation}
    \hat{K}^{(r+1,\ell)}_H(k)=\Bigg(\frac{\sin\big(\frac{kH}{2}\big)}{\frac{kH}{2}}\Bigg)^{\ell}\Bigg(c_{\ceil{r/2}+1}+2\sum_{\gamma=1}^{\ceil{r/2}}c_{\gamma}\cos(\gamma k H)\Bigg).
\end{equation}
Plotted in Figure \ref{fig:1D_Fourier_Kernel} we see the effects of various parameterizations of the SIAC kernel on its Fourier response. Specifically, varying $H$ dilates the response, while increasing the smoothness of the spline $\ell$ leads to a faster decay of the Fourier transform as the underlying kernel is smoother. Lastly, increasing the number of moments sharpens the response by making it flatter near zero, specifically the response has its first $r$ derivatives equal to $0$ at $k=0$.

\input{figures/files/1D_Fourier_Kernel}

\input{figures/files/1D_eigspectrum}

The truncation of the continuous filtering procedure required to produce a discrete filtering matrix may influence the spectral properties and stability of the resulting discrete filtering procedure. To that end, we consider for varying parameterizations and corrections the eigenvalues of the filtering matrices. As these matrices are non-symmetric, they may still exhibit temporary transient growth even when their eigenvalues are bound above in magnitude by unity \cite{trefethen2005spectra}. To address this, we also consider the singular values of the resulting filtering matrices which are plotted for varying moments, spline orders, and scalings in Figure \ref{fig:1D_spectrum_Kernel}. In all cases the maximum singular values are greater than 1 so transient growth may occur if the filter were applied in an iterative fashion, but that is not a concern in the correction procedure. Similar to the continuous case, increasing the number of moments, i.e. polynomials reproduced by the filter, results in eigenvalues moving closer to $1$. For periodic BCs, only the "periodic" constant polynomial is reproduced exactly. Increasing the spline order and thus the smoothness of the kernel increases the damping capabilities of the filter and causes eigenvalues to gravitate towards 0. Similarly, increasing the scaling also increases the damping capabilities of the filter and moves eigenvalues away from 1. All that is said for the eigenvalues also applies to the singular values, except that they are strictly non-negative and not bounded above by unity.

%% file: figures/files/LSIAC_Example.tex
\begin{figure}
    \centering
    \begin{tabular}{c c c}
       \includegraphics[width=0.4\linewidth]{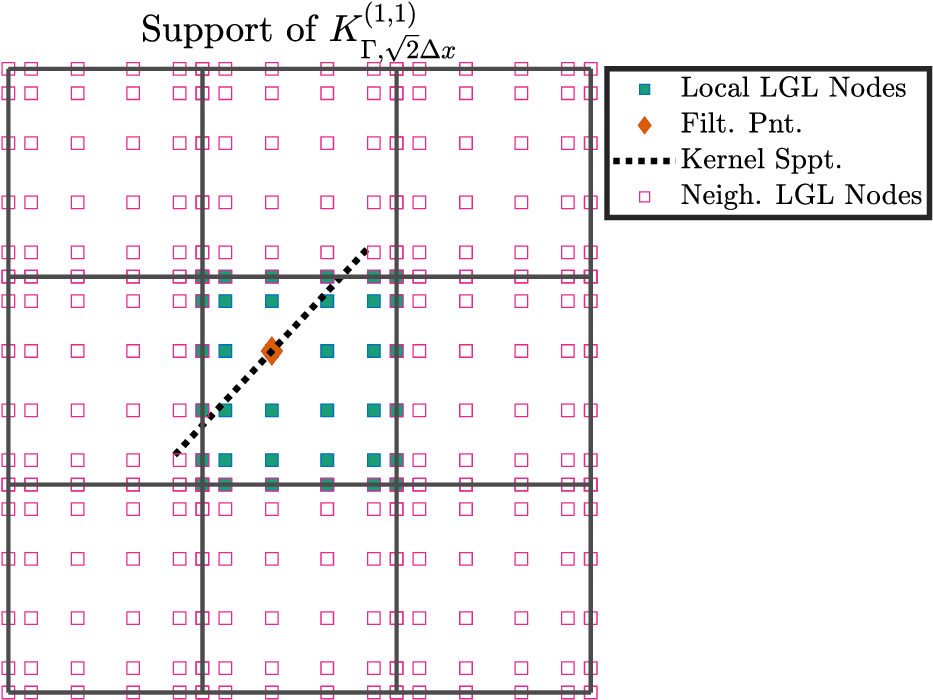}    &    \includegraphics[width=0.27\linewidth]{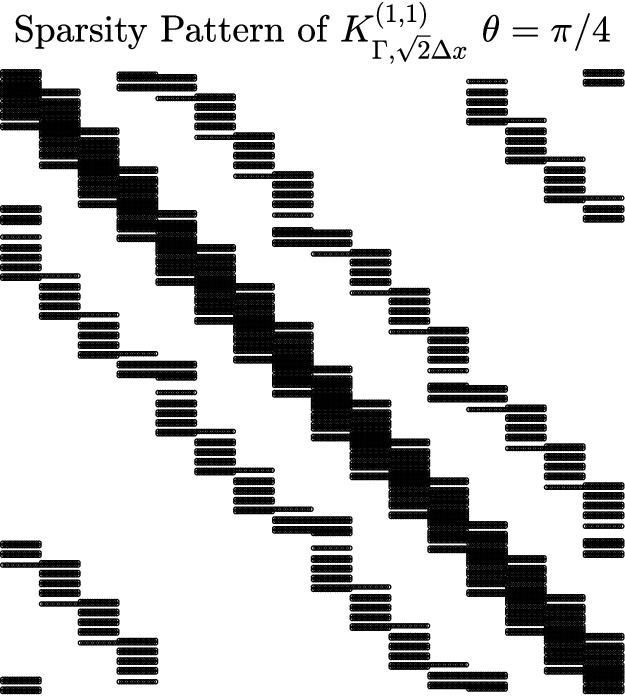}   
        \end{tabular}

    \caption{Left: Depiction of the LSIAC filter procedure for the $K^{(1,1)}_{\Gamma,\Delta x}$ filter with $\theta=\pi/4$ applied to $p=5$. Depiction of the sparsity pattern of the resultant $K^{(1,1)}_{\Gamma,\Delta x}$ filtering matrix for $p=3$ on a $4\times 4$ uniform mesh with $\theta=\pi/4$ (Right).}
    \label{fig:LSIAC_Example}
\end{figure}

%% file: figures/files/Global_filt_corr_comp.tex
\begin{figure}
    \centering
    \includegraphics[width=0.75\linewidth]{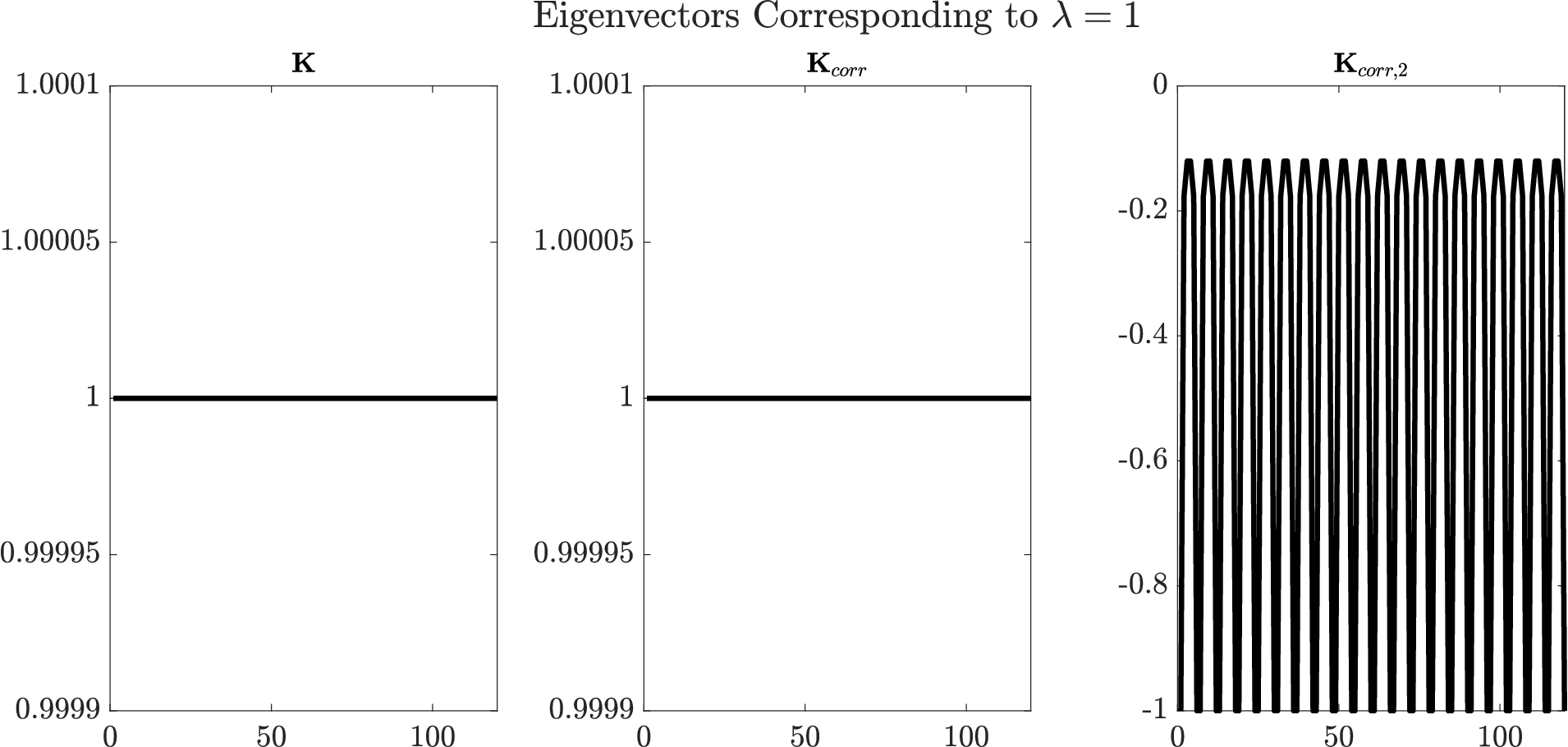}
    \caption{Eigenvector corresponding to $\lambda=1$ for varying corrections of $\mathbf{K}_{glob.avg}$. Note here a $20$ element mesh using a 6 node per element LGL quadrature matrix $\mathbf{M}$ was used. Note the eigenvectors $\mathbf{v}$ have been normalized so that $\norm{\mathbf{v}}_{\ell^{\infty}}=1$.}
    \label{fig:eigenvctor_corr}
\end{figure}

%% file: figures/files/1D_Fourier_Kernel.tex
\begin{figure}
    \centering
    \begin{tabular}{c c c}
    \multicolumn{3}{c}{Kernel Fourier Response}\\
       \includegraphics[width=0.33\linewidth]{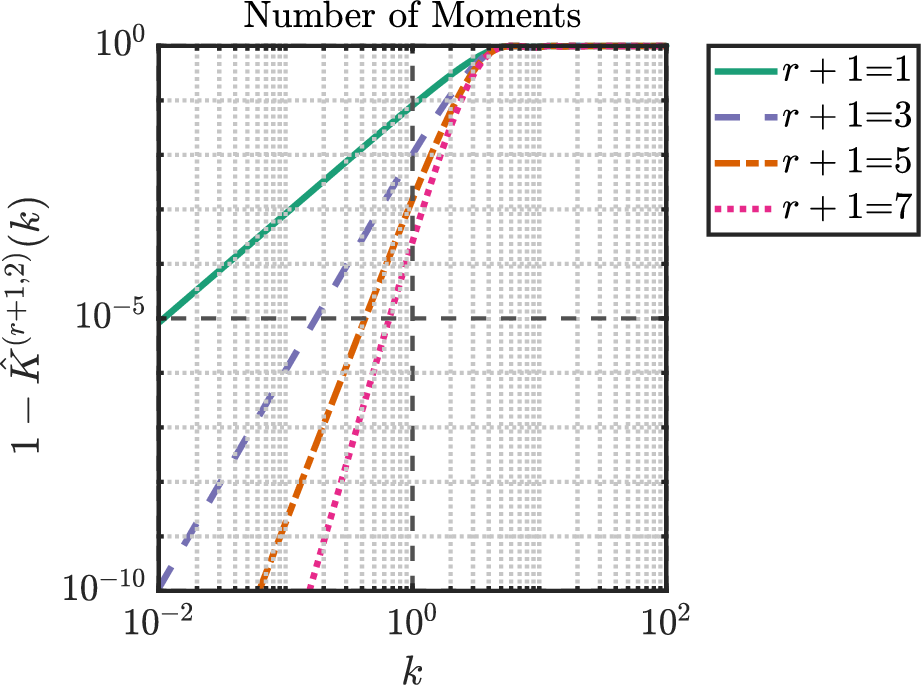}    &    \includegraphics[width=0.33\linewidth]{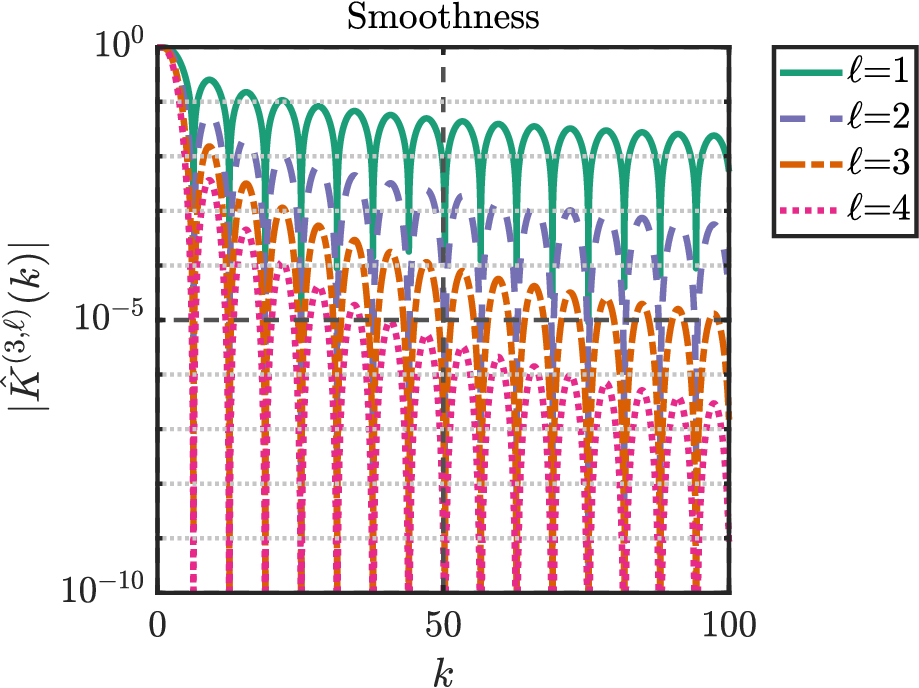}  &    \includegraphics[width=0.33\linewidth]{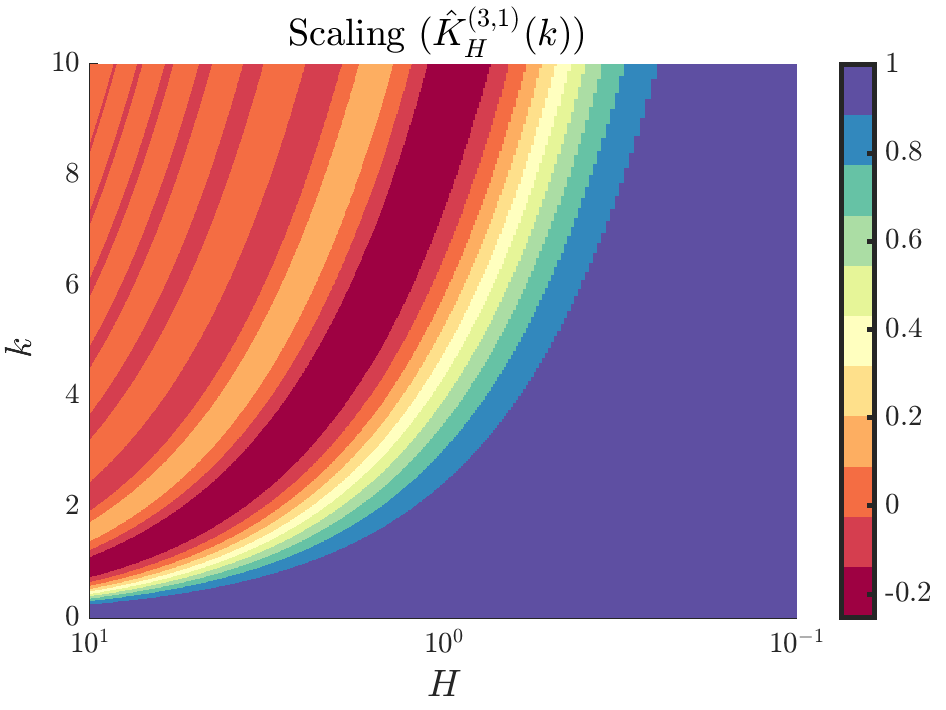}
        \end{tabular}

    \caption{Comparison of the effect of number of moments, spline order, and kernel scaling on the Fourier response of the SIAC filter.}
    \label{fig:1D_Fourier_Kernel}
\end{figure}

%% file: figures/files/1D_eigspectrum.tex
\begin{figure}
    \centering
    \begin{tabular}{c c c}
      Moments & Spline Order & Scaling\\ \includegraphics[width=0.33\linewidth]{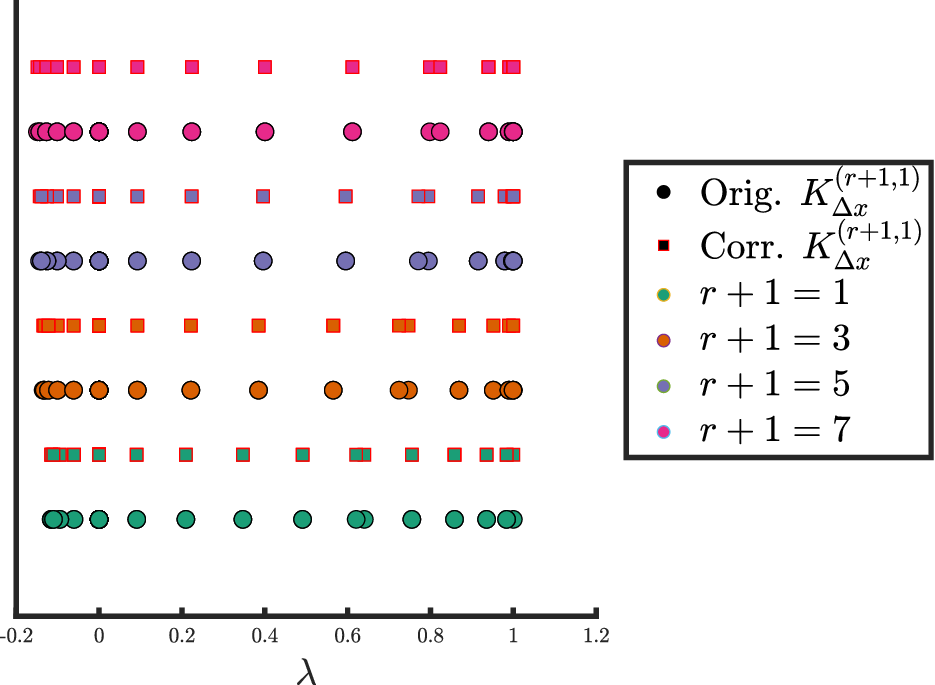}    &    \includegraphics[width=0.33\linewidth]{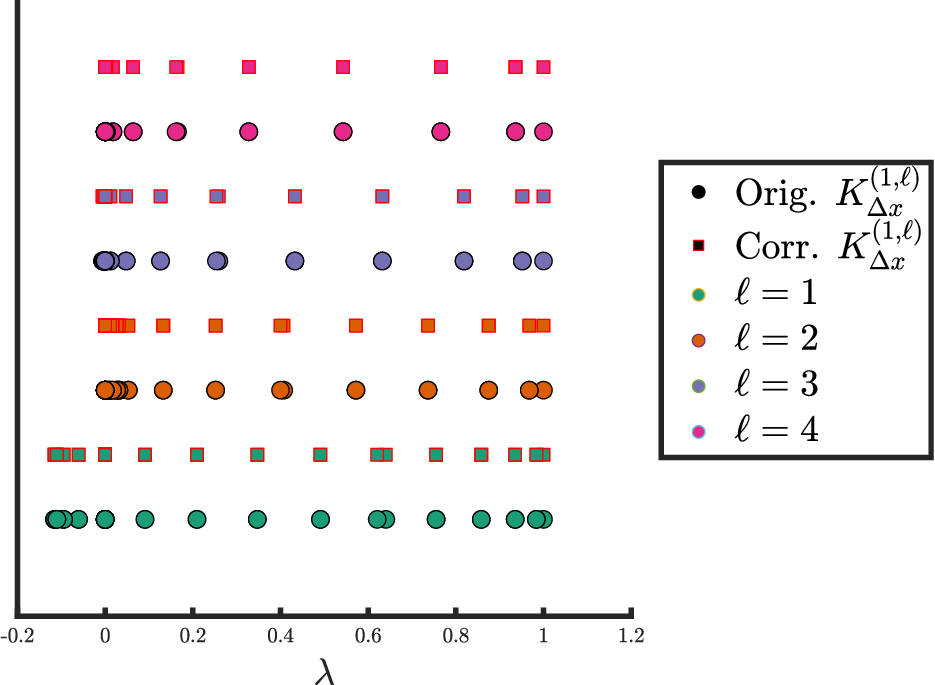}  &    \includegraphics[width=0.33\linewidth]{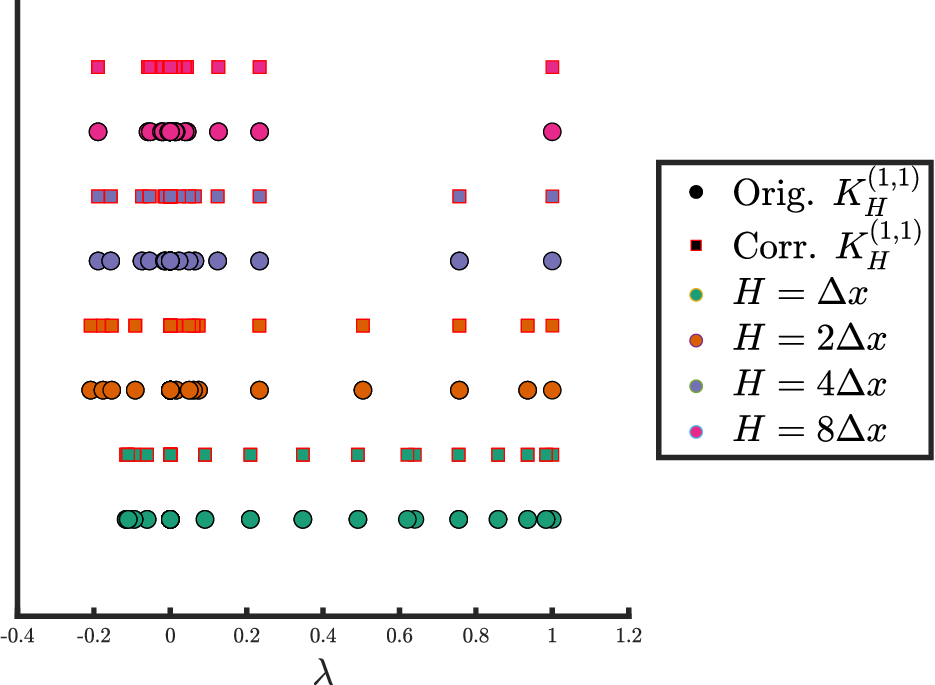}  \\
         \includegraphics[width=0.33\linewidth]{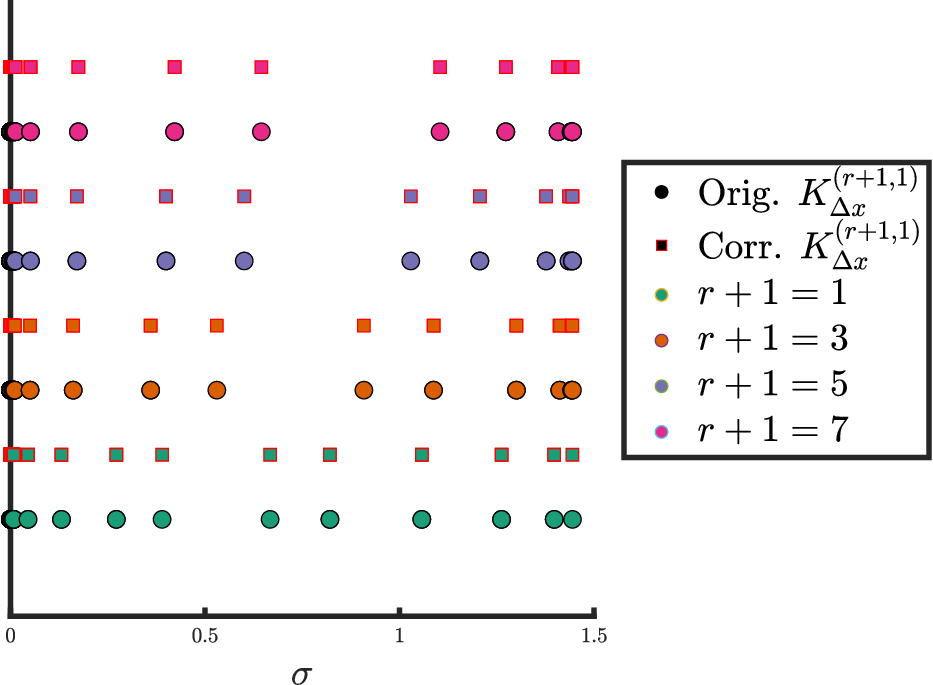}    &    \includegraphics[width=0.33\linewidth]{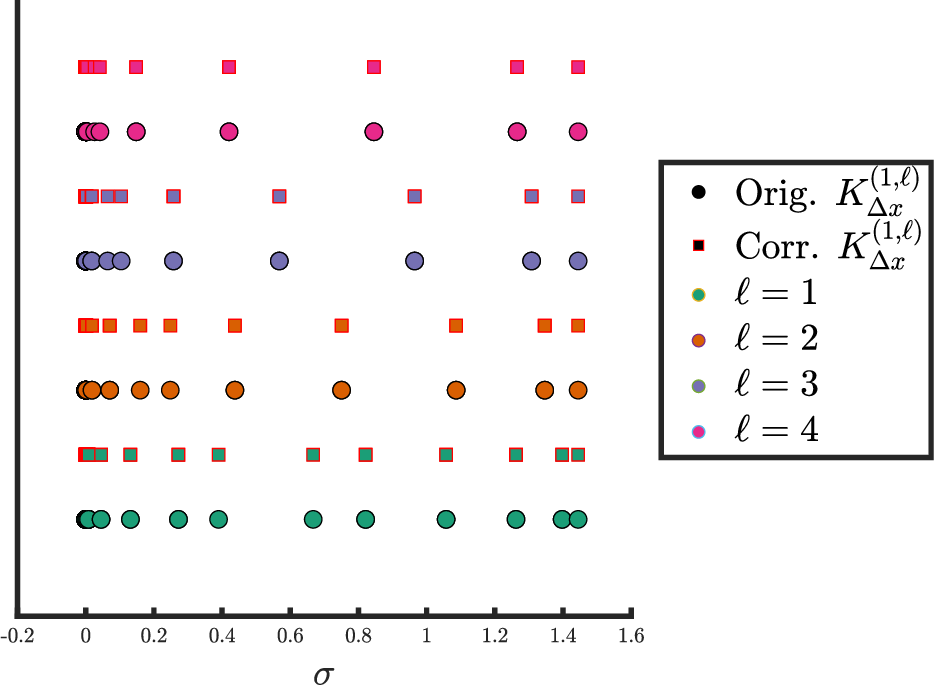}  &    \includegraphics[width=0.33\linewidth]{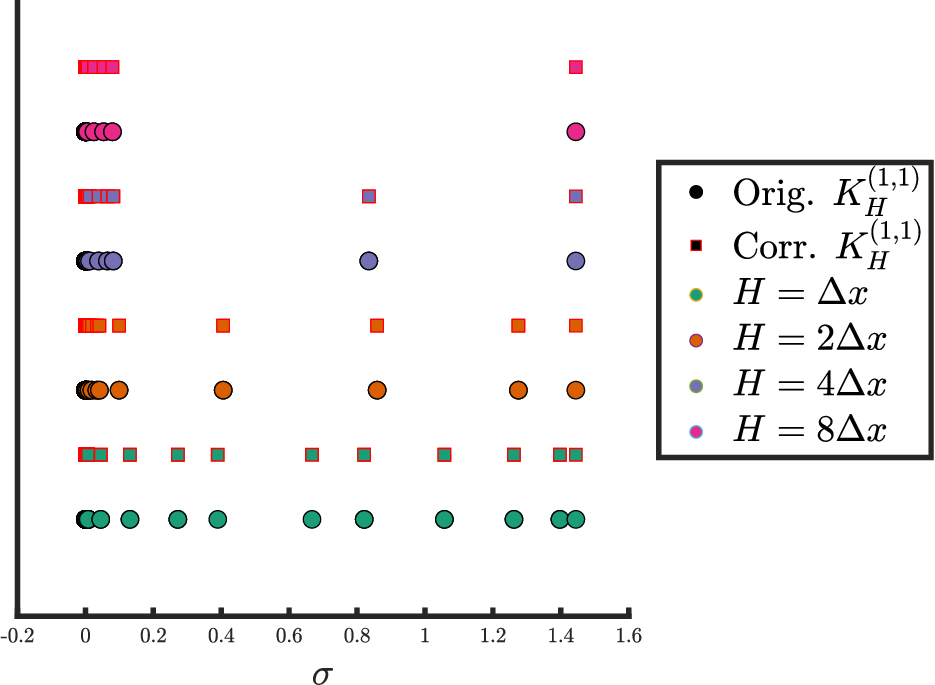}  \\
        \end{tabular}

     \caption{Comparison of the effect of number of moments, spline order, and kernel scaling on the spectrum and singular values of the discrete filtering matrix for a $p=3$ approximation on a $N_{elm}=10$ mesh. Note only real values of the eigenvalues are plotted here as the imaginary components are, if not 0, on the order of machine precision.}
    \label{fig:1D_spectrum_Kernel}
\end{figure}

%% file: sections/results.tex
\section{Results}\label{sec:num_res}
In this section we demonstrate the applicability of SIAC-based entropy-correction terms for energy conservation of the Burgers equation, as well as solution regularization. We also investigate the impact of different SIAC correction parameterizations. Additionally,  the SIAC corrections are shown to preserve the original order of accuracy of the underlying method.

\subsection{One-dimensional entropy conservation}
 In order to compare the performance of different correction methods, we consider their application to the inviscid Burgers problem $u_t+(u^2/2)_x=0$ on the periodic domain $\Omega=[0,2]$ subject to initial data $u(x,0)=\sin(\pi x)+0.01$. For this problem, the characteristics of the PDE cross and a shock forms at $x=1$ at $t=1/\pi$. 
 
 First, the choice in numerical {interface} flux used in the method is evaluated as a baseline. A non-dissipative central flux $f_{cent.,\;e+1/2}^{num}=\frac{1}{2}(f(u_{p}^e)+f(u_0^{e+1}))$ is compared to a slightly-more stabilizing Local Lax-Friedrichs (LLF) formulation $f^{num}_{LLF,e+1/2}=\frac{1}{2}(f(u_{p}^e)+f(u_0^{e+1}))+\frac{1}{2}\alpha_{LLF}(u_{p}^e-u_{0}^{e+1})$. Results are depicted in Figure \ref{fig:1D_cons_cent_vs_LLF} for a $p=5$ approximation over a $N_{elm}=21$ element mesh. In terms of solution quality we see that the approximations using the central flux are less desirable, and also less stable. In fact, the uncorrected central flux scheme crashes before reaching $t=1.5/\pi$, {at $t\approx 1.365/\pi$ while the uncorrected LLF scheme crashes immediately after ($t\approx1.507/\pi$)}. In terms of the subcell entropy error metric, the LLF approaches also perform better than their central flux counterparts. 
 Since the LLF does not appear to introduce numerical stiffness into the rRK-LS procedure, the LLF DGSEM scheme is exclusively used as a the baseline scheme  in the experiments to follow, although the witnessed trends are applicable to alternative baseline flux choices.

\input{figures/files/1D_centr_vs_LLF_flux}

Having fully specified the baseline scheme, the effects of different correction formulations are assessed for a $p=5$ and $N_{elm}=21$ element approximation. In Figure \ref{fig:1D_cons} the performance of the $K^{(1,1)}_{\Delta x}$ filter correction is contrasted with the Local approach and the LS implementation that weighs adaptively between a convex combination of the SIAC filter and the Local average formulation. These approaches are using a LLF flux in the DGSEM method and are also compared with the Skew-Symmetric Entropy-Conserving method of Gassner \cite{Gassner_Skew} which is an energy conserving {method making use of energy conserving interface fluxes and a balancing of volumetric flux terms}. As depicted, similar performance is observed with respect to the subcell entropy metric\footnote{Note that the Skew-Symmetric EC method does conserve energy at the subcell level; however, the form of this local conservation is not equivalent to the current subcell error metric $D_h$. Significant flexibility in the choice of subcell entropy metrics exists, and other options are certainly feasible. For all tests here the same subcell metric was applied to enable comparison of performance.}. Only after shock formation do any of the methods differentiate in their behavior. As expected, after this time the uncorrected scheme begins to blow up.

With respect to the relative correction magnitude, the SIAC approach has a lower relative magnitude then the other formulation prior to shock formation. The unconstrained LS approach has a much larger magnitude owing to $\theta$ being unconstrained which can increase the magnitude of $\mathbf{c}$. Looking at how the weighting parameter $\theta$ changes over time, we see that initially neither correction is favored in reducing $\norm{\mathbf{D}_h}_{\mathbf{M}}$, but eventually the SIAC correction becomes favorable just prior to shock formation and remains so after shock formation until eventually oscillations between the two corrections occur. In terms of the solution quality, we see that at shock formation the different approximations performs similarly, but afterwards the conservation of energy is no longer mathematically required and so this lack of decay introduces oscillations into the solution. These oscillations develop on the right side of the shock before eventually spreading throughout the profile. The energy-conserving Skew-Symmetric formulation produces oscillations on both sides of the shock profile, which over time also disperse through the profile. We wish to note that the uncorrected DGSEM method using the LLF flux (not plotted) begins to blow up following shock formation and crashes prior to the final time of $t=2/\pi$, further justifying the need for augmentation of the baseline scheme.

\input{figures/files/1D_Cons_Figure}

\newpage 
\subsection{One-dimensional regularization}
Taking the one-dimensional shock problem and adding a dissipative term to the correction as describe in section \ref{subsec:reg}, we next investigate the performance of different filter parameterizations: {moments, spline order, scaling}, and compare the local correction approach to SIAC-based corrections.

For all regularization tests, a $p=5$ approximation on $N_{elm}=21$ elements is evolved past shock formation to a final time of $t=5/\pi$. In comparing the qualitative performance of different SIAC-correction parameterizations the tuning parameters are held constant at $c_E=10$ and $c_{max}=1$. In Figure \ref{fig:1D_param_sol} (Left) increasing the number of moments, which is the number of splines, leads to degradation of solution quality at sharp gradients. Meanwhile increasing the spline order (i.e., increasing the smoothness of the kernel) slightly increases the energy dissipation of the regularization, which can be observe in Figure \ref{fig:1D_param_energy_mass} (Center). The kernel scaling has the greatest effect on dissipation, with greater scalings leading to more energy dissipation as depicted in the solution profiles Figure \ref{fig:1D_param_sol} (Right) and the normalized energy plots \ref{fig:1D_param_energy_mass} (Right). Lastly, we note that all parameterizations conserve the domain integral of $u$ up to floating point arithmetic as expected  (see Figure \ref{fig:1D_param_energy_mass} , row 1).

\input{figures/files/1D_param_sol_comp}

\input{figures/files/1D_param_mass_eng_comp}

\input{figures/files/FV_ref}

Comparing the local approach of Abgrall to a SIAC-based correction, we investigate the solution quality as well as dissipation magnitudes for different scalings. Here the tuning parameters are chosen so as to ensure qualitatively good solution profiles as compared to a $N_{elm}=2000$ fine mesh finite-volume reference solution depicted in Figure \ref{fig:FV reference}. Owing to the difficulty in finding quantitative measures of qualitative agreement, the performance of the two approaches is distinguished by applying the same tuning parameters to each. Each correction type is tuned for ``best" qualitative performance. Then those same tuning parameters are applied to the other correction in order to then compare the numerical solutions. {The choice of tuning parameters could perhaps be automated by machine learning type approaches, but that is outside the scope of this work.}

In tuning the local correction approach, it was found that $c_E=4$ and $c_{max}=0.15$ produced the best early time profiles. These dissipation parameters were then applied to both the SIAC and local corrections and metrics collected. Note that in terms of solution quality, parameter combinations that damped oscillations while preventing eventual flattening of the solution in the element containing the shock were elusive for the local approximation. The solution profiles in Figure \ref{fig:1D_local_vs_SIAC_local_params} (Left) depict this, with the SIAC correction performing better in damping oscillations and preventing their spread into other elements near the shock. The local approach does not fair as well. As discussed previously, the solution quality degrades eventually as the solution in the element containing the shock flattens and oscillations appear in its neighbors. The gradient-based sensor performs well in both cases. Figure \ref{fig:1D_local_vs_SIAC_local_params} (Center) shows how $\phi_{diss}$ varies spatially. Noting that the local approach requires a single dissipation magnitude for all elements, we see that the magnitude of dissipation increases when the shock begins to form at $t=1/\pi$, and in the local case, the greatest dissipation occurs in the elements containing the shock, and its neighbors. Both methods conserve the integral of $u$, and have similar magnitudes of energy dissipation over time Figure \ref{fig:1D_local_vs_SIAC_local_params} (Right), with both being more dissipative than the finite-volume reference solution.

\input{figures/files/1D_local_vs_SIAC_local_params}

 The SIAC correction is much less sensitive to the choice of tuning parameters and after minimal tuning, $c_E=10$ and $c_{max}=1$ are selected to produce solution profiles with minimal oscillations. Applying these dissipation parameters to both correction schemes, the solution profiles in Figure \ref{fig:1D_local_vs_SIAC_SIAC_params} (Left) were obtain. The SIAC approach performs better, with an oscillation free profile over time. The local approach eventually devolves into a stair-stepping  approximation, with flattened local approximation spreading out from the center element over time. This is caused by the increasing dissipation applied to those elements as shown in Figure \ref{fig:1D_local_vs_SIAC_SIAC_params} (Center). Comparatively, the dissipation magnitudes of the SIAC correction remain similar to the local-tuned case. While both methods still conserve the integral of $u$, the accuracy to which the local approach does so is inhibited, perhaps by floating point arithmetic effect. Similar magnitudes of energy dissipation over time are obtained, see Figure \ref{fig:1D_local_vs_SIAC_SIAC_params} (Right). The SIAC approach is more dissipative in later times than the local approach, and both remain more dissipative than the finite-volume reference.

\input{figures/files/1D_local_vs_SIAC_SIAC_params}

\subsection{Two-dimensional regularization}
Next the correction methodology is applied to two-dimensional inviscid Burgers problem:
\begin{equation*}
    \frac{d}{dt}u+\nabla \cdot \Big(\frac{u^2}{2}(1,1)\Big)=0,
\end{equation*}
on the domain $\Omega=[0,1]^2$ subject to periodic boundary conditions and initial data $u(x,y,0)=\sin(2\pi(x+y))+0.01$. This problem develops a shock in finite time at $t=\frac{1}{2\pi}$ along the lines $x+y=\frac{1}{4},\frac{3}{4}$ which introduce spurious oscillations in high-order polynomial approximations.

In solving this problem, we wish to stress that the correction methodology does not necessitate a specific formulation of the DG method. As such, the current PDE is discretized in space by a modal DG-style approximation \cite{shu_cockburn} using a Local Lax-Friedrichs flux, though time-stepping is still handled by the SSPRK(3,3) rRK method.

\input{figures/files/2D_reg_tests}

\input{figures/files/2D_reg_cross_section}

Figure \ref{fig:2D_Reg_Example}plots the solutions at the time of shock formation $t=\frac{1}{2\pi}$ and then at a later time $t=\frac{1}{\pi}$ using a $N_{elm}=8\times8=64$ element mesh and a $p=5$ $Q_k$ tensor-product basis. For the regularized correction approach we consider the $K^{(1,1)}_{\Gamma,\sqrt{2}\Delta x}$ filter and choose dissipation parameters $c_E=0.1$ and $c_{max}=0.15$. Here the goal is to demonstrate improved solution quality under the regularization. Improved qualitative performance may be obtained by a more exhaustive searching of the tuning parameter space. Regardless, the solutions with the energy-dissipative correction are of higher quality with oscillations of reduced magnitude as compared to the uncorrected case. This is most readily apparent in the $y=x$ cross-section plots in Figure \ref{fig:2D_Reg_cross}. Note that both approximations conserve the integral of $u$, and the energy-dissipative correction causes greater dissipation of energy (Figure \ref{fig:2D_Reg_Example} (Right)). Further note that the dissipation of energy observed for the uncorrected scheme could be a result of the exact integration of volume terms in the modal DG scheme, compared to under-integration in the DGSEM case.

\subsection{Confirmation of high-order accuracy}
In \cite{Abgrall2018} it was shown that the local correction procedure would maintain the order of accuracy of the underlying spatial discretization method provided sufficiently high-order quadratures were applied. Here the spatial accuracy of the one-dimensional DGSEM method is verified for the inviscid Burgers equation under the correction. We also demonstrate the retention of the order of accuracy for corrected simulation of two-dimensional linear advection. Also in one-dimension we demonstrate numerically the energy decay caused by non-rRK time-stepping methods paired with the correction. 

Focusing on the original one-dimensional test problem first, the same initial data is considered as before, but this time the solution is evolved to a final time of $t=1/(2\pi)$ which is prior to shock formation. Therefore the solution is still smooth. We compute the $L^2-$ and $L^{\infty}-$errors for the approximation at 17 Gauss-Legendre points per element, where the exact solution is computed via Newton's method applied to the implicit solution $u(x,t)=u_0(x-u t)$. In Figure \ref{fig:spatial_conv_1D} and Table \ref{tab:spatial_orders} we compare the order of accuracy of the uncorrected scheme against the local and SIAC corrected schemes. Only the orders of accuracy of the $L^2$-errors are provided as the $L^{\infty}$-error orders are similar. In general the uncorrected scheme exhibits order of accuracy between $h^p$ and $h^{p+1}$. The SIAC corrected schemes maintain this order of accuracy. The local approach performs similarly except for the $p=1$ case when the method does not appear to converge even at order $h$. Furthermore, the SIAC correction does not necessitate a low moment and order formulation as the $K^{(3,2)}_{\Delta x}$ filter performs equally well as the $K^{(1,1)}_{\Delta x}$ filter.

\input{figures/files/Tabs_1D_Spatial_error}

\input{figures/files/1D_Spatial_Conv_tests}
In two dimensions we consider the linear advection problem of $u(x,0)=(\sin(\pi x )+0.01)(\sin(\pi y )+0.01)$ on the periodic domain $\Omega=[0,2]^2$ with constant velocity field $(1,1)$. Here we compute the errors from the exact solution $u(x,t)=u_0(x-t,y-t)$ at 36 Gauss-Legendre nodes per element. As shown in Figure \ref{fig:spatial_conv_2D} the order of accuracy of the scheme is again maintained, while energy (not depicted) under the $K^{(1,1)}_{\Gamma,\sqrt{2}\Delta x}$ correction is conserved in a manner analogous to the one-dimensional case.
\input{figures/files/2D_Spatial_Conv_tests}

Consider now the effects of the time-stepping discretization on the relative energy error. As the rRK time-stepping method ensures energy conservation, we instead focus on the time--stepping errors incurred by non-energy conserving baseline RK methods. We consider here Forward Euler, the optimal 2- and 3-stage Strong Stability Preserving Runge-Kutta (SSPRK) methods, SSPRK(2,2) and SSPRK(3,3) respectively, and the standard 4-stage 4th order Runge-Kutta method RK(4,4) (see \cite{SPMforTimeDep} for explicit formulations). In Figure \ref{fig:temporal_conv} we plot the relative energy conservation errors with respect to the time step size $\Delta t$. Under the uncorrected method, energy is not conserved by either the spatial or temporal discretization, and so no convergence is obtained under reductions in $\Delta t$. With respect to the spatially energy conserving discretizations, we observe convergence in relative energy error at the order of the time discretization, except for the SSPRK(2,2) method which exhibits 3rd order error. This  phenomena was previously observed in \cite{CAPUANO2017}.

\input{figures/files/1D_Temporal_Conv}

%% file: figures/files/1D_centr_vs_LLF_flux.tex
\begin{figure}
    \centering

\includegraphics[width=1.0\textwidth]{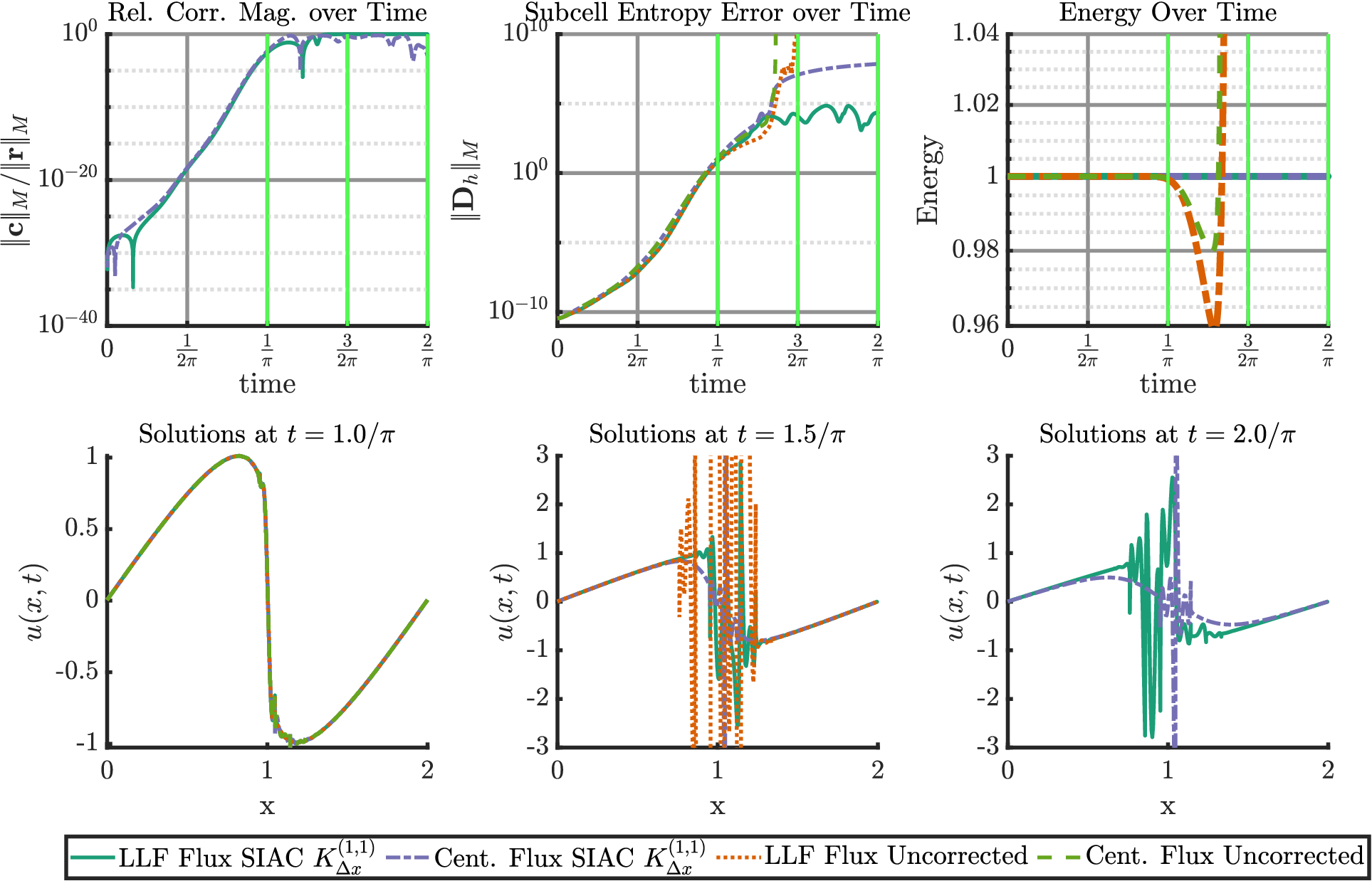}
    \caption{Comparison of SIAC corrected and uncorrected approaches using central fluxes and LLF fluxes. First row: Relative magnitude of the correction to the RHS, subcell entropy metric, and energy over time of the method. Second Row: Solutions at time of shock formation $t=\frac{1}{\pi}$ and then afterwards at $t=\frac{3}{2\pi}$ (the unplotted uncorrected central flux solution has crashed) and $t=\frac{2}{\pi}$ (Only corrected solutions remain). Note that the range of solutions is capped at $[-3,3]$ for $t=1.5/\pi$ and $t=2/\pi$ to better display solution behavior at the shock.}
    \label{fig:1D_cons_cent_vs_LLF}
\end{figure}

%% file: figures/files/1D_Cons_Figure.tex
\begin{figure}
    \centering
    \begin{tabular}{c c c}
       \includegraphics[width=0.33\linewidth]{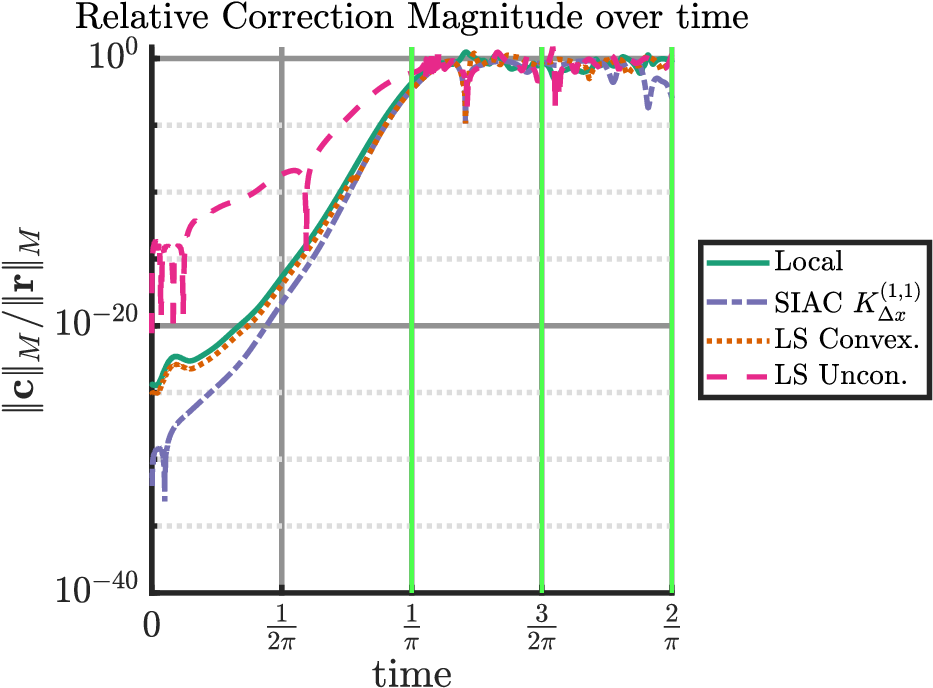}    &    \includegraphics[width=0.33\linewidth]{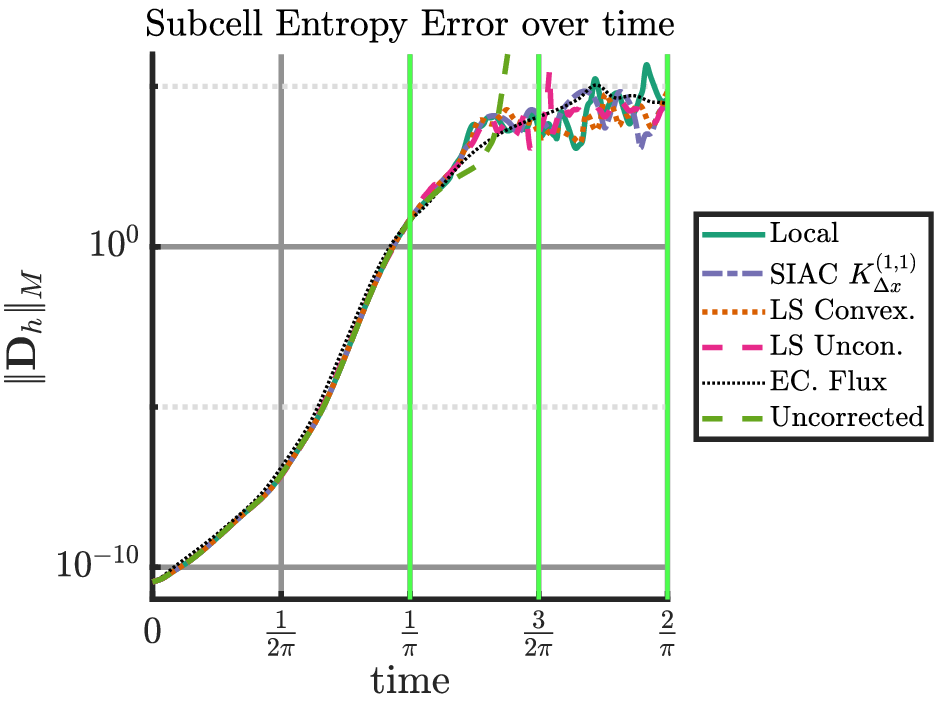}  &    \includegraphics[width=0.33\linewidth]{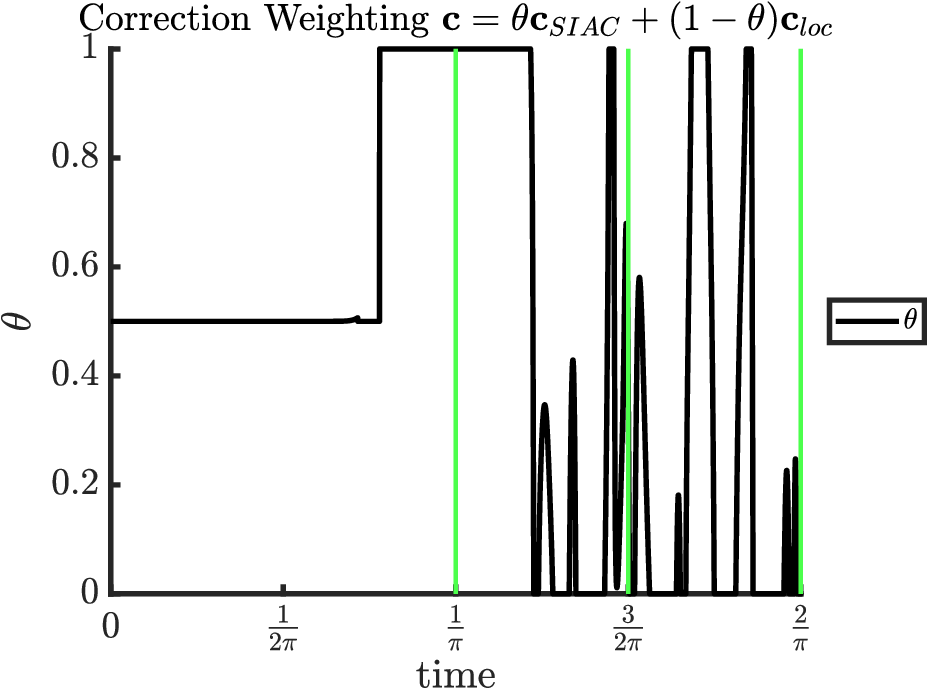} \\
         \includegraphics[width=0.33\linewidth]{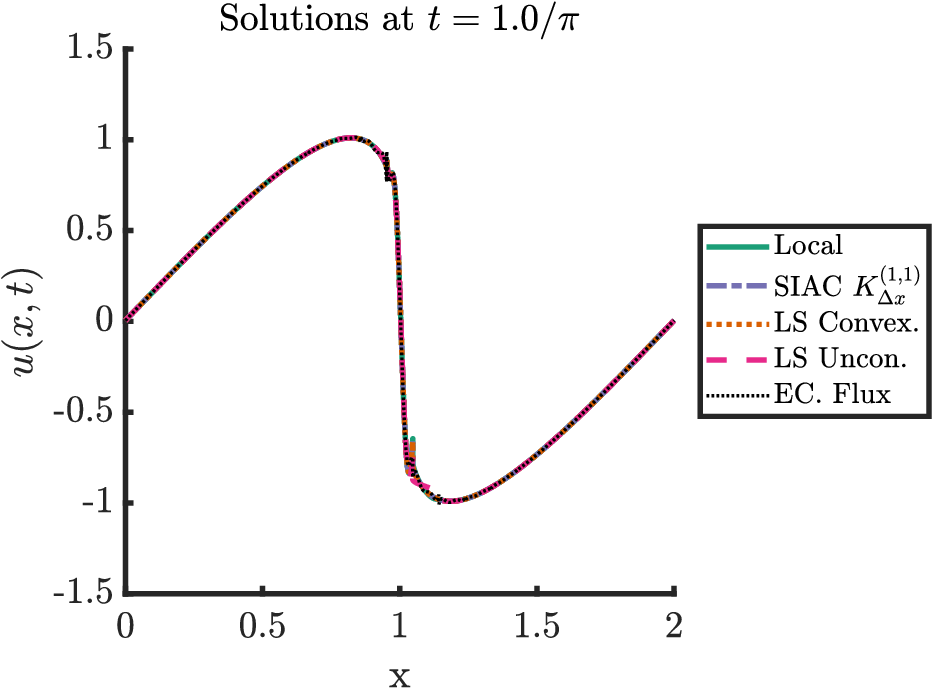}    & \includegraphics[width=0.33\linewidth]{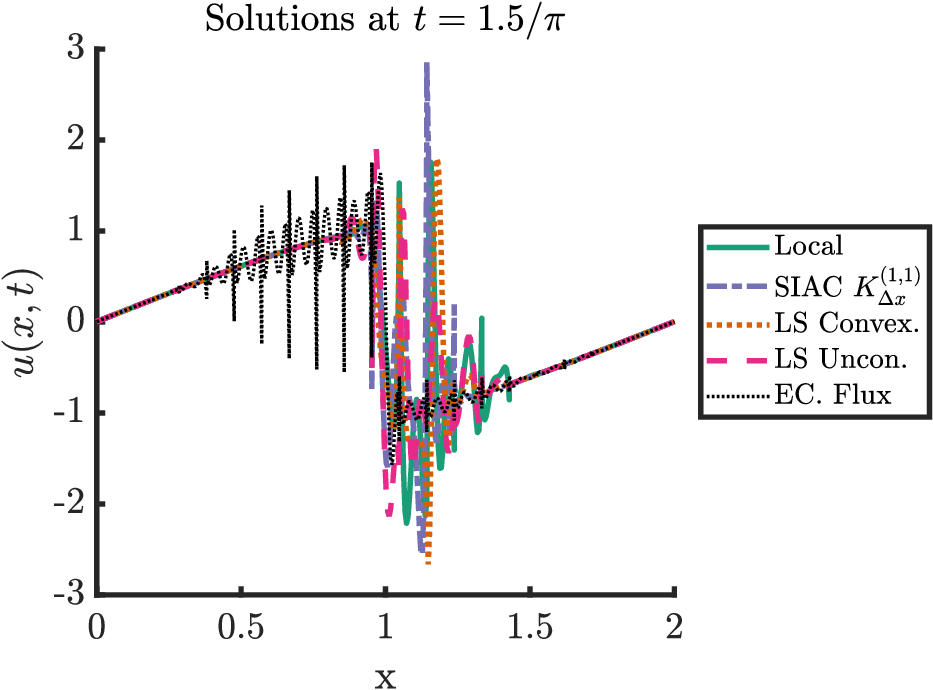} & \includegraphics[width=0.33\linewidth]{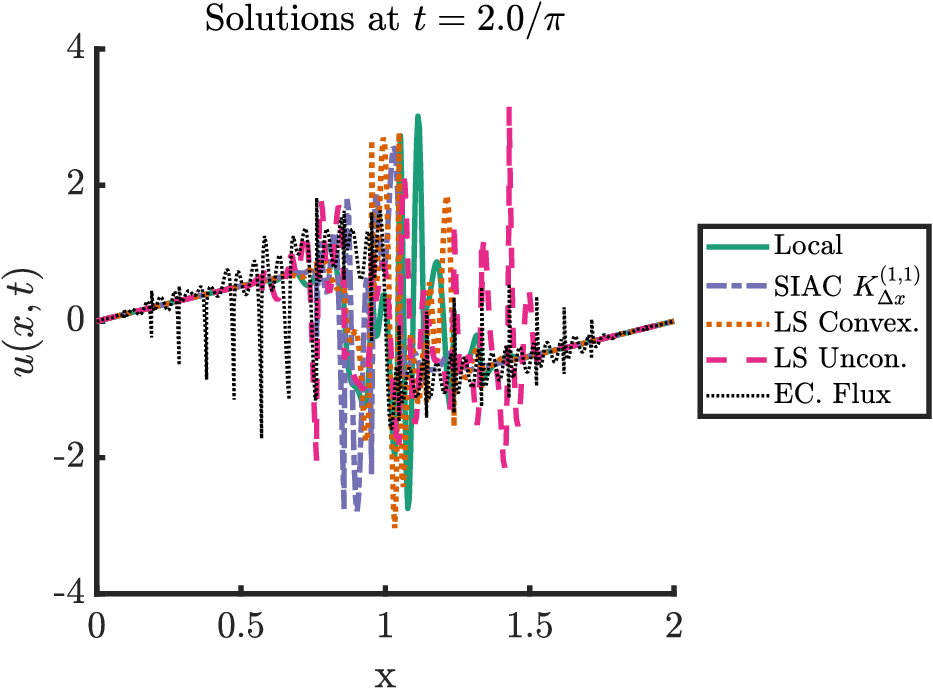} 
        \end{tabular}

    \caption{Comparison of SIAC and Local EC approaches, and their LS combinations with the EC Skew-Symmetric DGSEM formulation. First row: Relative magnitude of the correction to the rhs, subcell entropy metric, and weighing of constrained LS combination of corrections with $\theta=1$ being only a SIAC correction and $\theta=0$ being only a local correction). Second Row: Solutions at time of shock formation $t=\frac{1}{\pi}$ and then afterwards at $t=\frac{3}{2\pi}$ and $t=\frac{2}{\pi}$.}
    \label{fig:1D_cons}
\end{figure}

%% file: figures/files/1D_param_sol_comp.tex
\begin{figure}
    \centering
    \begin{tabular}{c c c}
       Moments $r+1=(1,3,5)$ & Spline Order $\ell=(1,2,3)$ & Kernel Scaling $H=(\Delta x,2\Delta x,4\Delta x)$\\
    \includegraphics[width=0.33\linewidth]{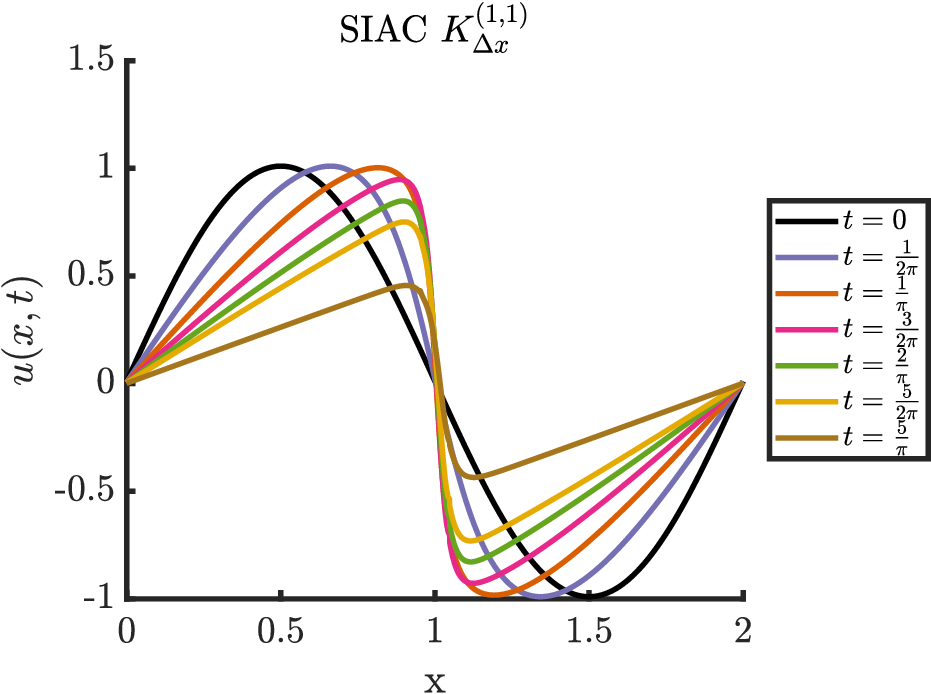}    &   \includegraphics[width=0.33\linewidth]{figures/images/1D_reg/reg_test_K11_dx_over_time.eps}     &   \includegraphics[width=0.33\linewidth]{figures/images/1D_reg/reg_test_K11_dx_over_time.eps}  \\
      \includegraphics[width=0.33\linewidth]{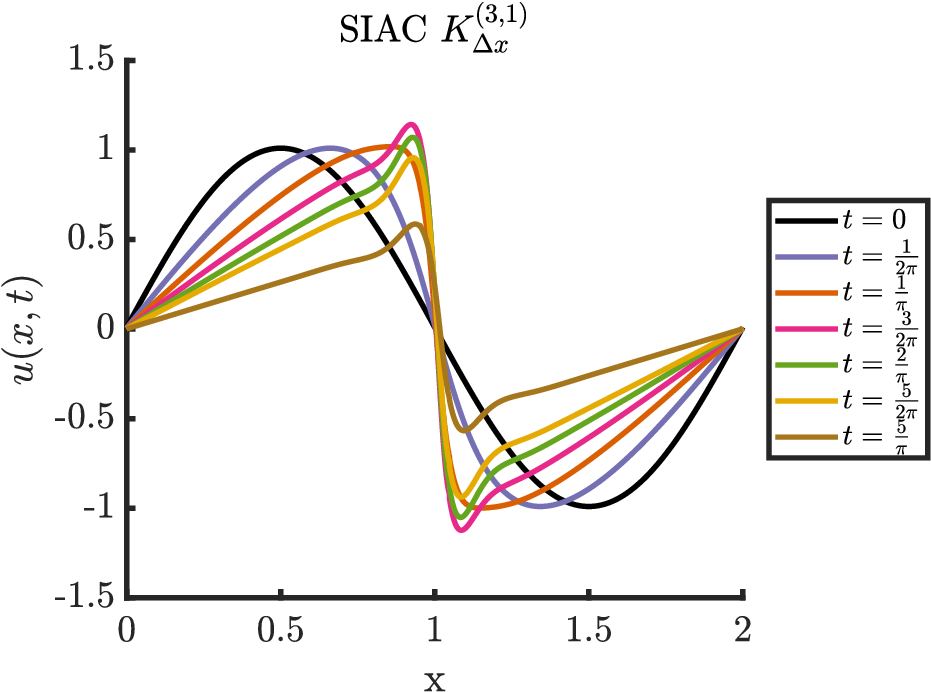}    &  \includegraphics[width=0.33\linewidth]{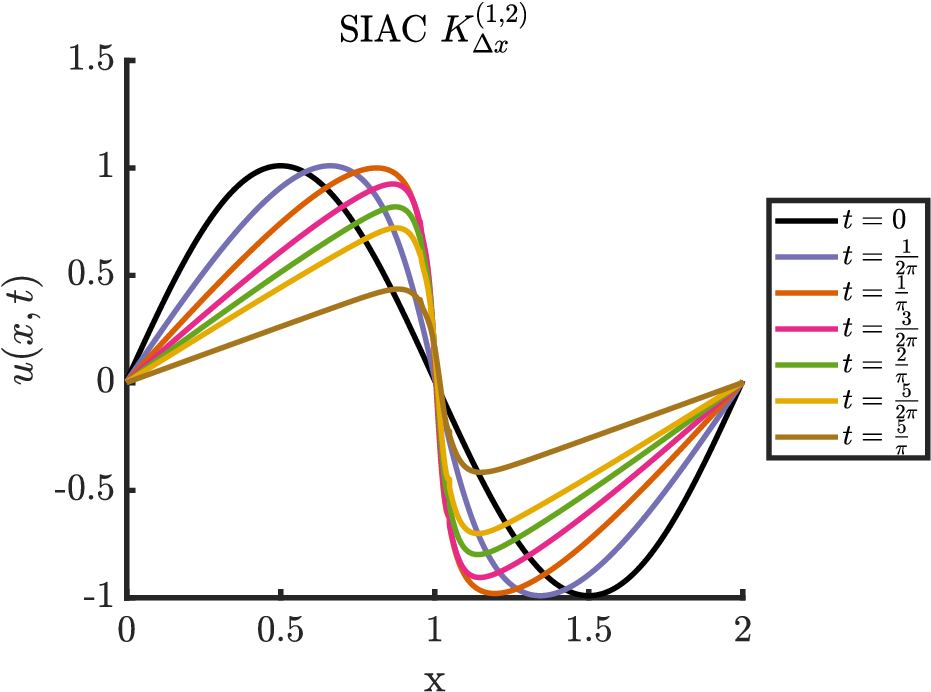}     &  
      \includegraphics[width=0.33\linewidth]{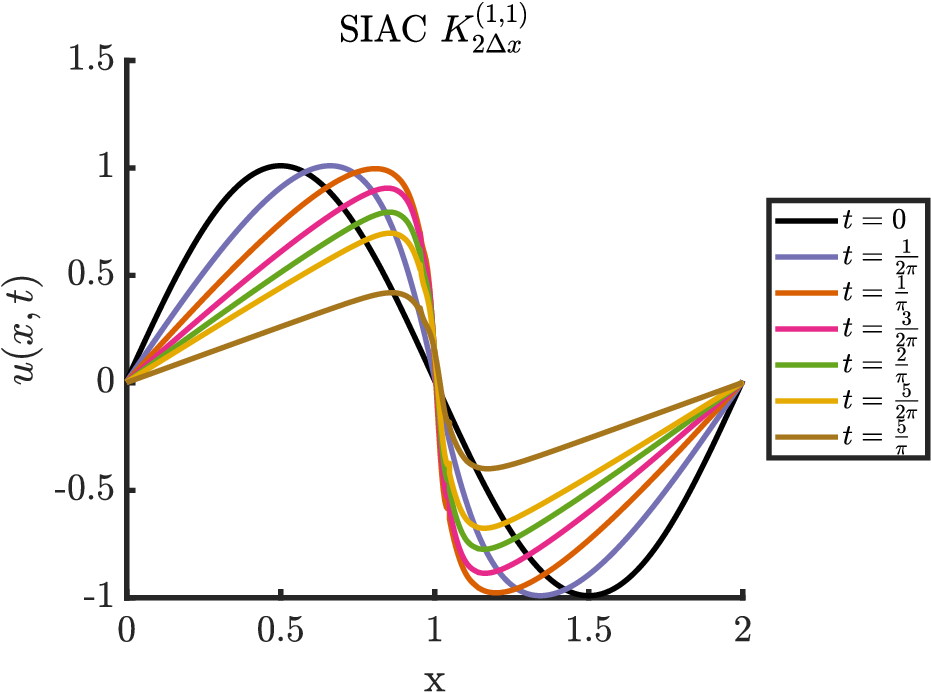} \\
          \includegraphics[width=0.33\linewidth]{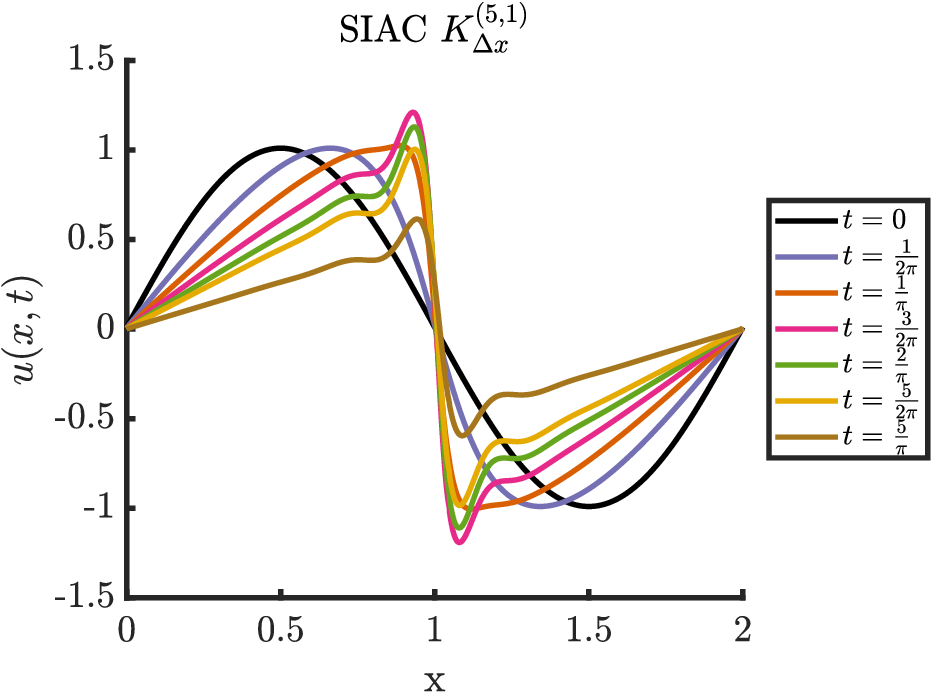}    &  \includegraphics[width=0.33\linewidth]{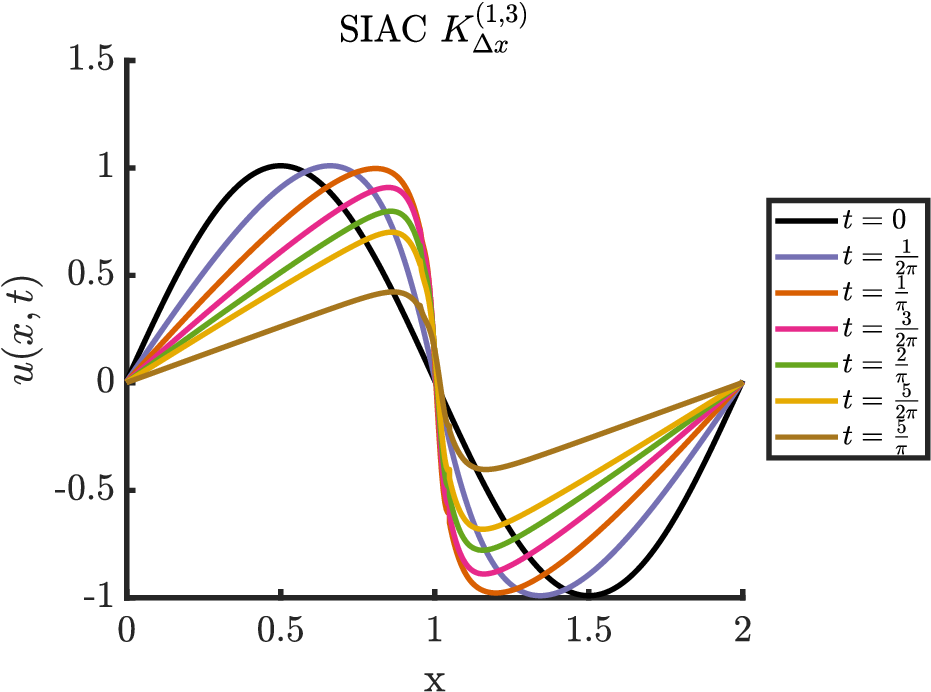}     &  
      \includegraphics[width=0.33\linewidth]{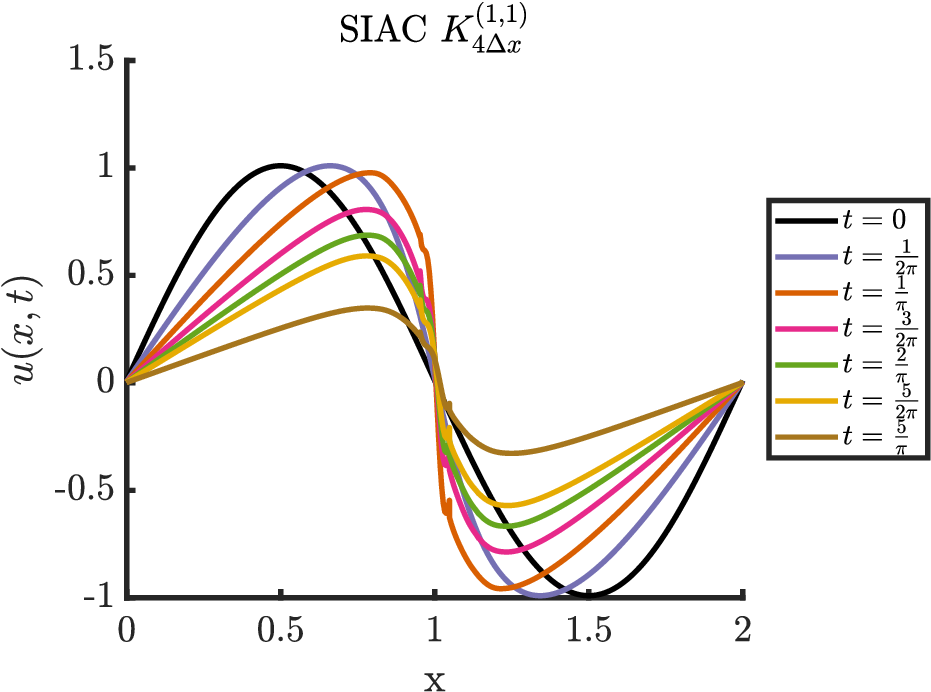} 
        \end{tabular}

    \caption{Comparison of SIAC correction regularization results for varying kernel parameterizations: number of moments/splines (Left), spline order/smoothness (Center), and kernel scaling (Right). Here $K^{(1,1)}_{\Delta x}$ is the baselines kernel with changed parameters as indicated. For all tests $c_e=10$ and $c_{max}=1$ was used. Increasing moments degrades solution quality when the approximation loses smoothness. Increasing the order and scaling makes the regularization more dissipative.}
    \label{fig:1D_param_sol}
\end{figure}

%% file: figures/files/1D_param_mass_eng_comp.tex
\begin{figure}
    \centering
    \begin{tabular}{c c c}
    \multicolumn{3}{c}{Energy and Relative Conservation Error over Time}\\
       Moments & Spline Order & Kernel Scaling \\
    \includegraphics[width=0.33\linewidth]{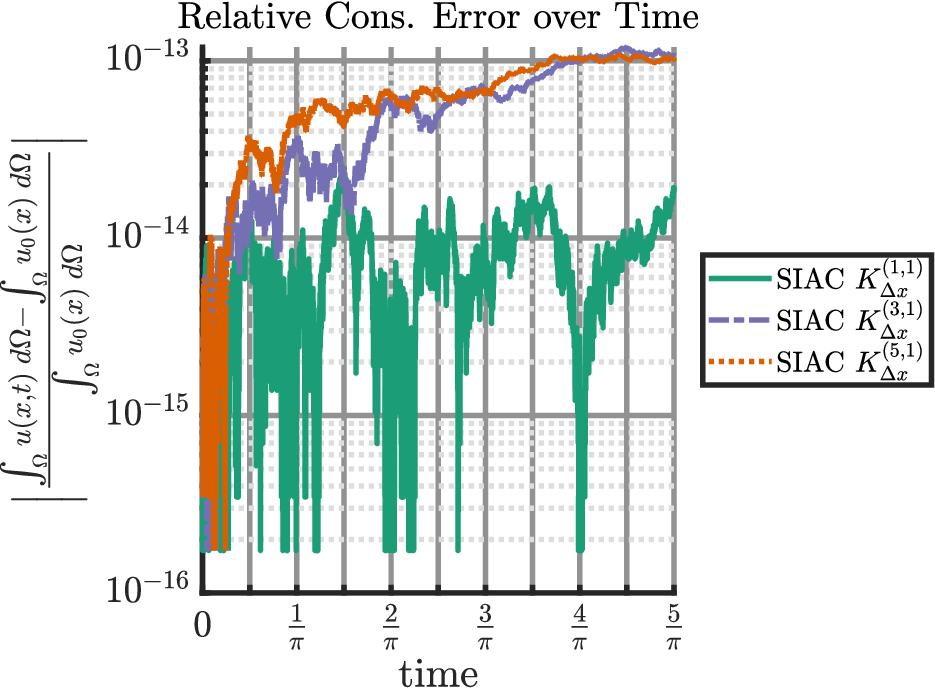}    &   \includegraphics[width=0.33\linewidth]{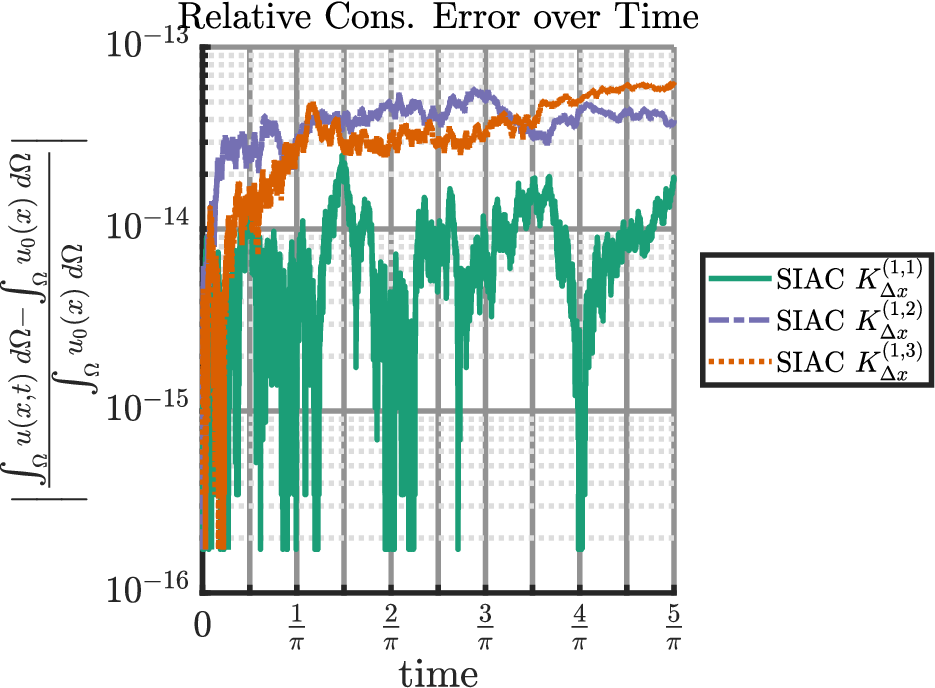}       &   \includegraphics[width=0.33\linewidth]{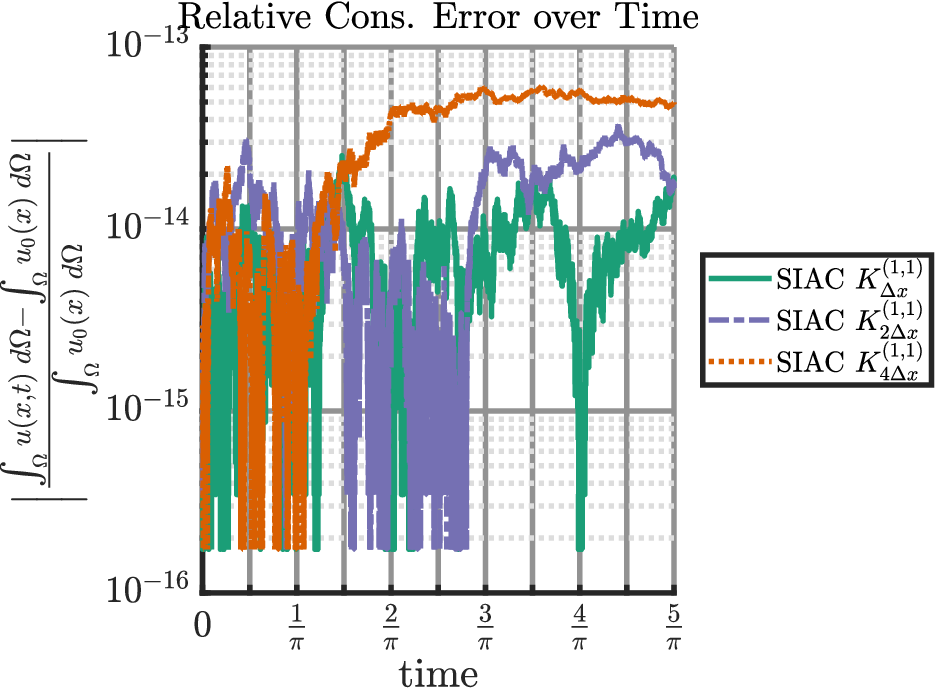}    \\
       \includegraphics[width=0.33\linewidth]{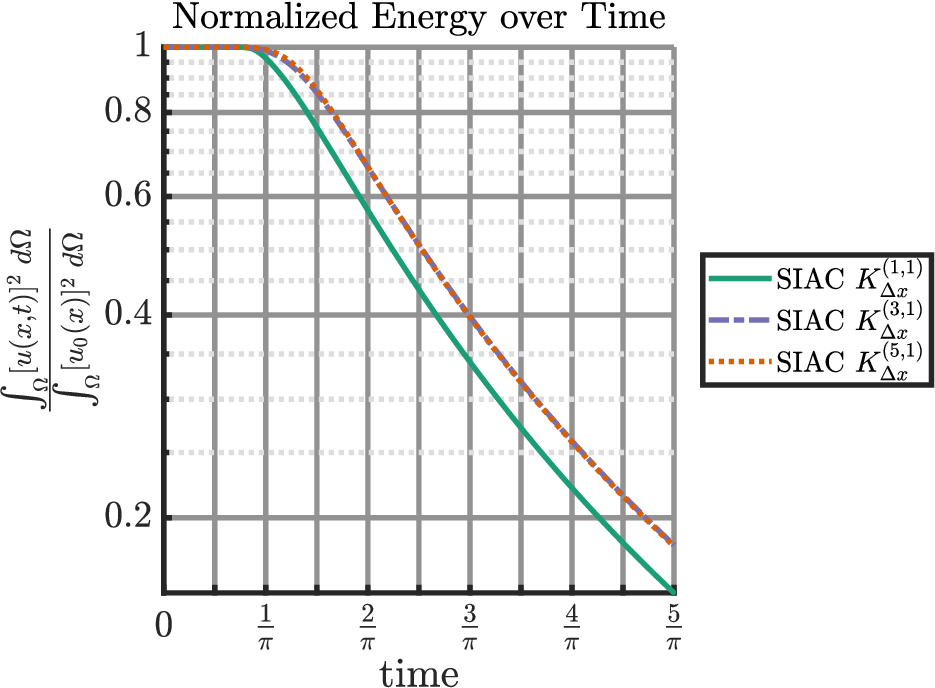}    &   \includegraphics[width=0.33\linewidth]{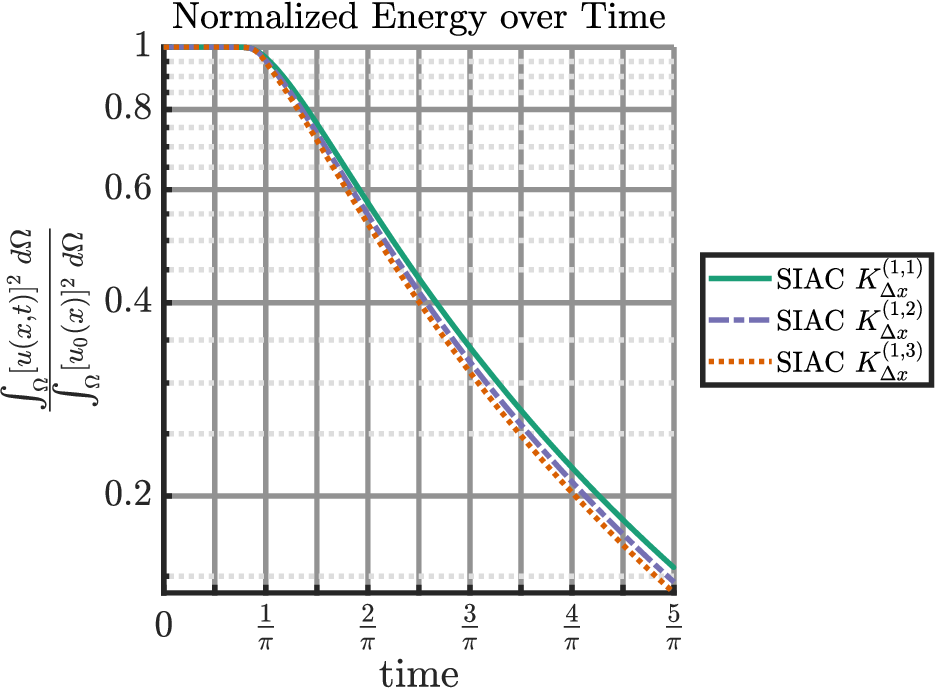}       &   \includegraphics[width=0.33\linewidth]{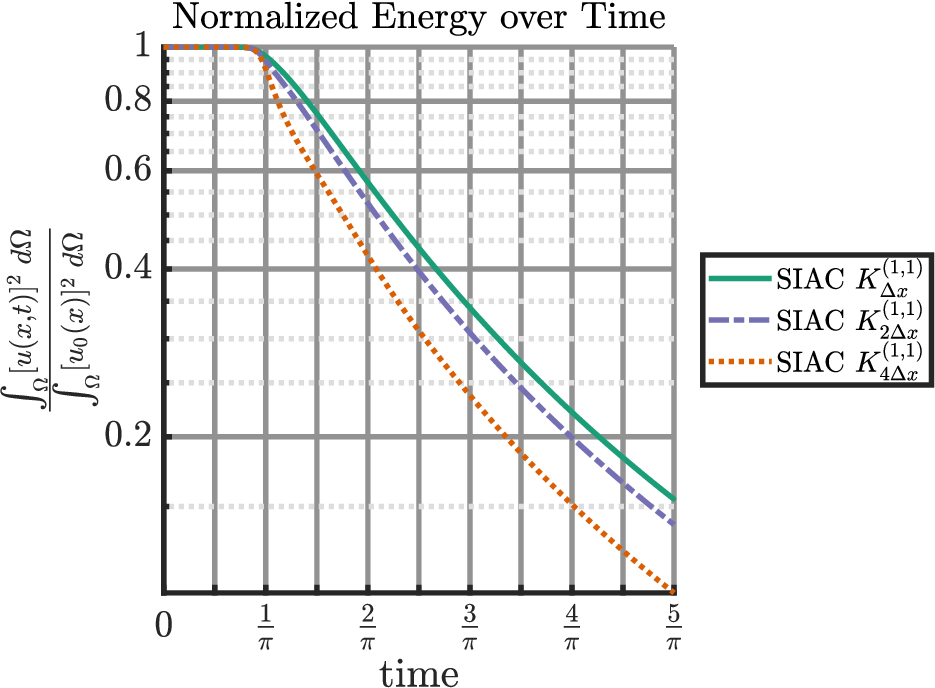}    \\
        \end{tabular}

    \caption{Comparison of the relative conservation error (Row 1) and energy dissipation (Row 2) of the SIAC correction regularization for varying kernel parameterizations. Here $K^{(1,1)}_{\Delta x}$ is the baselines kernel with changed parameters as indicated. For all tests $c_e=10$ and $c_{max}=1$ was used.}
    \label{fig:1D_param_energy_mass}
\end{figure}

%% file: figures/files/FV_ref.tex
\begin{figure}
    \centering
       \includegraphics[width=0.33\linewidth]{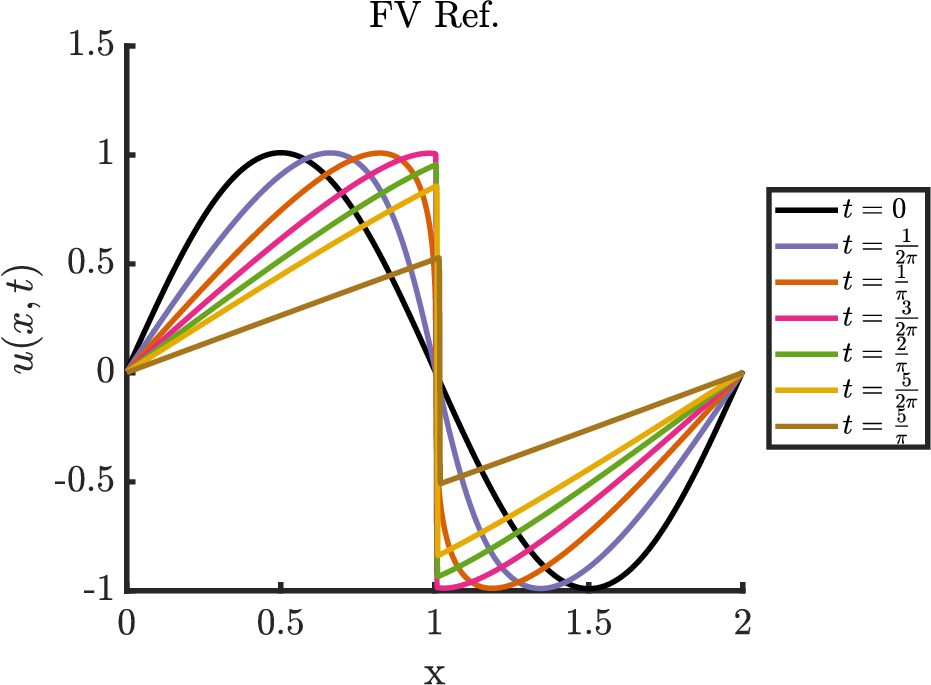}    

    \caption{Reference finite-volume solution to Burgers equation with $N_{elm}=2000$ elements and a Local Lax-Friedrichs flux.}
    \label{fig:FV reference}
\end{figure}

%% file: figures/files/1D_local_vs_SIAC_local_params.tex
\begin{figure}
    \centering
    \begin{tabular}{c c c}
    \multicolumn{3}{c}{Local tuning $c_E=4$, $c_{max}=0.15$ }\\
  
\includegraphics[width=0.25\linewidth]{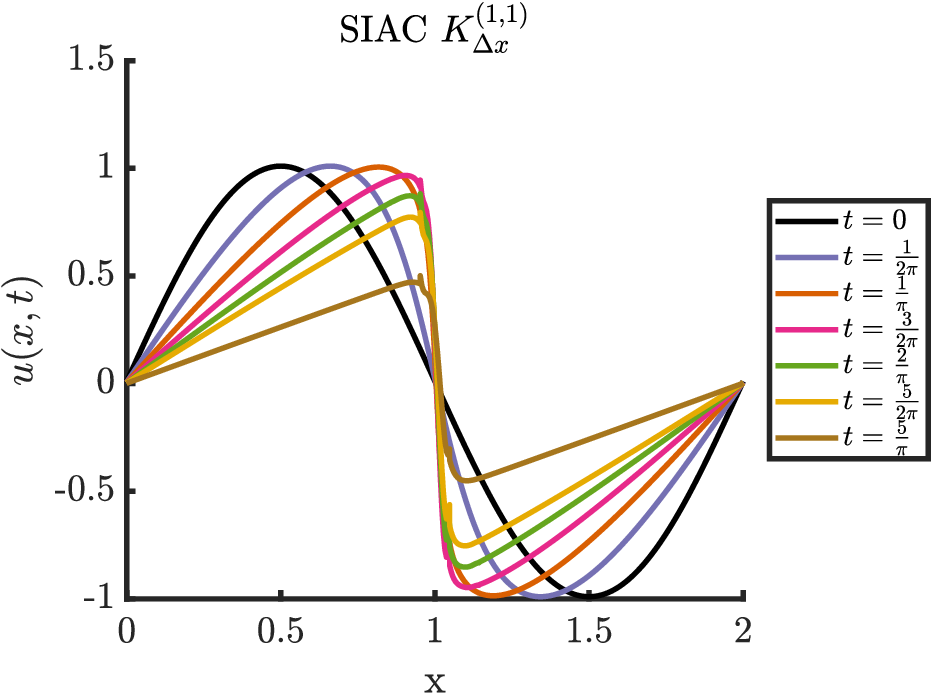}

 &
    \includegraphics[width=0.25\linewidth]{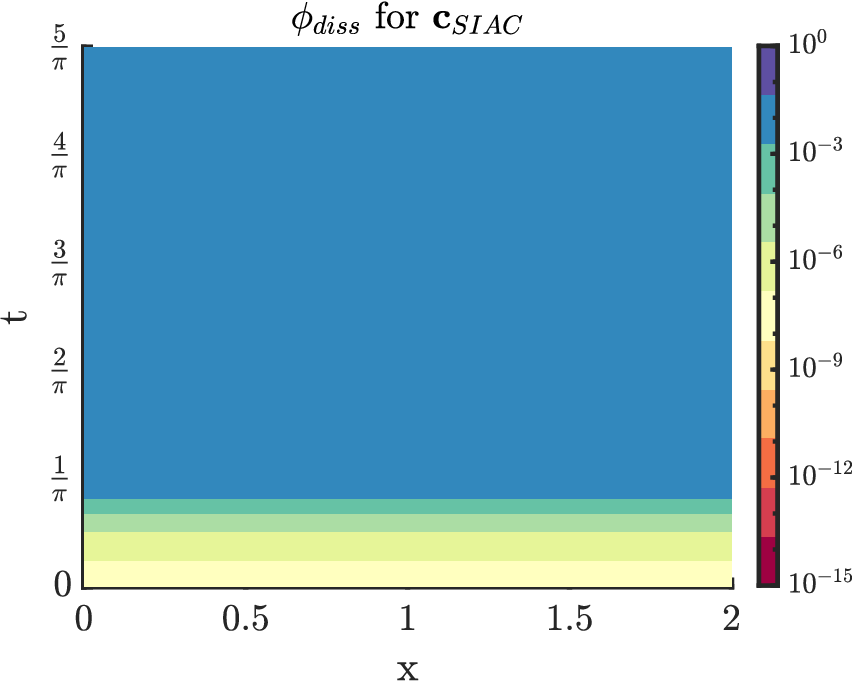}  
       &
   \includegraphics[width=0.25\linewidth]{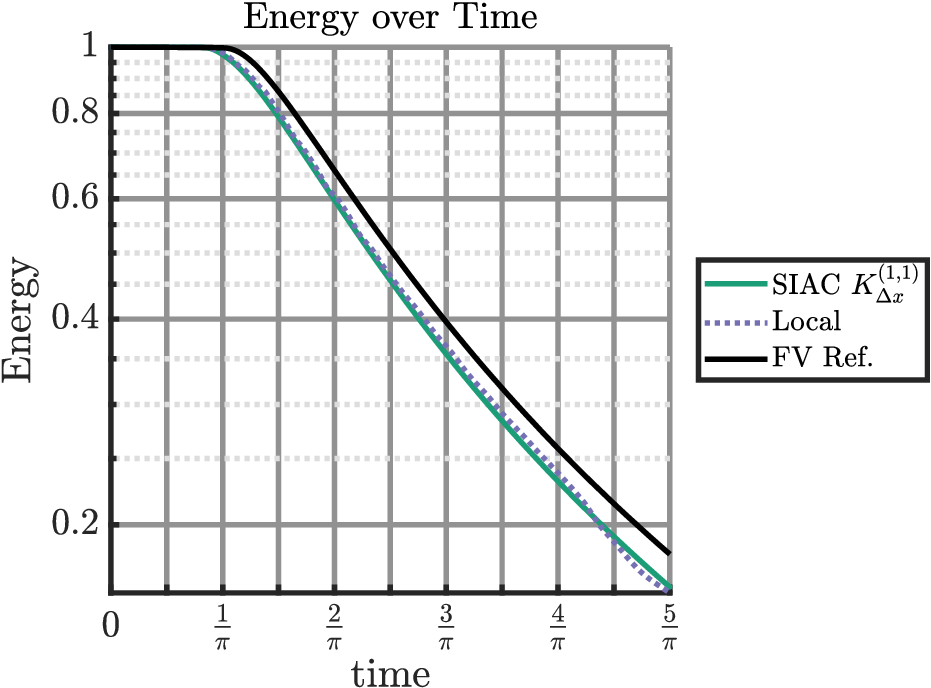} 
   \\
\includegraphics[width=0.25\linewidth]{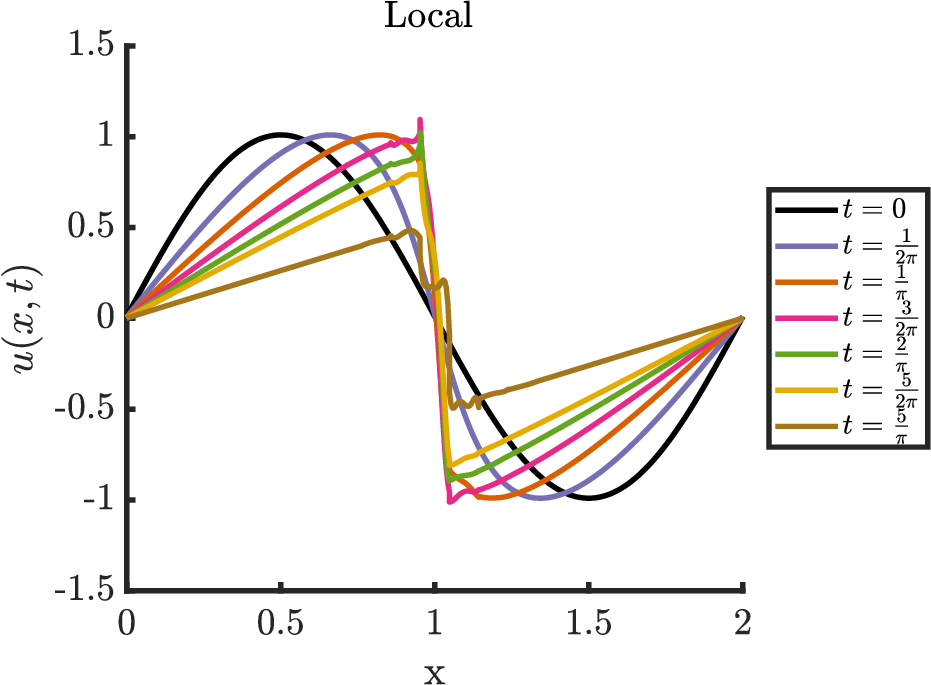}

   &
    \includegraphics[width=0.25\linewidth]{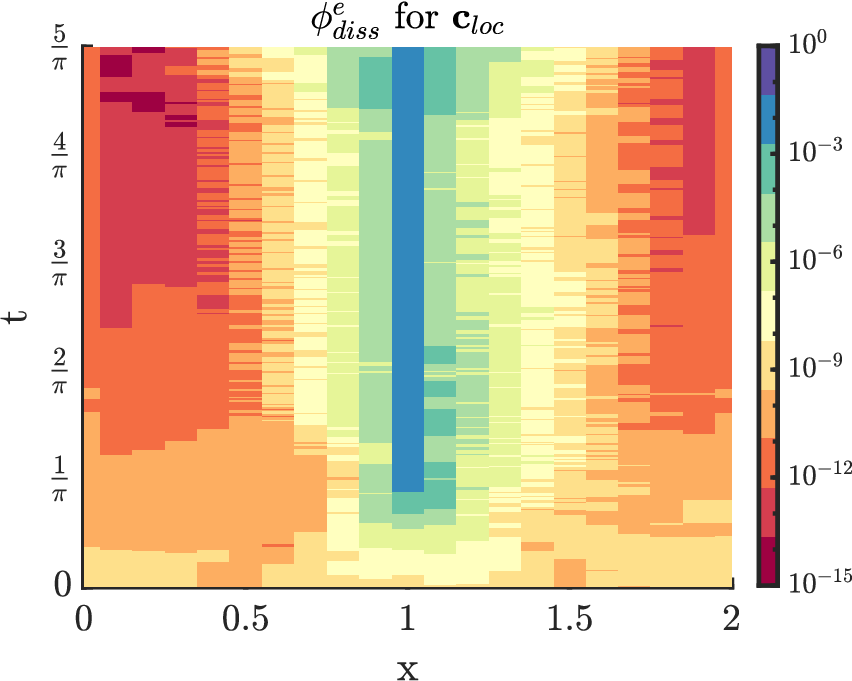} 
       &
   \includegraphics[width=0.25\linewidth]{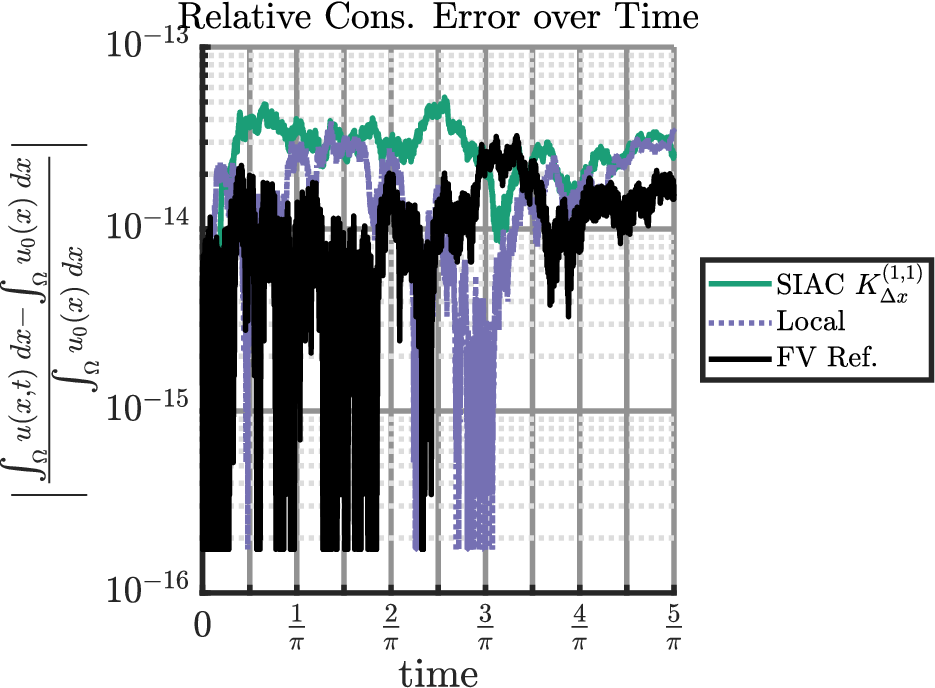} 
      
        \end{tabular}

    \caption{Comparison of local and SIAC correction approaches to regularization. Here the tuning parameters used were those that qualitatively gave the best performance for the local correction.}
    \label{fig:1D_local_vs_SIAC_local_params}

\end{figure}

%% file: figures/files/1D_local_vs_SIAC_SIAC_params.tex
\begin{figure}
    \centering
    \begin{tabular}{c c c}
    \multicolumn{3}{c}{SIAC tuning $c_E=10$, $c_{max}=1$ }\\
  
\includegraphics[width=0.25\linewidth]{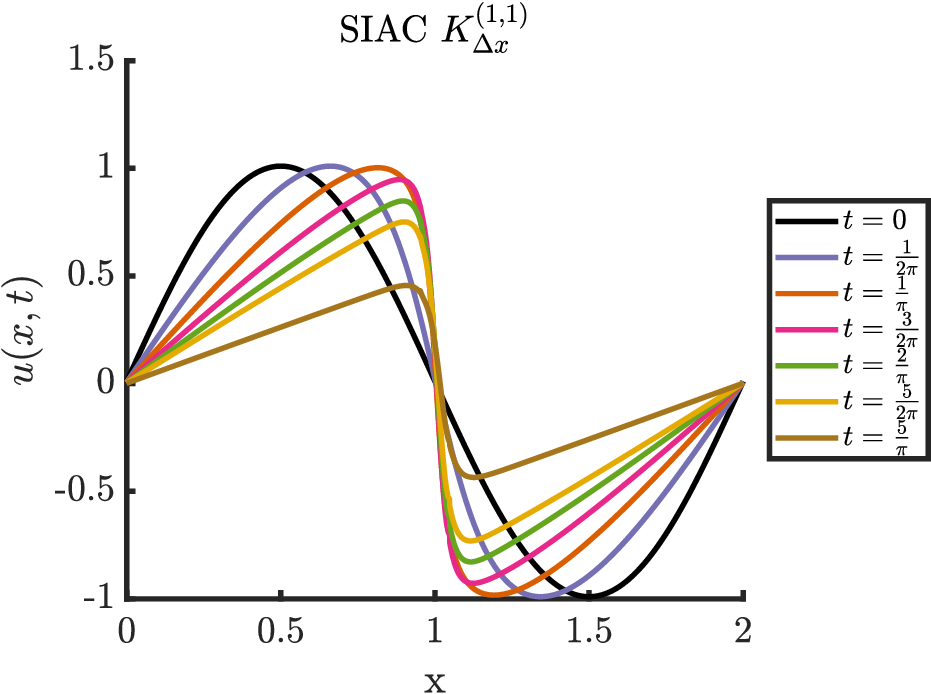}
 &
    \includegraphics[width=0.25\linewidth]{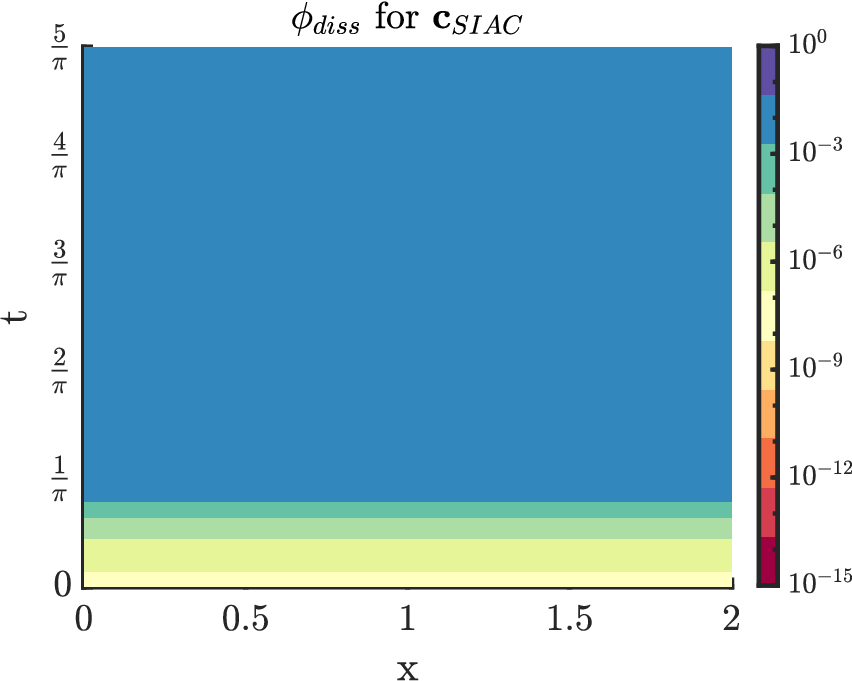}  
       &
   \includegraphics[width=0.25\linewidth]{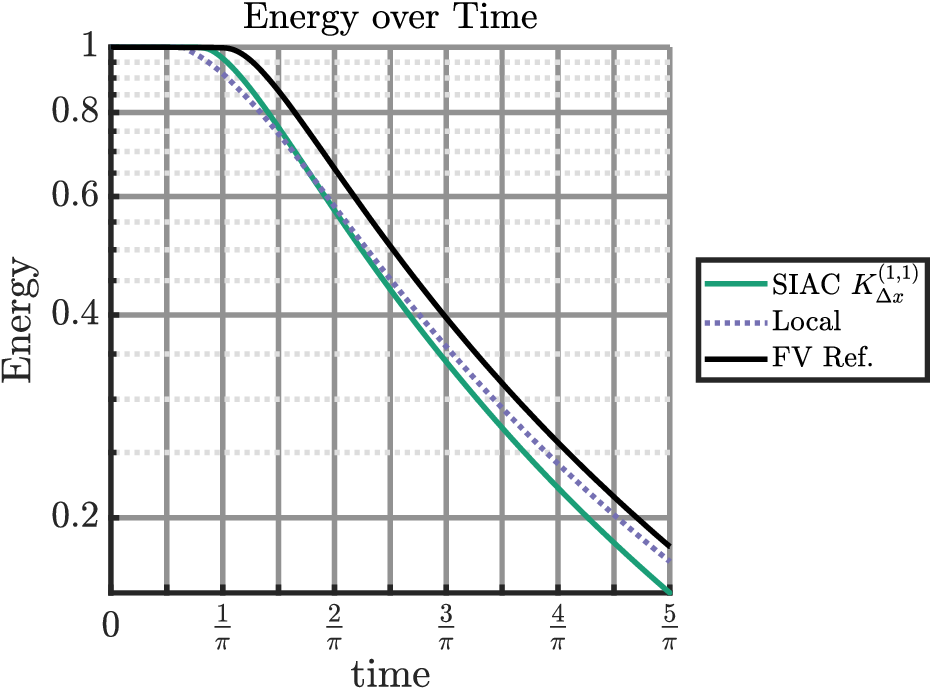} 
   \\
\includegraphics[width=0.25\linewidth]{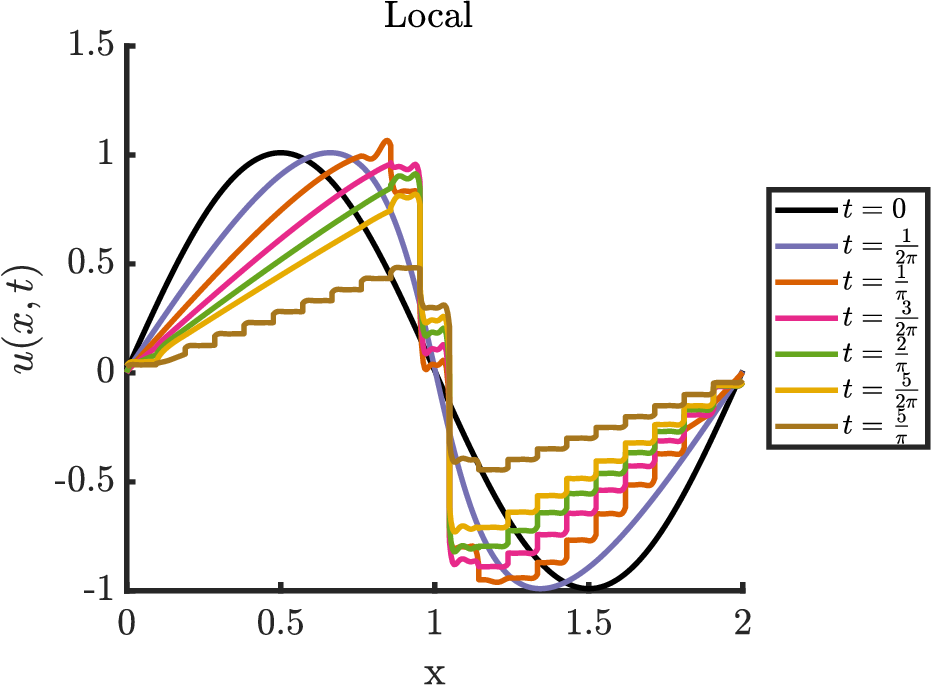}
       &
         
    \includegraphics[width=0.25\linewidth]{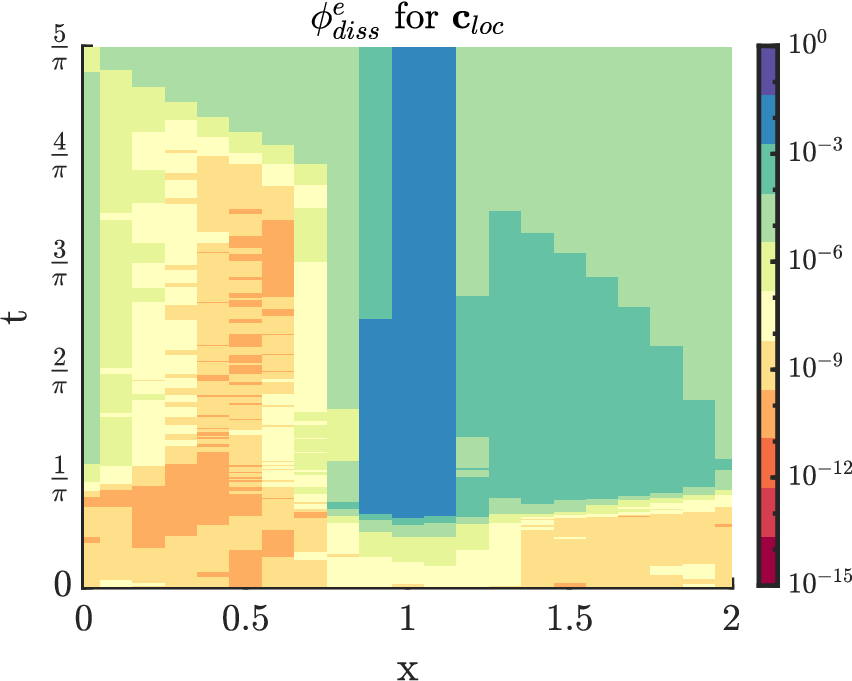}   
   &\includegraphics[width=0.25\linewidth]{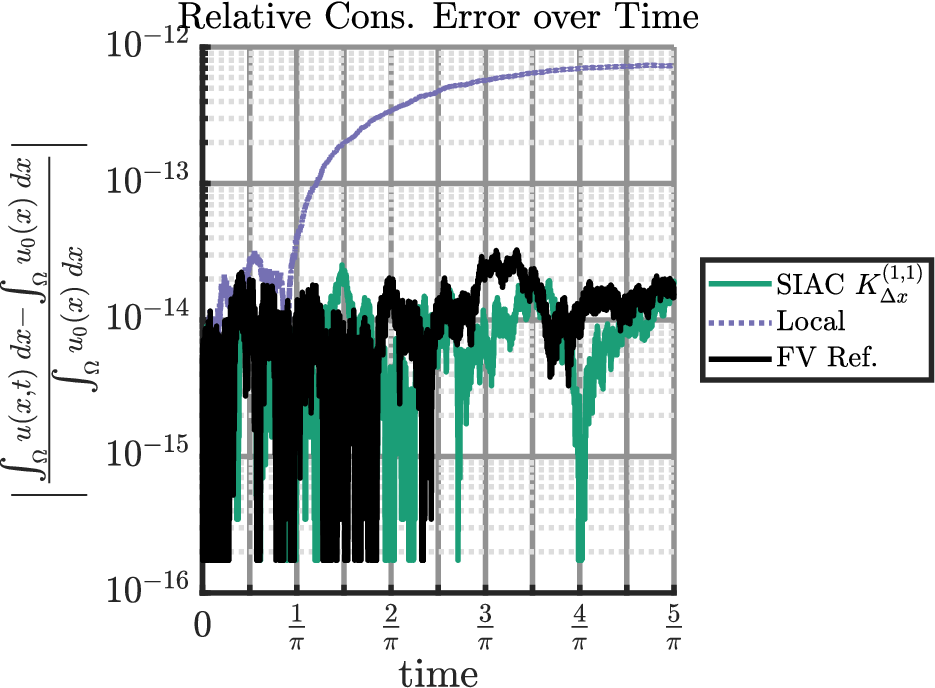} 
     
        \end{tabular}

    \caption{Comparison of local and SIAC correction approaches to regularization. Here the tuning parameters used were those that qualitatively gave the best performance for the SIAC correction.}
    \label{fig:1D_local_vs_SIAC_SIAC_params}
\end{figure}

%% file: figures/files/2D_reg_tests.tex
\begin{figure}
    \centering
    \begin{tabular}{c c c}
    \multicolumn{3}{c}{Two-Dimensional Regularization Example}\\
  
\includegraphics[width=0.33\linewidth]{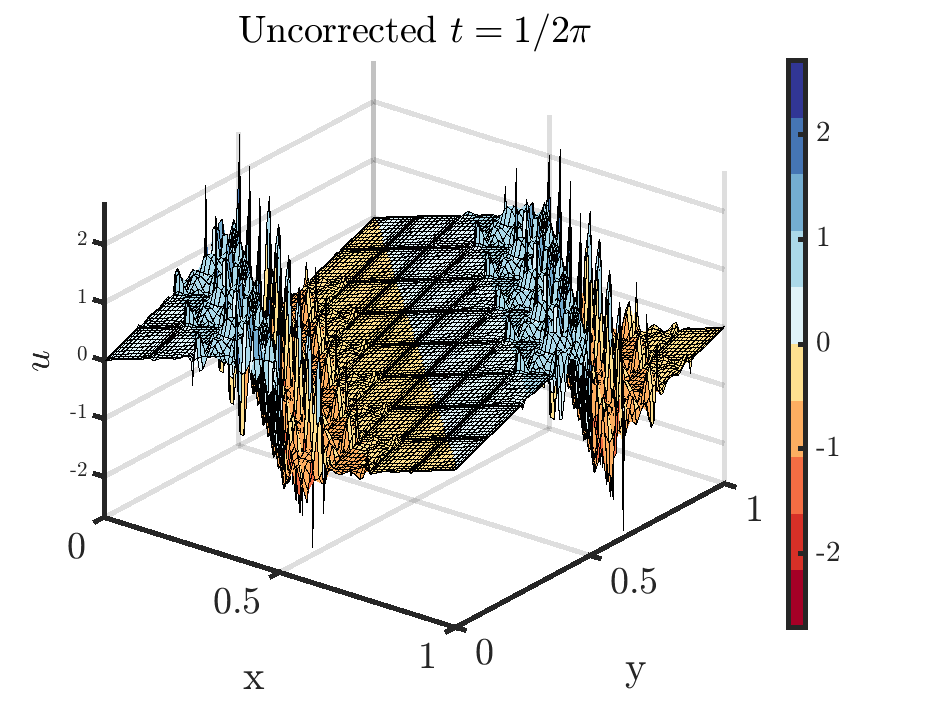}

 &
\includegraphics[width=0.33\linewidth]{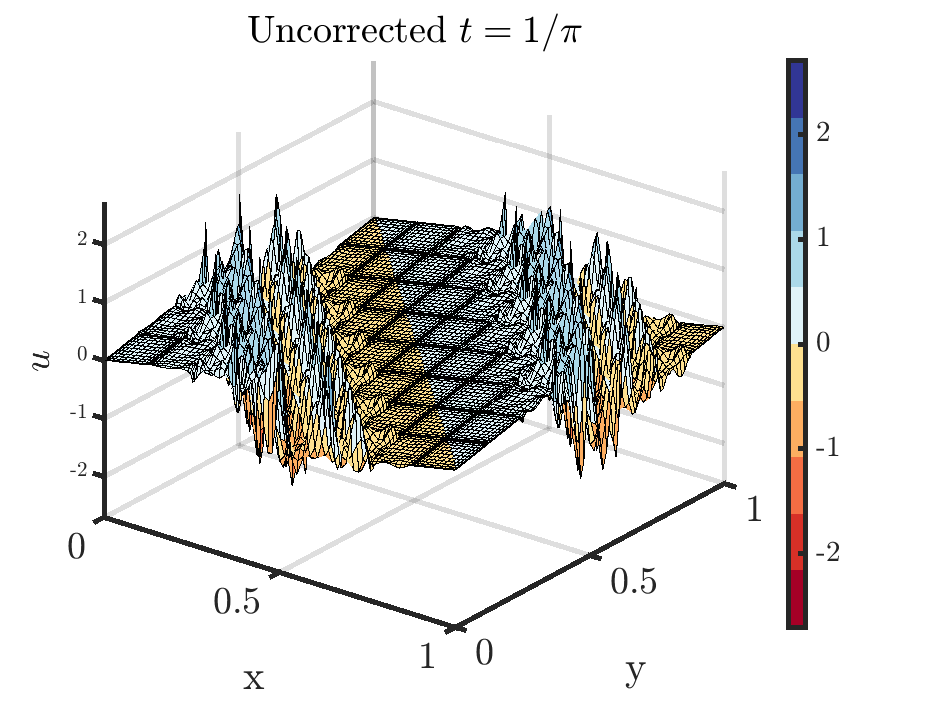} 
       &
\includegraphics[width=0.33\linewidth]{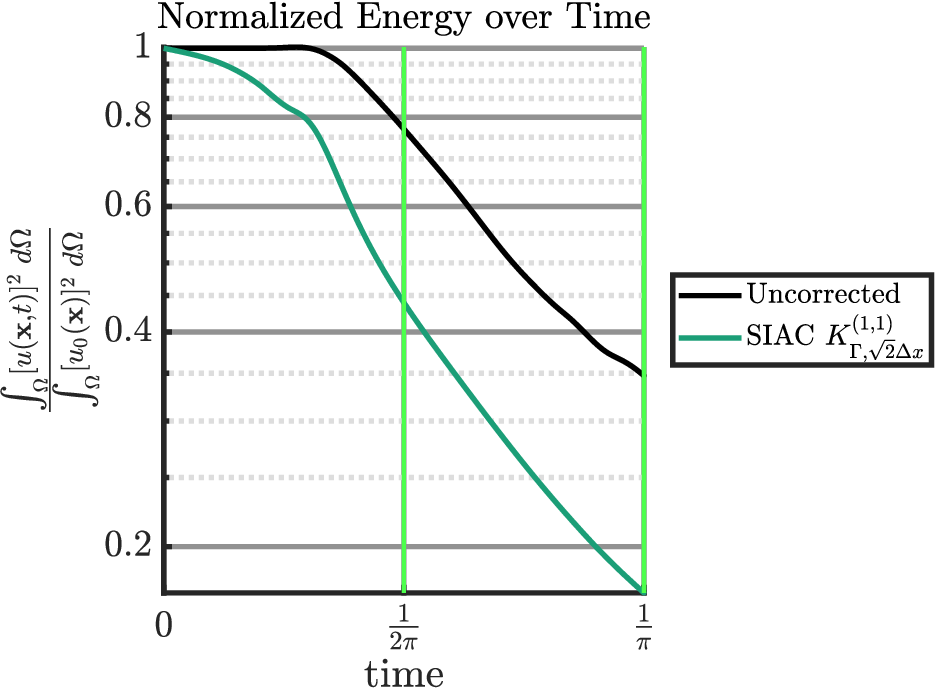}
   \\
\includegraphics[width=0.33\linewidth]{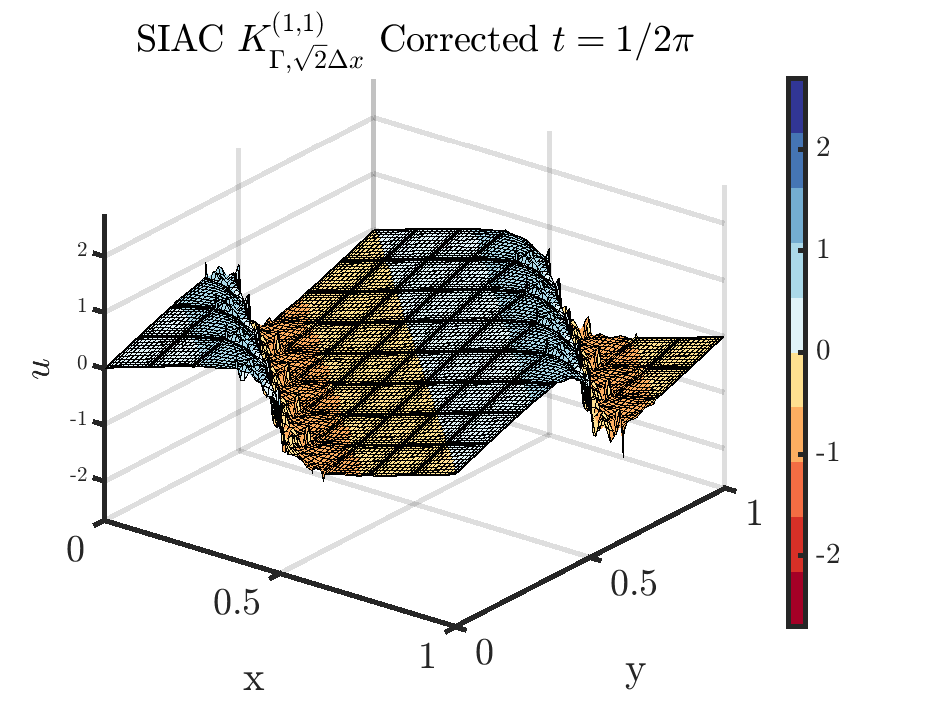}

 &
\includegraphics[width=0.33\linewidth]{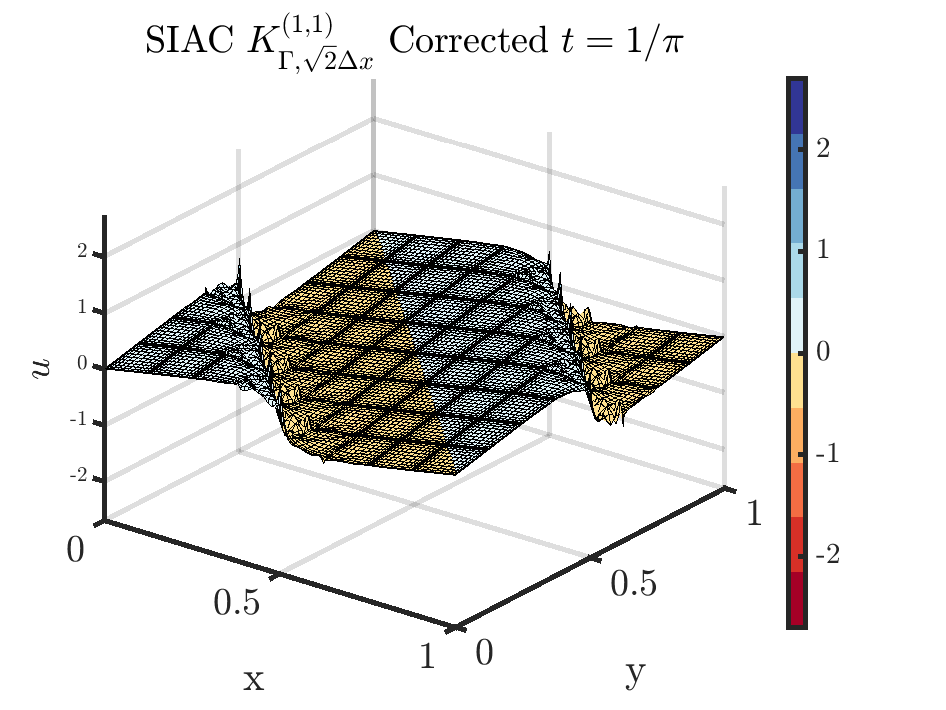} 
       &
\includegraphics[width=0.33\linewidth]{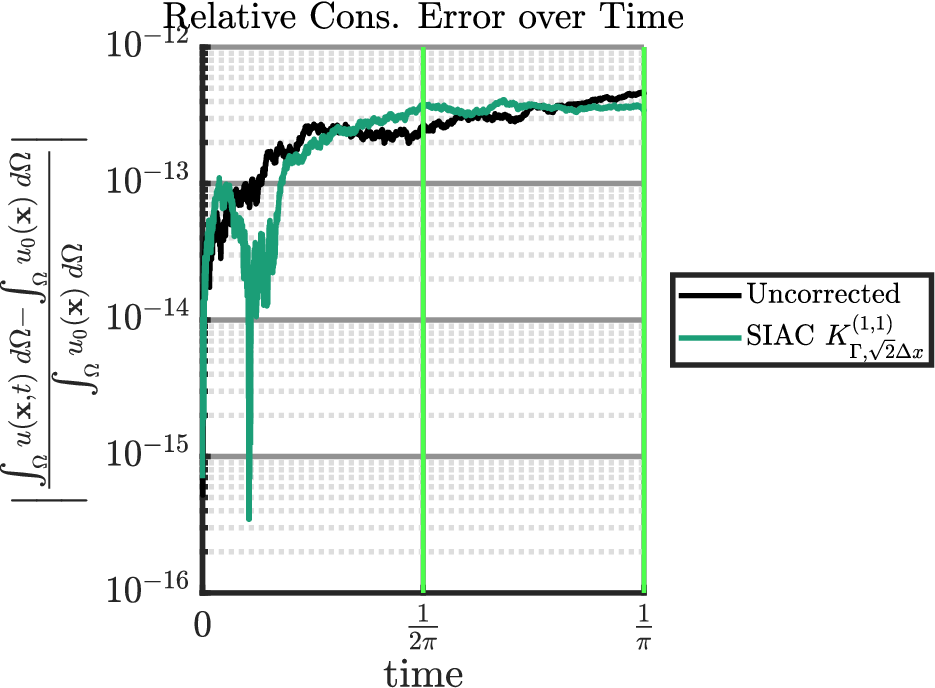}
        \end{tabular}

    \caption{Comparison of uncorrected and corrected 2D Modal DG approximations of the periodic inviscid Burgers equation with $u_0(x,y)=\sin(2\pi(x+y))+0.01$ on $[0,1]^2$ with $p=5$ and $N_{elm}=8\times8=64$ elements. For regularization we used $c_E=0.1$ and $c_{max}=0.15$.}
    \label{fig:2D_Reg_Example}
\end{figure}

%% file: figures/files/2D_reg_cross_section.tex
\begin{figure}
    \centering
    \begin{tabular}{c c}
    \multicolumn{2}{c}{Two-Dimensional Regularization $y=x$ Cross-Sections}\\
  
\includegraphics[width=0.33\linewidth]{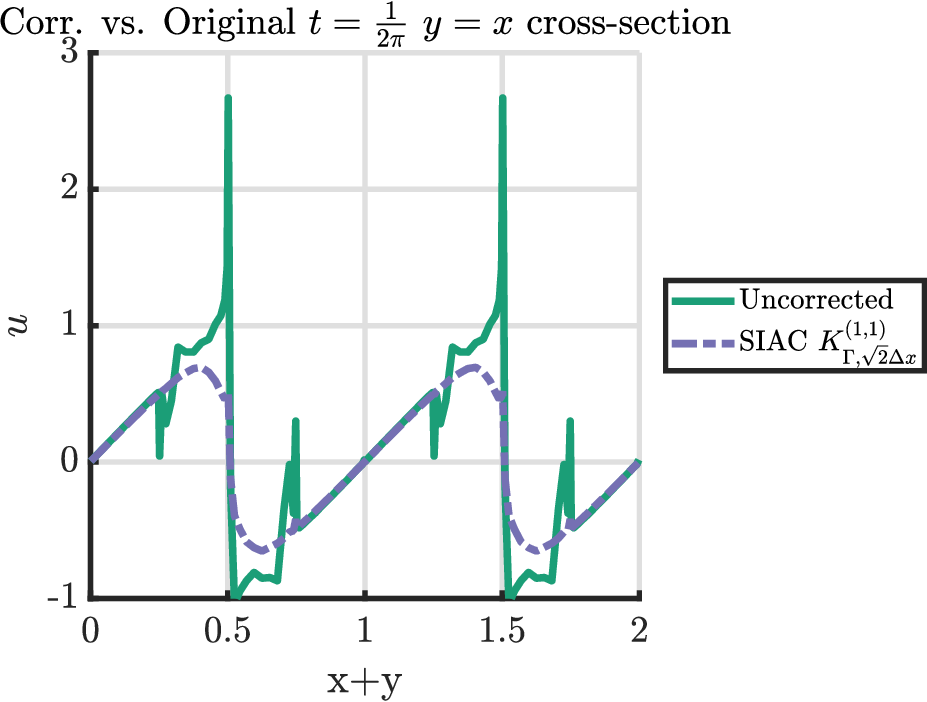}

 &
\includegraphics[width=0.33\linewidth]{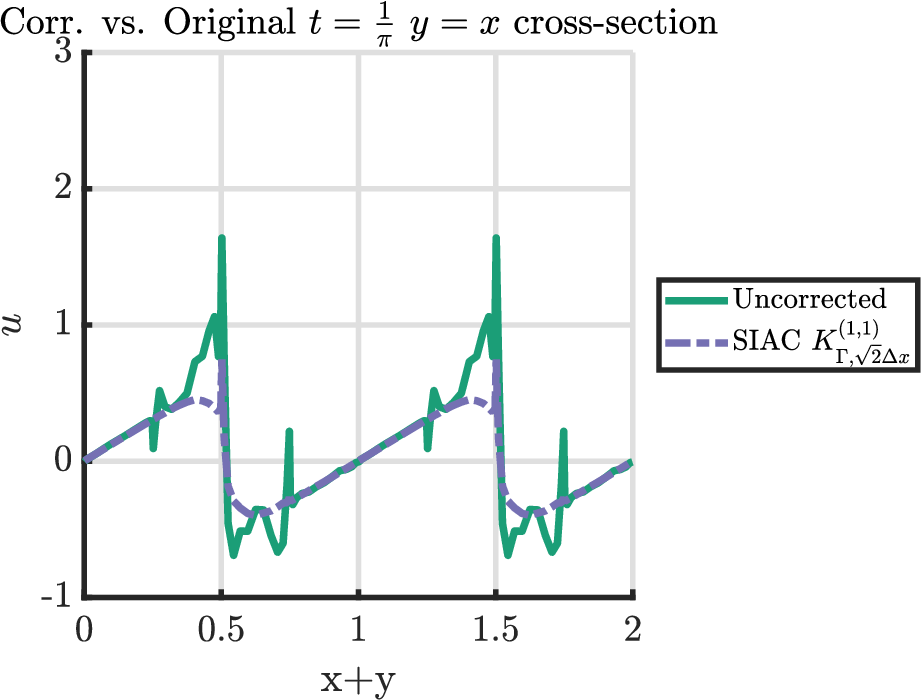} 
        \end{tabular}

    \caption{Cross-sections $(y=x)$ of uncorrected and corrected 2D Modal DG approximations of the periodic inviscid Burgers equation with $u_0(x,y)=\sin(2\pi(x+y))+0.01$ on $[0,1]^2$ with $p=5$ and $N_{elm}=8\times8=64$ elements. Solutions depicted at shock formation $t=\frac{1}{2\pi}$ and afterwards at $t=\frac{1}{\pi}$.}
    \label{fig:2D_Reg_cross}
\end{figure}

%% file: figures/files/Tabs_1D_Spatial_error.tex
\begin{table}
    \centering
    \begin{tabular}{|c|c|c|c|c|}\hline 
    \multicolumn{5}{|c|}{Uncorrected $L^2-$Error Order}\\ \hline
    $N_0\rightarrow2N_0$  & $p=1$ & $p=2$ &$p=3$& $p=4$\\ \hline
   $N_0=20$ &     1.3664 &   2.4936   & 2.9906  &  4.5180 \\ \hline
   $N_0=40$ &   1.4108  &  2.5765  &  3.2344  &  4.4176  \\ \hline
  $N_0=80$ &    1.4557  &  2.6864 &   3.2991  &  4.5704  \\ \hline
   $N_0=160$ &   1.4749  &  2.7806  &  3.3551 &   4.7000 \\ \hline

    \multicolumn{5}{|c|}{Local Correction $L^2-$Error Order}\\ \hline
    $N_0\rightarrow2N_0$   & $p=1$ & $p=2$ &$p=3$& $p=4$\\ \hline
   $N_0=20$ &     0.8521&   2.4927& 2.9913&  4.5180\\ \hline
   $N_0=40$ &       0.4958&   2.5745&  3.2343&  4.4176\\ \hline
  $N_0=80$ &        0.2446&  2.6851&   3.2990&  4.5704\\ \hline
   $N_0=160$ &       0.1321&  2.7798&  3.3551&   4.6996\\ \hline
    \multicolumn{5}{|c|}{$K^{(3,2)}_{\Delta x}$ $L^2-$Error Order}\\ \hline
    $N_0\rightarrow2N_0$ & $p=1$ & $p=2$ &$p=3$& $p=4$\\ \hline
   $N_0=20$ & 1.3289&    2.4992 &   2.9907 &   4.5180\\ \hline
 $N_0=40$ &   1.3901 &   2.5847 &   3.2344 &   4.4176\\ \hline
  $N_0=80$ &  1.4481 &   2.6938 &   3.2991  &  4.5704\\ \hline
  $N_0=160$ &  1.4743 &   2.7861 &   3.3551   & 4.7000\\ \hline
     \multicolumn{5}{|c|}{$K^{(1,1)}_{\Delta x}$ $L^2-$Error Order}\\ \hline  
    $N_0\rightarrow2N_0$   & $p=1$ & $p=2$ &$p=3$& $p=4$\\ \hline
   $N_0=20$ &1.3389 &   2.4948 &   2.9906 &   4.5180\\ \hline
   $N_0=40$ & 1.3941 &   2.5768&    3.2344  &  4.4176\\ \hline
  $N_0=80$ &  1.4485 &   2.6864 &   3.2991 &   4.5704\\ \hline
   $N_0=160$ & 1.4729&    2.7806 &   3.3551  &  4.7000\\ \hline
    \end{tabular}
    \caption{$L^2$ orders of convergence for DGSEM solution to Burgers equation with initial data $\sin(\pi x)+0.01$ on $[0,2]$ run to time $t=1/2\pi$ and correction type as indicated. }
    \label{tab:spatial_orders}
\end{table}

%% file: figures/files/1D_Spatial_Conv_tests.tex
\begin{figure}
    \centering
    \begin{tabular}{c c }
       \includegraphics[width=0.33\linewidth]{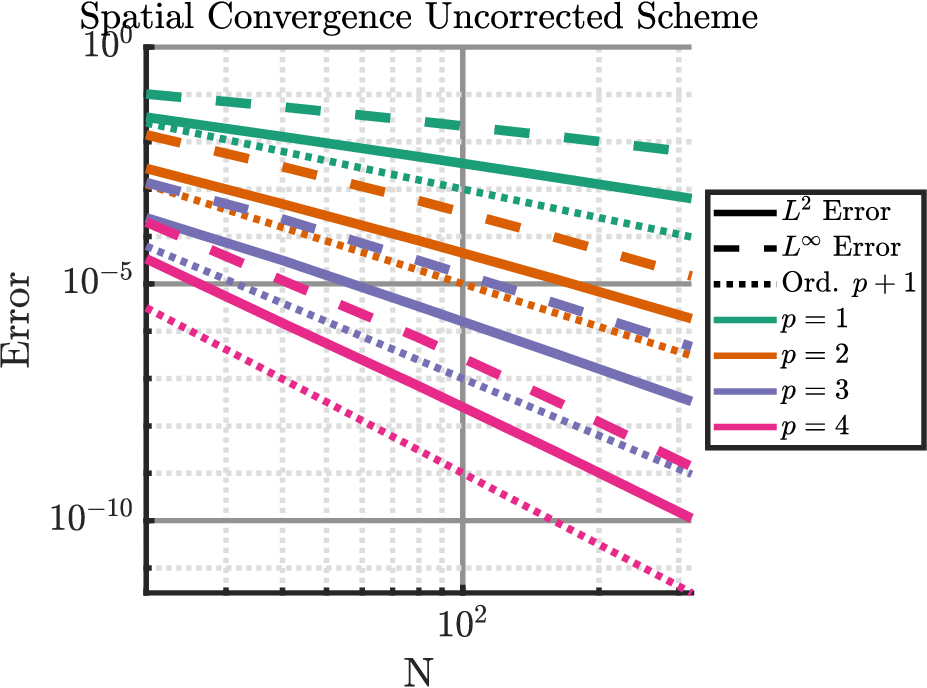}    &    \includegraphics[width=0.33\linewidth]{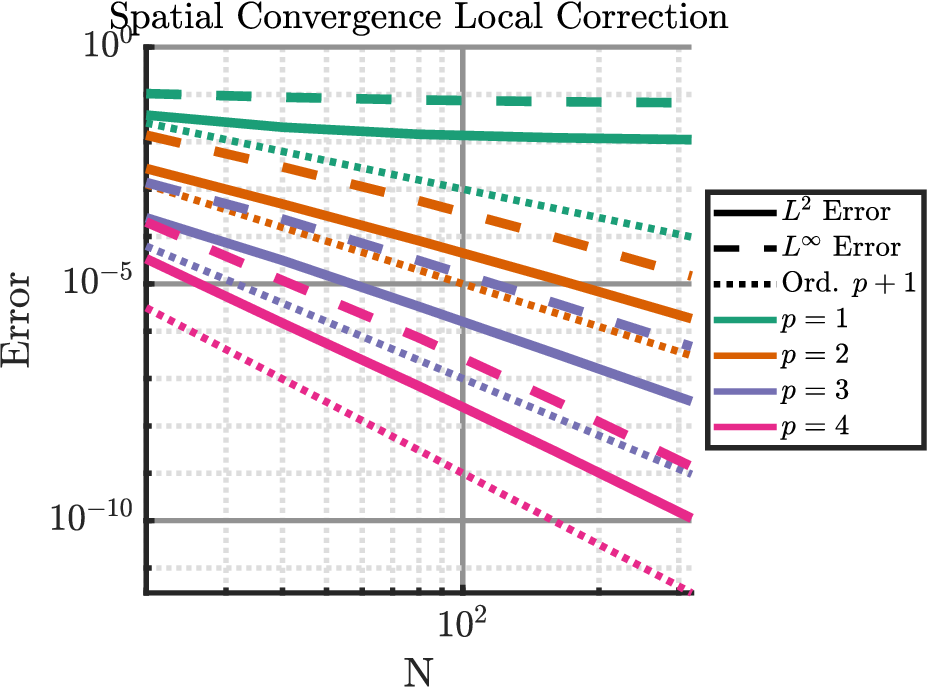} \\
           \includegraphics[width=0.33\linewidth]{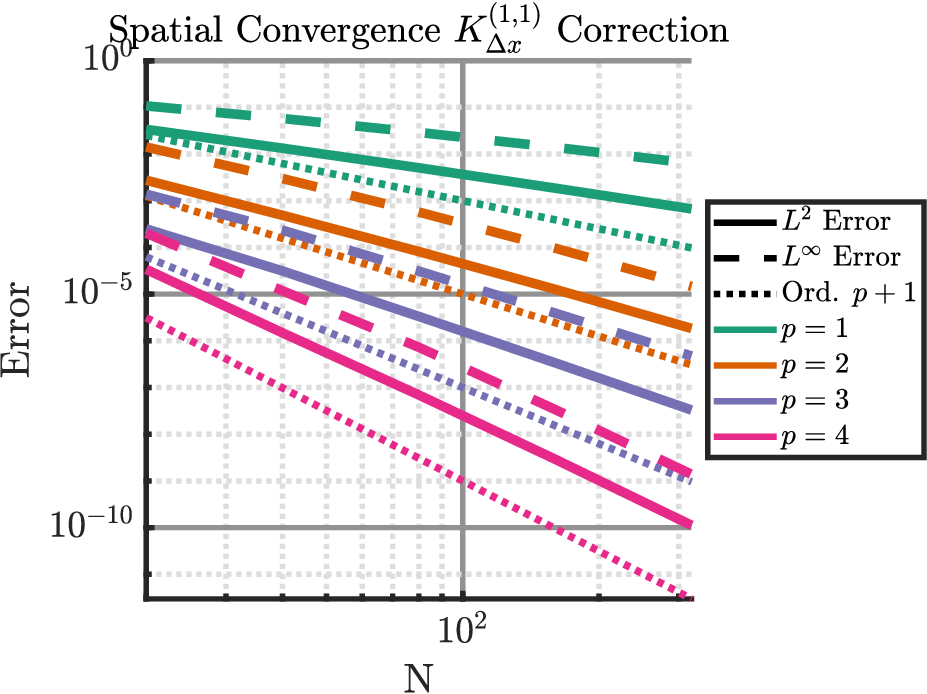} & \includegraphics[width=0.33\linewidth]{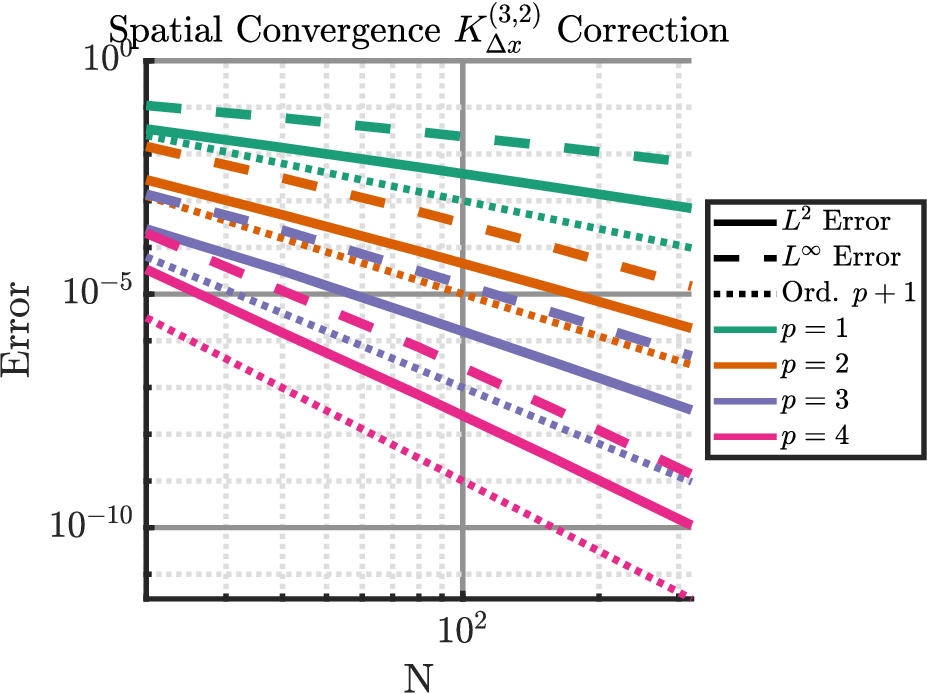}
        \end{tabular}

    \caption{Spatial order of convergence results for the periodic one-dimensional Burgers problem with $u(x,0)=\sin(\pi x )$ on $\Omega=[0,2]$ prior to shock formation. Exact solution obtained by Newton iteration.}
    \label{fig:spatial_conv_1D}
\end{figure}

%% file: figures/files/2D_Spatial_Conv_tests.tex
\begin{figure}
    \centering
      \begin{tabular}{c c }
       \includegraphics[width=0.33\linewidth]{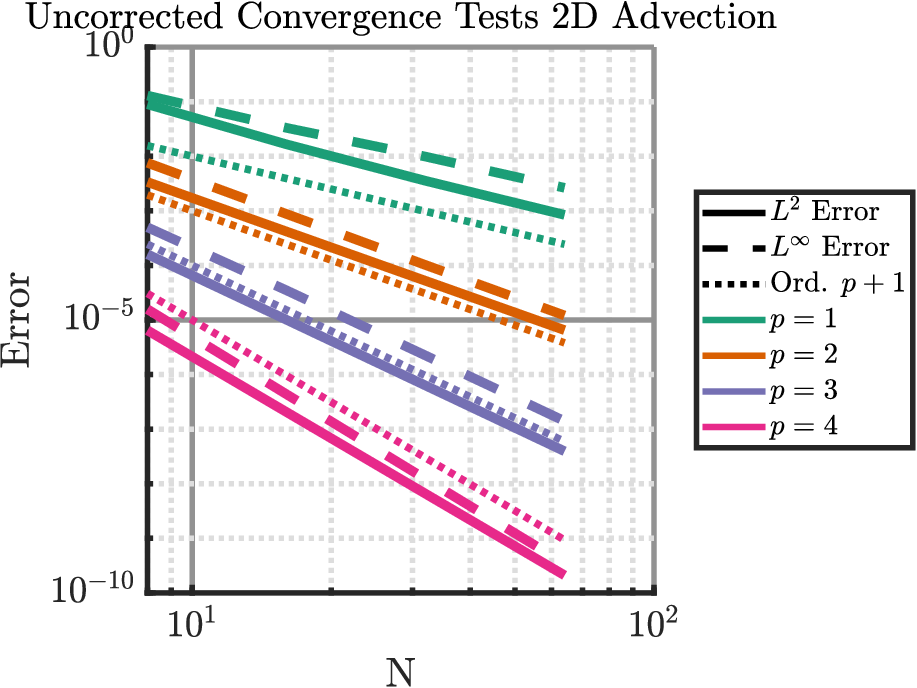}    &    \includegraphics[width=0.33\linewidth]{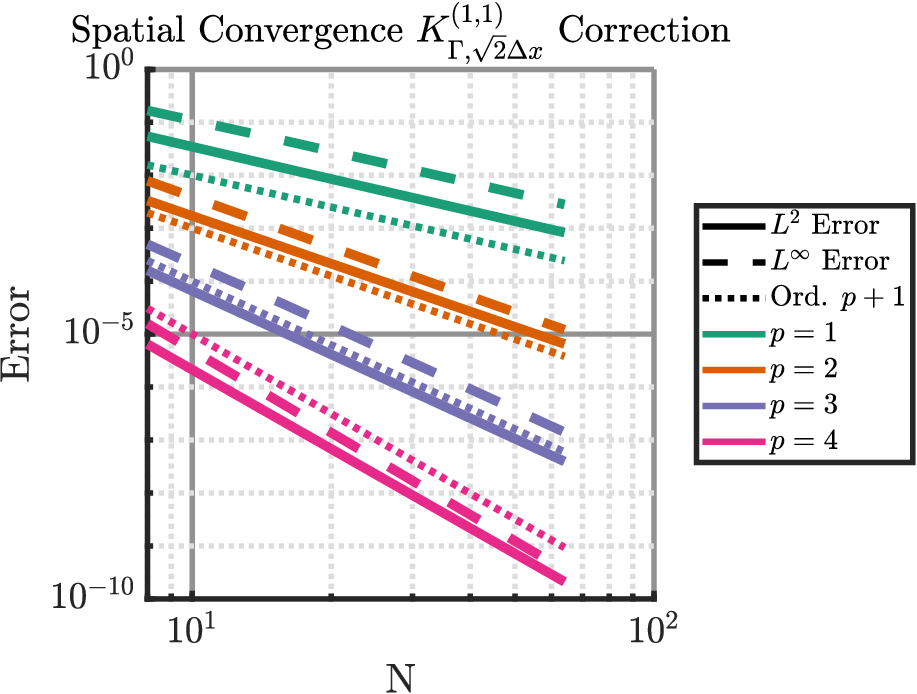}
        \end{tabular}

    \caption{Spatial order of convergence results for the periodic two-dimensional constant advection problem. Here $N$ is the number of elements in a given Cartesian direction.}
    \label{fig:spatial_conv_2D}
\end{figure}

%% file: figures/files/1D_Temporal_Conv.tex
\begin{figure}
    \centering
    \begin{tabular}{c c c}
       \includegraphics[width=0.3\linewidth]{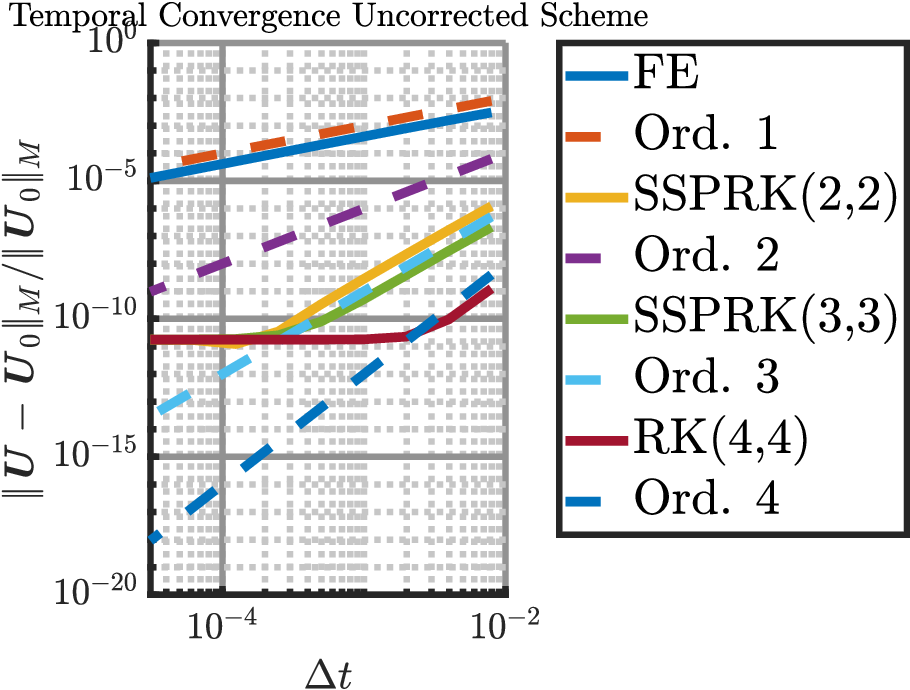}    &    \includegraphics[width=0.3\linewidth]{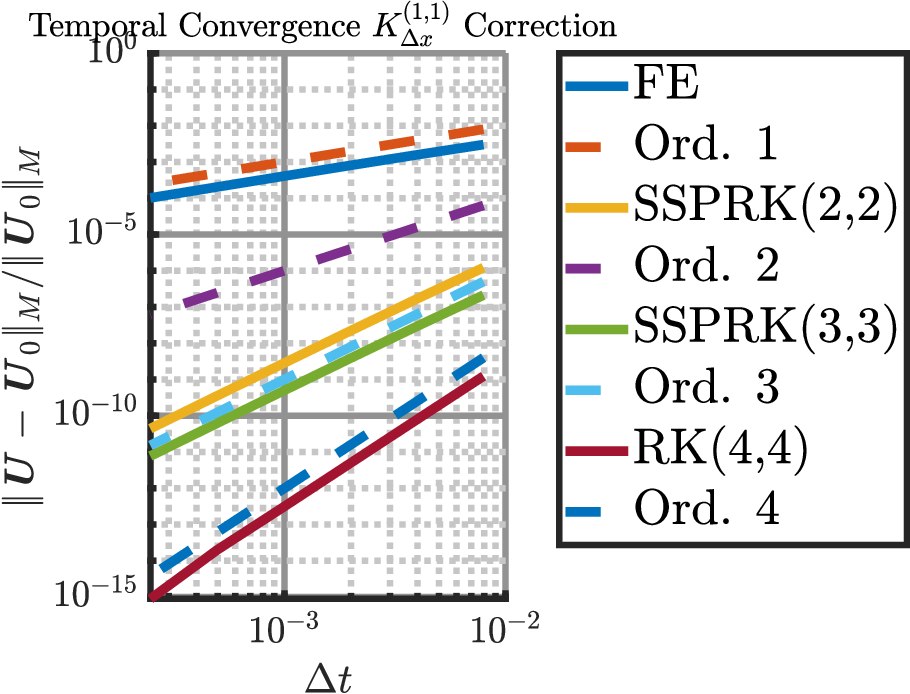}  &    \includegraphics[width=0.3\linewidth]{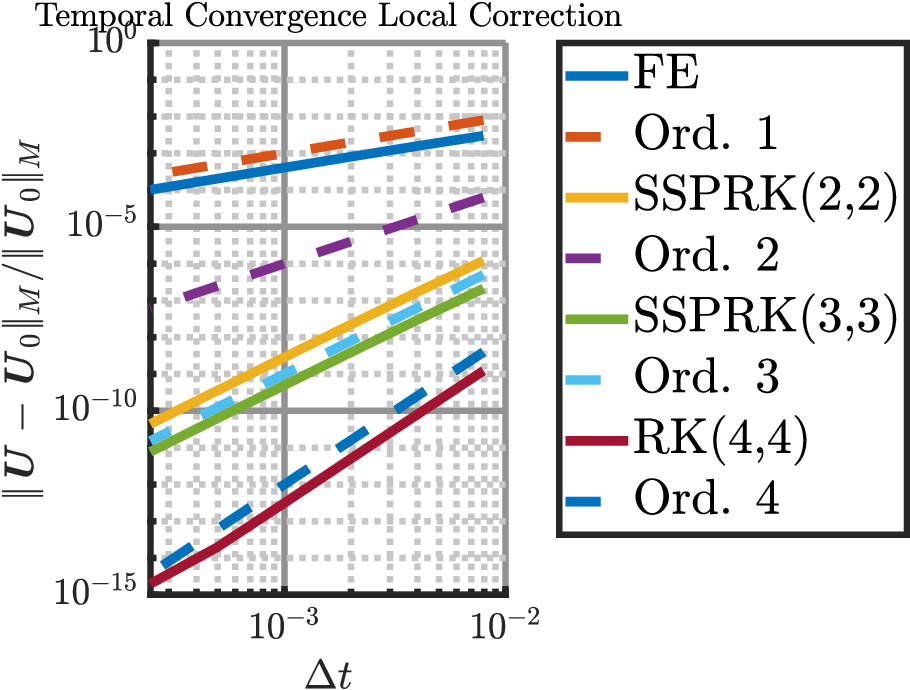} 
        \end{tabular}

    \caption{Relative energy change over time using non-energy conserving time-stepping methods.}
    \label{fig:temporal_conv}
\end{figure}

%% file: sections/conclusions.tex
\section{Conclusions}\label{sec:conc}

We have demonstrated in this work that the local entropy correction approach of Abgrall can be viewed from the perspective of a filtering procedure, and generalized to include arbitrary conservative filters. The significant flexibility in the choice of filter enables one to apply multi-element (L)SIAC filters for satisfying auxiliary dynamics. The present study has focused on energy conservation of one- and two-dimensional inviscid Burgers equation, where energy conservative timestepping was guaranteed by application of relaxation Runge-Kutta methods. The SIAC-based correction was seen to have a smaller magnitude correction relative to the scheme residual compared to the original single-element correction approach, and was shown to have similar performance in terms of a subcell energy conservation metric. Furthermore, upon adding a gradient based regularization term for shock capture, the SIAC approach was shown to produce qualitatively better solutions than the local correction implementation. Error of accuracy experiments demonstrate that the SIAC-based correction preserves the original order of accuracy of the scheme, even in the $p=1$ case while the original approach did not. In the future, we expect that SIAC-style correction procedure will prove useful in more complex scenarios such as under-resolved compressible flows, enforcement of multiple target dynamics, and alternative entropies, as well for stabilization of meshless methods such as radial basis function schemes.

%% file: sections/acknowledgments.tex
\section{Acknowledgments}\label{sec:ackno}

The authors would like to acknowledge Professor Jennifer Ryan for discussions regarding the method development as well as the current manuscript.

\subsection*{Financial disclosure}
The first author received no funding for the preparation of this manuscript. The second author recognizes funding and support for this work by the Air Force Office of Scientific Research (AFOSR) (program officers Drs. Chiping Li and Fariba Farhoo), as well as Jacobs Engineering Inc. under contract No. FA9300-20-F-9801.

\subsection*{Conflict of interest}

The authors declare no potential conflict of interests. Any opinions, findings, and conclusions or recommendations expressed in this material are those of the author(s) and do no necessarily reflects views of the United States Air Force.